\documentclass[11pt]{article}
\usepackage{amsfonts}
\usepackage{amsfonts,mathrsfs,amssymb,amsthm,mathptm}
\usepackage{amsmath,amscd}
\usepackage{mathptm,pslatex}
\usepackage{color}
\usepackage[dvips]{graphicx}
\usepackage[all]{xy}
\usepackage{graphicx}
\usepackage{fancyhdr}

\begin{document}
\newtheorem{Def}{Definition}[section]
\newtheorem{Bsp}[Def]{Example}
\newtheorem{Prop}[Def]{Proposition}
\newtheorem{Theo}[Def]{Theorem}
\newtheorem{Lem}[Def]{Lemma}
\newtheorem{Koro}[Def]{Corollary}
\theoremstyle{definition}
\newtheorem{Rem}[Def]{Remark}

\newcommand{\add}{{\rm add}}
\newcommand{\con}{{\rm con}}
\newcommand{\gd}{{\rm gl.dim}}
\newcommand{\dm}{{\rm dom.dim}}
\newcommand{\tdim}{{\rm dim}}
\newcommand{\E}{{\rm E}}
\newcommand{\Mor}{{\rm Morph}}
\newcommand{\End}{{\rm End}}
\newcommand{\ind}{{\rm ind}}
\newcommand{\rsd}{{\rm res.dim}}
\newcommand{\rd} {{\rm rep.dim}}
\newcommand{\ol}{\overline}
\newcommand{\overpr}{$\hfill\square$}
\newcommand{\rad}{{\rm rad}}
\newcommand{\soc}{{\rm soc}}
\renewcommand{\top}{{\rm top}}
\newcommand{\pd}{{\rm proj.dim}}
\newcommand{\id}{{\rm inj.dim}}
\newcommand{\fld}{{\rm flat.dim}}
\newcommand{\Fac}{{\rm Fac}}
\newcommand{\Gen}{{\rm Gen}}
\newcommand{\fd} {{\rm fin.dim}}
\newcommand{\Fd} {{\rm Fin.dim}}
\newcommand{\Pf}[1]{{\mathscr P}^{<\infty}(#1)}
\newcommand{\DTr}{{\rm DTr}}
\newcommand{\cpx}[1]{#1^{\bullet}}
\newcommand{\D}[1]{{\mathscr D}(#1)}
\newcommand{\Dz}[1]{{\mathscr D}^+(#1)}
\newcommand{\Df}[1]{{\mathscr D}^-(#1)}
\newcommand{\Db}[1]{{\mathscr D}^b(#1)}
\newcommand{\C}[1]{{\mathscr C}(#1)}
\newcommand{\Cz}[1]{{\mathscr C}^+(#1)}
\newcommand{\Cf}[1]{{\mathscr C}^-(#1)}
\newcommand{\Cb}[1]{{\mathscr C}^b(#1)}
\newcommand{\Dc}[1]{{\mathscr D}^c(#1)}
\newcommand{\K}[1]{{\mathscr K}(#1)}
\newcommand{\Kz}[1]{{\mathscr K}^+(#1)}
\newcommand{\Kf}[1]{{\mathscr  K}^-(#1)}
\newcommand{\Kb}[1]{{\mathscr K}^b(#1)}
\newcommand{\modcat}{\ensuremath{\mbox{{\rm -mod}}}}
\newcommand{\Modcat}{\ensuremath{\mbox{{\rm -Mod}}}}

\newcommand{\stmodcat}[1]{#1\mbox{{\rm -{\underline{mod}}}}}
\newcommand{\pmodcat}[1]{#1\mbox{{\rm -proj}}}
\newcommand{\imodcat}[1]{#1\mbox{{\rm -inj}}}
\newcommand{\Pmodcat}[1]{#1\mbox{{\rm -Proj}}}
\newcommand{\Imodcat}[1]{#1\mbox{{\rm -Inj}}}
\newcommand{\opp}{^{\rm op}}
\newcommand{\otimesL}{\otimes^{\rm\mathbb L}}
\newcommand{\rHom}{{\rm\mathbb R}{\rm Hom}\,}
\newcommand{\projdim}{\pd}
\newcommand{\Hom}{{\rm Hom}}
\newcommand{\Coker}{{\rm Coker}}
\newcommand{ \Ker  }{{\rm Ker}}
\newcommand{ \Cone }{{\rm Con}}
\newcommand{ \Img  }{{\rm Im}}
\newcommand{\Ext}{{\rm Ext}}
\newcommand{\StHom}{{\rm \underline{Hom}}}

\newcommand{\gm}{{\rm _{\Gamma_M}}}
\newcommand{\gmr}{{\rm _{\Gamma_M^R}}}

\def\vez{\varepsilon}\def\bz{\bigoplus}  \def\sz {\oplus}
\def\epa{\xrightarrow} \def\inja{\hookrightarrow}

\newcommand{\lra}{\longrightarrow}
\newcommand{\lraf}[1]{\stackrel{#1}{\lra}}
\newcommand{\ra}{\rightarrow}
\newcommand{\dk}{{\rm dim_{_{k}}}}

\newcommand{\colim}{{\rm colim\, }}
\newcommand{\limt}{{\rm lim\, }}
\newcommand{\Add}{{\rm Add }}
\newcommand{\Tor}{{\rm Tor}}
\newcommand{\Cogen}{{\rm Cogen}}
\newcommand{\Tria}{{\rm Tria}}
\newcommand{\tria}{{\rm tria}}

{\Large \bf
\begin{center}
Recollements of derived categories I: Exact contexts
\end{center}}

\medskip
\centerline{\textbf{Hongxing Chen} and \textbf{Changchang Xi}$^*$}

\renewcommand{\thefootnote}{\alph{footnote}}
\setcounter{footnote}{-1} \footnote{ $^*$ Corresponding author.
Email: xicc@cnu.edu.cn; Fax: 0086 10 68903637.}
\renewcommand{\thefootnote}{\alph{footnote}}
\setcounter{footnote}{-1} 
\renewcommand{\thefootnote}{\alph{footnote}}
\setcounter{footnote}{-1} 

\begin{abstract}
Recollements were introduced originally by Beilinson, Bernstein and Deligne to study the derived categories of
perverse sheaves, and nowadays become
very powerful in understanding relationship among three algebraic,
geometric or topological objects. The purpose of this series of papers is to study
recollements in terms of derived module categories and homological ring epimorphisms, and then to apply our results to both representation
theory and algebraic $K$-theory.

In this paper we present a new and systematic method to
construct recollements of derived module categories. For this aim,
we introduce a new ring structure, called the noncommutative tensor product, and give necessary and sufficient conditions for
noncommutative localizations which appears often in representation
theory, topology and $K$-theory, to be homological. The
input of our machinery is an exact context which can be easily
obtained from a rigid morphism that exists in very general
circumstances. The output is a recollement of derived module
categories of rings in which the
noncommutative tensor product of an exact context plays a crucial role. Thus we obtain a large variety of new recollements from commutative and noncommutative localizations, ring epimorphisms and extensions.
\end{abstract}

{\footnotesize \tableofcontents}

\section{Introduction\label{sect1}}

Recollements were first introduced by Beilinson, Berstein and
Deligne in 1982 in order to describe the derived categories of
perverse sheaves over singular spaces, by using derived versions of
Grothendieck's six functors (see \cite{Gr, BBD}). Later,
recollements of derived categories were employed to study
stratifications of the derived categories of modules over blocks of
the Bernstein-Gelfand-Gelfand category $\mathcal O$ (see
\cite{CPS}). Further, recollements were used by Happel to establish
a relationship among finitistic dimensions of finite-dimensional algebras (see
\cite{Happel}). Recently, they become of great interest in
understanding the derived categories of the endomorphism rings of
infinitely generated tilting modules (see \cite{Bz, AKL, xc1}). It
turns out that recollements are actually a very useful framework for
investigating relationships among three algebraic, geometric or
topological objects (see \cite{BBD, Ranicki}).

Recollements of derived module categories have an intimate connection
with homological ring epimorphisms (\cite{GL, Koen, xc1, NS}) which play a crucial role in many branches of mathematics.
Recall that a ring epimorphism $R\to S$ is said to be
\emph{homological} if $\Tor^R_i(S, S)=0$ for all $i>0$. In commutative algebra, homological ring epimorphisms often appear as
localizations which are one of the fundamental tools in algebraic geometry. In
representation theory, homological ring epimorphisms have been used
to study perpendicular categories, sheaves and stratifications of
derived module categories of rings (see \cite{GL, CPS, xc1}), and to construct infinitely generated tilting modules
(see \cite{HJ}). In algebraic $K$-theory, Neeman and Ranicki have
employed homological noncommutative localizations, a special class of
homological ring epimorphisms, to establish a useful long exact
sequence of algebraic $K$-groups (see \cite{nr1}), which generalizes
many earlier results in the literature (see \cite{ne2}). Also, in
Banach algebra, homological ring epimorphisms have been
topologically modified to investigate the analytic functional
calculus (see \cite{taylor}), where they were called
``localizations".

In this paper, we shall provide a systematic study of recollements of derived module categories through homological ring
epimorphisms, especially, those arising from noncommutative localizations
which have been used widely in topology and geometry (see, for
instance, \cite{nr1} and \cite{Ranicki}). For this purpose, we introduce the notion of exact contexts and define their noncommutative tensor products which not only generalise the usual tensor products
over commutative rings, but also cover some well-known constructions in the literature: coproducts of rings, dual extensions and endomorphism rings. Under a Tor-vanishing condition, we give a constructive method
to produce new homological noncommutative localizations and
recollements of derived categories of rings.  Roughly speaking, the
input of our machinery is a quadruple consisting of two ring homomorphisms, a
bimodule and a special element of the bimodule, such that they are linked by an
exact sequence. The output is a recollement of derived
module categories of rings in which the noncommutative tensor products play an essential role. As a consequence, we apply our general results to ring epimorphisms, (commutative and noncommutative) localizations and extensions, and get a large class of new recollements of derived module categories. This kind of recollements was already applied to study the Jordan-H\"older theorem for stratifications of derived module categories in \cite{xc1} and will be used to investigate relationships among homological or $K$-theoretical properties of three algebras
(see  \cite{xc5, xc6}).

Now, let us explain our results more explicitly. First of all, we introduce some
notation.

Let $R$, $S$ and $T$ be associative rings with identity, and let $\lambda:R\to S$
and $\mu:R\to T$ be ring homomorphisms.
Suppose that $M$ is an $S$-$T$-bimodule together with an element $m\in M$.
We say that the quadruple $(\lambda, \mu, M, m)$ is an \emph{exact context} if the following sequence
$$
0\lra R\lraf{(\lambda,\,\mu)}S\oplus T
\lraf{\left({\cdot\,m\,\atop{-m\,\cdot}}\right)}M\lra 0
$$
is an exact sequence of abelian groups, where $\cdot m$ and $m \cdot$ denote the right and
left multiplication by $m$ maps, respectively. If $M=S\otimes_RT$ and $m=1\otimes 1$
in an exact context $(\lambda,\mu, M,m)$, then we simply say that the pair $(\lambda,\mu)$ is exact.
Exact contexts can be easily constructed from rigid morphisms in an additive category (see Section \ref{sect3.0} below).

Given an exact context $(\lambda,\mu, M, m)$, we introduce, in Section \ref{sect3}, a new
multiplication $\circ$ on the abelian group $T\otimes_RS$, so that
$T\otimes_RS$ becomes an associative ring with identity and that
the following two maps
$$\rho: S\to T\otimes_RS,\; s\mapsto 1\otimes s\;\;\mbox{for}\;\;s\in S,\;\;\mbox{and}\;\;
\phi: T\to T\otimes_RS,\; t\mapsto t\otimes 1\;\;\mbox{for}\;\; t\in
T$$ are ring homomorphisms (see Lemma \ref{ring}). Furthermore, if both $S$ and $T$ are $R$-algebras over a commutative ring $R$ and if the pair $(\lambda,\mu)$ is exact, then this new ring structure on $T\otimes_RS$ coincides with the usual tensor product of the $R$-algebras $T$ and $S$ over $R$. Due to this reason, the new ring $(T\otimes_RS,\circ)$ is called the \emph{noncommutative tensor product} of the exact context $(\lambda,\mu, M, m)$, and denoted by $T\boxtimes_RS$ in this paper.
Note that if $(\lambda,\mu)$ is an exact pair, then the ring $T\boxtimes_RS$, together with $\rho$ and $\phi$, is actually the coproduct of the $R$-rings $S$ and $T$ (via the ring homomorphisms $\lambda$ and $\mu$) over $R$, and further, if $\lambda$ is a ring epimorphism, then $T\boxtimes_RS$ is isomorphic to the endomorphism ring of the $T$-module $T\otimes_RS$ (see Remark \ref{coproduct}).

\smallskip
Let $$B:=\left(\begin{array}{lc} S & M\\
0 & T\end{array}\right),\quad  C:=\left(\begin{array}{lc} T\boxtimes_RS & T\boxtimes_RS\\
T\boxtimes_RS& T\boxtimes_RS\end{array}\right).$$
Let $\beta: M\ra T\otimes_RS$ be the unique $R$-$R$-bimodule homomorphism such that $\phi=(m\cdot)\beta$ and $\rho=(\cdot m)\beta$ (see Section \ref{sect3.1}). We define a ring homomorphism
$$\theta:=\left(\begin{array}{cc} \rho & \beta\\
0 & \phi\end{array}\right): \; B \lra C. $$

First of all, this ring homomorphism is of particular interest in representation theory: The map $\theta$ can be regarded as the noncommutative localization of $B$ at a homomorphism between finitely
generated projective $B$-modules, and therefore it is a ring
epimorphism with $\Tor_1^B(C,C)=0$ (see Section \ref{sect5.1} and
\cite{Sch}), and yields a fully faithful exact functor $\theta_*:
C\Modcat\ra B\Modcat$, called the restriction functor, between the
category of all left $C$-modules and the one of all left
$B$-modules. Moreover, the map $\theta$ plays a fundamental role in
stratifications of derived categories and in algebraic $K$-theory
(see \cite{xc1, nr1, Ranicki}).

Generally speaking, $\theta$ is not always homological in the sense of Geigle and Lenzing (see \cite{GL}). In
$\cite{xc1}$, there is a sufficient condition for $\theta$ to be
homological. Concisely, if $\lambda:R\to S$ is an injective ring epimorphism with $\Tor_1^R(S,S)=0$ and if $T$ is the endomorphism ring of the $R$-module $S/R$ with $\mu: R\to T$  the ring homomorphism defined by
$r\mapsto(x\mapsto xr)$ for $r\in R$ and $x\in S/R$, then $B$ is
isomorphic to the endomorphism ring of the $R$-module $S\oplus
S\!/\!R$. For $\theta$ to be homological, we assume in \cite{xc1} that $_RS$ has projective dimension at most $1$.
In general context, it seems not much to be known about the map $\theta$ being homological.
So, the following general questions arise:

\smallskip
{\bf Questions.} {\it Let $(\lambda,\mu,M,m)$ be an exact context.

$(1)$ When is $\theta:B\to C$ homological, or equivalently, when is
the derived functor $D\big(\theta_*\big):\D{C}\to\D{B}$ fully
faithful ?
\smallskip

$(2)$ If $\theta$ is homological, is the Verdier quotient of $\D{B}$ by $\D{C}$ equivalent to the derived module category of a ring? or does $\D{B}$ admit a recollement
of derived module categories of rings $R$ and $C$?}

\smallskip The present paper will provide necessary and sufficient conditions to
these questions. Here, we will assume neither that $\lambda$ is injective,
nor that $_RS$ has projective
dimension at most $1$, nor that $\lambda$ is homological (compare with \cite{HJ, xc1}). Furthermore, we allow some flexibilities for
the choice of the ring homomorphism $\mu:R\to T$ and the bimodule
$M$. Our main result in this paper can be formulated as follows.

\begin{Theo}\label{th1}
Let $(\lambda, \mu, M, m)$ be an exact context. Then:

$(1)$ The following assertions are equivalent:

\quad $(a)$ The ring homomorphism $\theta:B\to C$  is homological.

\quad $(b)$ $\Tor^R_i(T,S)=0$ for all $i\geq 1$.

Moreover, if the pair $(\lambda,\mu)$ is exact and $\lambda$ is homological, then each of the above is equivalent to

\quad $(c)$ The ring homomorphism $\phi: T\ra T\boxtimes_RS$ is homological.

\medskip
$(2)$ If one of the above assertions in $(1)$ holds, then
there exists a recollement among the derived module categories of rings:

$$
\xymatrix@C=1.2cm{\D{T\boxtimes_RS}\ar[r]
&\D{B}\ar[r]\ar@/^1.2pc/[l]\ar@/_1.2pc/[l] &\D{R}.
\ar@/^1.2pc/[l]\ar@/_1.2pc/[l]}\vspace{0.2cm}$$

\medskip
\noindent If, in addition, the projective dimensions of $_RS$ and $T_R$ are finite, then the
above recollement can be restricted to a recollement of bounded
derived module categories:

$$
\xymatrix@C=1.2cm{\Db{T\boxtimes_RS}\ar[r]
&\Db{B}\ar[r]\ar@/^1.2pc/[l]\ar@/_1.2pc/[l] &\Db{R}.
\ar@/^1.2pc/[l]\ar@/_1.2pc/[l]}$$
\end{Theo}

\medskip
Note that $\D{B}$ is always a recollement of $\D{T}$ and $\D{S}$, in
which the derived category $\D{R}$ of the given ring $R$ is missing. However, Theorem
\ref{th1} provides us with a different recollement for $\D{B}$. A
remarkable feature of this recollemnt is that it contains $\D{R}$ as
a member, and thus provides a way to understand properties of the
ring $R$ through those of the rings closely related to $S$ and $T$. This idea will
be discussed in detail in the forthcoming papers \cite{xc5, xc6} of this
series.

The homological condition (b) in Theorem \ref{th1} can be satisfied in many cases. For instance, in commutative algebra, we may take $\lambda: R\ra S$ to be a localization, and in non-commutative case, we refer to the general examples in Section \ref{sect4.2}.

\smallskip
A realization of Theorem \ref{th1} occurs in noncommutative
localizations which have played an important role in topology (see \cite{Ranicki}).

\medskip
Given a ring homomorphism $\lambda:R\to S$, we may consider
$\lambda$ as a complex $\cpx{Q}$ of left $R$-modules with $R$ and
$S$ in degrees $-1$ and $0$, respectively. Then there is a
distinguished triangle $ R\lraf{\lambda} S\lraf{\pi}\cpx{Q}\lra
R[1]$ in the homotopy category $\K{R}$ of the category of all
$R$-modules. This triangle induces a canonical ring homomorphism
from $R$ to the endomorphism ring of $\cpx{Q}$ in $\K{R}$, and
therefore yields a ring homomorphism $\lambda'$ from $R$ to the
endomorphism ring of $\cpx{Q}$ in $\D{R}$, which depends on
$\lambda$ (see Section \ref{sect4} for details). Let $S':=\End_{\D
{R}}(\cpx{Q})$. Observe that if $\lambda$ is injective, then
$\cpx{Q}$ can be identified in $\D{R}$ with the $R$-module $S/R$,
and consequently, the map $\lambda':R\to S'$ coincides with the induced map
$R\to\End_R(S/R)$ by the right multiplication.

Further, let $\Lambda:=\End_{\D R}\big(S\oplus \cpx{Q}\big)$, and
let $\pi^*$ be the following induced map $$\Hom_{\D{R}}(S\oplus \cpx{Q},\,
\pi):\Hom_{\D{R}}(S\oplus \cpx{Q},\,S)\lra\Hom_{\D{R}}(S\oplus
\cpx{Q},\,\cpx{Q})$$ which is a homomorphism of finitely generated
projective $\Lambda$-modules. Let $\lambda_{\pi^*}: \Lambda\to
\Lambda_{\pi^*}$ stand for the noncommutative localization  of $\Lambda$
at $\pi^*$ (``universal localization" in terminology of Cohn and Schofield
\cite{cohenbook1,Sch}).

If $\lambda$ is a ring epimorphism such that $\Hom_R\big(S,\Ker(\lambda)\big)=0$, then we show in Section \ref{sect4} that the pair $(\lambda,\lambda')$ is exact. So, applying Theorem \ref{th1} to $(\lambda,\lambda')$, we get the following corollary.

\begin{Koro} \label{new}
If $\lambda: R\to S$ is a homological ring epimorphism such that
$\Hom_R\big(S,\Ker(\lambda)\big)=0$, then the following assertions
are equivalent:

 $(1)$ The noncommutative
localization $\lambda_{\pi^*}:\Lambda\to \Lambda_{\pi^*}$ of
$\Lambda$ at $\pi^*$ is homological.

$(2)$ The ring homomorphism $\phi: S'\ra S'\boxtimes_RS$ is
homological.

$(3)$ $\Tor^R_i(S',S)=0$ for any $i\geq 1$.

In particular, if one of the above assertions holds, then there
exists a recollement of derived module categories:

$$\xymatrix@C=1.2cm{\D{\End_{S'}(S'\otimes_RS)}\ar[r]&\D{\Lambda}\ar[r]
\ar@/^1.2pc/[l]\ar@/_1.2pc/[l]
&\D{R}\ar@/^1.2pc/[l]\ar@/_1.2pc/[l]}.\vspace{0.3cm}$$
\end{Koro}

As an application of Corollary \ref{new}, we obtain the following
result which not only generalizes the first statement of \cite[Corollary
6.6 (1)]{xc1} since we do not require that the ring epimorphism
$\lambda$ is injective, but also gives a way to get derived equivalences of rings (see \cite{Rk} for definition).

\begin{Koro}\label{tilting}
Let $\lambda: R\to S$ be a homological ring epimorphism such that
$\Hom_R\big(S,\Ker(\lambda)\big)=0$. Then we have the following:

$(1)$ If $_RS$ has projective dimension at most $1$, then
$\lambda_{\pi^*}:\Lambda\to \Lambda_{\pi^*}$ is homological.

$(2)$ The ring $\Lambda_{\pi^*}$ is zero if and only if there is an
exact sequence $ 0 \ra P_1\ra P_0\ra {}_RS \ra 0$ of $R$-modules
such that $P_i$ is finitely generated  and projective for $i=0,1$.
In this case, the rings $R$ and $\Lambda$ are derived equivalent.
\end{Koro}

As another application of Corollary \ref{new}, we have the following
result in which we do not impose any restriction
on the projective dimension of $_RS$.

\begin{Koro} \label{com1}
Suppose that $R\subseteq S$ is an extension of rings, that is, $R$
is a subring of the ring $S$ with the same identity. Let $S'$ be the
endomorphism ring of the $R$-module $S/R$ and
$B:=\left(\begin{array}{lc} S& \Hom_R(S,S/R)\\ 0 &
S'\end{array}\right).$

$(1)$ If the left $R$-module $S$ is flat, then there exists a recollement
of derived module categories:
$$
\xymatrix@C=1.2cm{\D{S'\boxtimes_RS}\ar[r]
&\D{B}\ar[r]\ar@/^1.2pc/[l]\ar@/_1.2pc/[l]
&\D{R}\ar@/^1.2pc/[l]\ar@/_1.2pc/[l]}$$ where $S'\boxtimes_RS$ is
the noncommutative tensor product of an exact context.

$(2)$ If $S$ is commutative and the inclusion $ R\ra S$ is
homological, then the ring $S'$ is commutative and
there exists a recollement of derived module categories:
$$\xymatrix@C=1.2cm{\D{S'\otimes_RS}\ar[r]&\D{B}\ar[r] \ar@/^1.2pc/[l]\ar@/_1.2pc/[l]
&\D{R}\ar@/^1.2pc/[l]\ar@/_1.2pc/[l]}\vspace{0.3cm}$$ where
$S'\otimes_RS$ is the usual tensor product of $R$-algebras.
\end{Koro}

\smallskip
Let us remark that, in commutative algebra, there is a lot of ring extensions satisfying the
`homological' assumption of Corollary \ref{com1} (2). For example, if $R$ is a commutative ring and $\Phi$ is a multiplicative subset of $R$ (that is, $\emptyset\ne \Phi$ and
$st\in \Phi$ whenever $s,t\in \Phi$), then the ordinary localization $R\ra
\Phi^{-1}R$ of $R$ at $\Phi$ is always homological. Further, if $f: R\ra R'$ is a homomorphism from the ring $R$ to another commutative ring $R'$, then the image of a multiplicative subset of $R$ under $f$ is again a multiplicative set in $R'$. So, as a consequence of
Corollary \ref{com1} (2), we obtain the following result which may
be of its own interest in commutative algebra.

\begin{Koro}\label{com2}
Suppose that $R$ is  a commutative ring with $\Phi$ a multiplicative
subset of $R$. Let $S$ be the localization $\Phi^{-1}R$ of $R$ at
$\Phi$, with $\lambda: R\ra S$ the canonical ring
homomorphism. If the map $\lambda$ is injective (for example, if $R$
is an integral domain), then there exists a recollement of derived
module categories:

$$\xymatrix@C=1.2cm{\D{\Psi^{-1}S'}\ar[r]&\D{\End_R(S\oplus S/R)}\ar[r]
\ar@/^1.2pc/[l]\ar@/_1.2pc/[l]
&\D{R}\ar@/^1.2pc/[l]\ar@/_1.2pc/[l]}\vspace{0.3cm}$$ where
$S':=\End_R(S/R)$, and $\Psi$ is the image of $\Phi$ under the induced map
$R\to S'$ given by the right multiplication.
\end{Koro}

Observe that the recollements in Corollaries \ref{com1} and \ref{com2} occur in the study of infinitely generated tilting modules (see \cite{HJ} and \cite{xc1}).

\medskip
The contents of this paper are outlined as follows. In Section 2, we
fix notation and recall some definitions and basic facts which will
be used throughout the paper. In particular, we shall recall the
definitions of noncommutative localizations, coproducts of rings and
recollements, and prepare several lemmas for our proofs. In Section
\ref{sect3.0}, we introduce the notion of exact contexts. To
construct exact contexts, we introduce rigid morphisms or
hypercyclic bimodules, and show that rigid morphisms exist almost
everywhere in representation theory. For example, all kinds of
approximations are rigid morphisms. Thus, exact contexts exist
rather abundantly. In Section \ref{sect3}, we define the so-called
noncommutative tensor products of exact contexts, which will
characterize the left parts of recollements constructed in Section
\ref{sect5}. Also, we provide examples to demonstrate that
noncommutative tensor products cover many well-known constructions
in noncommutative algebra. In Section \ref{sect5}, we prove Theorem
\ref{th1} and all of its corollaries mentioned in Section
\ref{sect1}. Finally, in Section \ref{sect6}, we give several
examples to explain the necessity of some assumptions in our
results.

In the second paper \cite{xc5} of this series, we shall consider the
algebraic $K$-theory of recollements, and establish a long Mayer-Vietoris sequence of higher
algebraic $K$-groups for homological Milnor squares. In the third paper \cite{xc6}, we shall study relationships among finitistic dimensions of three algebras involved in a recollement. This will extend an earlier result of Happel and a recent result by Xu.

\section{Preliminaries\label{sect2}}
In this section, we shall recall some definitions, notation and
basic results which are closely related to our proofs.

\subsection{Notation and basic facts on derived categories}\label{sect2.1}

Let $\mathcal C$ be an  additive category.

Throughout the paper, a full subcategory $\mathcal B$ of $\mathcal
C$ is always assumed to be closed under isomorphisms, that is, if
$X\in {\mathcal B}$ and $Y\in\cal C$ with $Y\simeq X$, then
$Y\in{\mathcal B}$.

Given two morphisms $f: X\to Y$ and $g: Y\to Z$ in $\mathcal C$, we
denote the composite of $f$ and $g$ by $fg$ which is a morphism from
$X$ to $Z$. The induced morphisms $\Hom_{\mathcal
C}(Z,f):\Hom_{\mathcal C}(Z,X)\ra \Hom_{\mathcal C}(Z,Y)$ and
$\Hom_{\mathcal C}(f,Z): \Hom_{\mathcal C}(Y, Z)\ra \Hom_{\mathcal
C}(X, Z)$ are denoted by $f^*$ and $f_*$, respectively.

We denote the composition of a functor $F:\mathcal {C}\to
\mathcal{D}$ between categories $\mathcal C$ and $\mathcal D$ with a
functor $G: \mathcal{D}\to \mathcal{E}$ between categories $\mathcal
D$ and $\mathcal E$ by $GF$ which is a functor from $\mathcal C$ to
$\mathcal E$. The kernel and the image of the functor $F$ are
denoted by $\Ker(F)$ and $\Img(F)$, respectively.

Let $\mathcal{Y}$ be a full subcategory of $\mathcal{C}$. By
$\Ker(\Hom_{\mathcal{C}}(-,\mathcal{Y}))$ we denote the full subcategory of $\mathcal{C}$ which is left
orthogonal to $\mathcal{Y}$, that is, the
full subcategory of $\mathcal{C}$ consisting of the objects $X$ such
that $\Hom_{\mathcal{C}}(X,Y)=0$ for all objects $Y$ in
$\mathcal{Y}$. Similarly, $\Ker(\Hom_{\mathcal{C}}(\mathcal{Y},-))$
stands for the right orthogonal subcategory in $\cal C$ with respect
to $\mathcal{Y}$.

Let $\C{\mathcal{C}}$ be the category of all complexes over
$\mathcal{C}$ with chain maps, and $\K{\mathcal{C}}$ the homotopy
category of $\C{\mathcal{C}}$. When $\mathcal{C}$ is abelian, the
derived category of $\mathcal{C}$ is denoted by $\D{\mathcal{C}}$,
which is the localization of $\K{\mathcal C}$ at all
quasi-isomorphisms. It is well known that both $\K{\mathcal{C}}$ and
$\D{\mathcal{C}}$ are triangulated categories. For a triangulated
category, its shift functor is denoted by $[1]$ universally.

If $\mathcal{T}$ is a triangulated category with small coproducts
(that is, coproducts indexed over sets exist in ${\mathcal T}$),
then, for each object $U$ in $\mathcal{T}$, we denote by ${\rm
Tria}(U)$ the smallest full triangulated subcategory of
$\mathcal{T}$ containing $U$ and being closed under small
coproducts. We mention the following properties related to
Tria$(U)$:

Let $F: \mathcal{T}\ra \mathcal{T}'$ be a triangle functor of
triangulated categories, and let $\mathcal Y$ be a full subcategory
of $\mathcal{T}'$. We define $F^{-1}\mathcal{Y}:=\{X\in
\mathcal{T}\mid F(X)\in \mathcal{Y}\}$. Then

(1) If $\mathcal Y$ is a triangulated subcategory, then
$F^{-1}\mathcal{Y}$ is a full triangulated subcategory of
$\mathcal{T}$.

(2) Suppose that $\mathcal{T}$ and $\mathcal{T}'$ admit small
coproducts and that $F$ commutes with coproducts. If $\mathcal Y$ is
closed under small coproducts in $\mathcal{T}'$, then
$F^{-1}\mathcal{Y}$ is closed under small coproducts in
$\mathcal{T}$. In particular, for an object $U\in \mathcal{T}$, we
have $F(\Tria(U))\subseteq \Tria(F(U))$.

\medskip
In this paper, all rings considered are assumed to be associative
and with identity, and all ring homomorphisms preserve identity.
Unless stated otherwise, all modules are referred to left modules.

Let $R$ be a ring. We denote by $R\Modcat$ the
category of all unitary left $R$-modules. By our convention of the
composite of two morphisms, if $f: M\ra N$ is a homomorphism of
$R$-modules, then the image of $x\in M$ under $f$ is denoted by
$(x)f$ instead of $f(x)$. The endomorphism ring of the $R$-module
$M$ is denoted by $\End_R(M)$.

As usual, we shall simply write $\C{R}$, $\K{R}$ and $\D{R}$ for
$\C{R\Modcat}$, $\K{R\Modcat}$ and $\D{R\Modcat}$, respectively, and
identify $R\Modcat$ with the subcategory of $\D{R}$ consisting of
all stalk complexes concentrated in degree zero. Further, we denote by $\Db{R}$
the full subcategory of $\D{R}$ consisting of all complexes which are isomorphic in $\D{R}$ to bounded complexes of $R$-modules¡£

Let $(\cpx{X},d_{\cpx{X}})$ and $(\cpx{Y},d_{\cpx{Y}})$ be two chain
complexes over $R\Modcat$. The mapping cone of a chain map
$\cpx{h}:\cpx{X}\to\cpx{Y}$ is usually denoted by
$\Cone{(\cpx{h})}$. In particular, we have a triangle
$\cpx{X}\lraf{\cpx{h}}\cpx{Y}\lra\Cone(\cpx{h})\lra\cpx{X}[1]$ in
$\K{R}$, called a \emph{distinguished triangle}. For each
$n\in\mathbb{Z}$, we denote by $H^n(-):\D{R}\to R\Modcat$ the $n$-th
cohomology functor. Certainly, this functor is naturally isomorphic
to the Hom-functor $\Hom_{\D{R}}(R,-[n])$.

The Hom-complex $\cpx{\Hom}_R(\cpx{X},\cpx{Y})$ of $\cpx{X}$ and
$\cpx{Y}$ is defined to be the complex
$\big(\Hom_R^n(\cpx{X},\cpx{Y}),
d^{\,n}_{\cpx{X},\cpx{Y}}\big)_{n\in\mathbb{Z}}$ with
$$\Hom_R^n(\cpx{X},\cpx{Y}):=\prod_{p\in\mathbb{Z}}\Hom_R(X^p,Y^{p+n})$$
and the differential $d^{\,n}_{\cpx{X},\cpx{Y}}$ of degree $n$ given
by
$$(f^p)_{p\in\mathbb{Z}}\mapsto \big(f^p
d_{\cpx{Y}}^{p+n}-(-1)^{n}d_{\cpx{X}}^p
f^{p+1}\big)_{p\in\mathbb{Z}}$$ for
$(f^p)_{p\in\mathbb{Z}}\in\Hom_R^n(\cpx{X},\cpx{Y})$. For example,
if $X\in R\Modcat$, then we have
$$\cpx{\Hom}_R(X,\cpx{Y})=\big(\Hom_R(X,Y^n),\Hom_R(X,d_{\cpx{Y}}^n)\,\big)_{n\in\mathbb{Z}};$$
if $Y\in R\Modcat$, then
$$\cpx{\Hom}_R(\cpx{X},\,Y)=\big(\Hom_R(X^{-n},\,Y),\,(-1)^{n+1}\Hom_R(d_{\cpx{X}}^{-n-1},\,Y)\,\big)_{n\in\mathbb{Z}}.$$
For simplicity, we denote $\cpx{\Hom}_R(X,\cpx{Y})$ and
$\cpx{\Hom}_R(\cpx{X},\,Y)$ by $\Hom_R(X,\cpx{Y})$ and
$\Hom_R(\cpx{X},\,Y)$, respectively.  Note that
$\Hom_R(\cpx{X},\,Y)$ is also isomorphic to the complex
$\big(\Hom_R(X^{-n},\,Y),\,\Hom_R(d_{\cpx{X}}^{-n-1},\,Y)\,\big)_{n\in\mathbb{Z}}.$

Moreover, it is known that $H^n(\cpx{\Hom}_R(\cpx{X},\cpx{Y}))\simeq
\Hom_{\K R}(\cpx{X},\cpx{Y}[n])$ for any $n\in\mathbb{Z}$.

Let $\cpx{Z}$ be a chain complex over $R^{op}\Modcat$. Then the
tensor complex $\cpx{Z}\cpx{\otimes}_R\cpx{X}$ of $\cpx{Z}$ and
$\cpx{X}$ over $R$ is defined to be the complex
$\big(\cpx{Z}{\otimes}^n_R\cpx{X},\partial^n_{\cpx{Z},\cpx{X}}\big)_{n\in\mathbb{Z}}$
with
$$\cpx{Z}{\otimes}^n_R\cpx{X}:=\bigoplus_{p\in\mathbb{Z}}Z^p\otimes_RX^{n-p}$$
and the differential $\partial_{\cpx{Z},\cpx{X}}$ of degree $n$
given by
$$z\otimes x \mapsto(z)d^p_{\cpx{Z}}\otimes x
+(-1)^p z\otimes(x)d_{\cpx{X}}^{n-p}$$ for  $z\in Z^p $ and $x\in
X^{n-p}$. For instance, if $X\in R\Modcat$, then
$\cpx{Z}\cpx{\otimes}_RX=\big(Z^n\otimes_RX, d_{\cpx{Z}}^n\otimes
1\big)_{n\in\mathbb{Z}}$. In this case, we denote
$\cpx{Z}\cpx{\otimes}_RX$ simply by $\cpx{Z}\otimes_RX$.

The following result establishes a relationship between
Hom-complexes and tensor complexes.

Let $S$ be an arbitrary ring. Suppose that
$\cpx{X}=(X^n,d^n_{\cpx{X}})$ is a bounded complex of
$R$-$S$-bimodules. If $_RX^n$ is finitely generated and projective
for all $n\in\mathbb{Z}$, then there is a natural isomorphism of
functors:
$$\Hom_R(\cpx{X}, R)\cpx{\otimes}_R-\,\lraf{\simeq}\cpx{\Hom}_R(\cpx{X},-):\C{R}\to\C{S}.$$

To prove this, we note that, for any $R$-$S$-bimodule $X$ and any
$R$-module $Y$, there is a homomorphism of $S$-modules:
$\delta_{X,Y}:\Hom_R(X,R)\otimes_RY\lra\Hom_R(X,Y)$ defined by
$f\otimes y\mapsto [\,x\mapsto (x)f y\,]$ for $f\in\Hom_R(X,R)$,
$y\in Y$ and $x\in X$, which is natural in both $X$ and $Y$.
Moreover, the map $\delta_{X,Y}$ is an isomorphism if $_RX$ is
finitely generated and projective. For any $\cpx{Y}\in \C{R}$ and
any  $n\in\mathbb{Z}$, it is clear that
$$\Hom_R(\cpx{X},\,R)\displaystyle{\otimes}^n_R\cpx{Y}=
\bigoplus_{p\in\mathbb{Z}}\Hom_R(X^{-p},R)\otimes_RY^{n-p}\;\,\mbox{and}
\;\,\Hom_R^n(\cpx{X},\cpx{Y})=\bigoplus_{p\in\mathbb{Z}}\Hom_R(X^p,Y^{p+n})$$
since $\cpx{X}$ is a bounded complex.  Now, we define
$\Delta_{\cpx{X},\,\cpx{Y}}^n:=\sum_{p\in\mathbb{Z}}{(-1)}^{p(n-p)}\delta_{X^{-p},\,Y^{n-p}}$,
which is a homomorphism of $S$-modules from
$\Hom_R(\cpx{X},\,R){\otimes}^n_R\cpx{Y}$ to
$\Hom_R^n(\cpx{X},\cpx{Y})$. Then, one can check that
$\cpx{\Delta_{\cpx{X},\,\cpx{Y}}}:=\big(\Delta_{\cpx{X},\,\cpx{Y}}^n\big)_{n\in\mathbb{Z}}$
is a chain map from $\Hom_R(\cpx{X}, R)\cpx{\otimes}_R\cpx{Y}$ to
$\cpx{\Hom}_R(\cpx{X},\cpx{Y})$. Since $_RX^{-p}$ is finitely
generated and projective for each $p\in\mathbb{Z}$, the map
$\delta_{X^{-p},\,Y^{n-p}}$ is an isomorphism, and so is the map
$\Delta_{\cpx{X},\,\cpx{Y}}^n$. This implies that
$$\cpx{\Delta_{\cpx{X},\,\cpx{Y}}}:\Hom_R(\cpx{X}, R)\cpx{\otimes}_R\cpx{Y}\lra \cpx{\Hom}_R(\cpx{X},\cpx{Y})$$
is an isomorphism in $\C{S}$.  Since the homomorphism $\delta_{X,Y}$
is natural in the variables $X$ and $Y$, it can be checked directly
that
$$\cpx{\Delta_{\cpx{X},\,-}}:\Hom_R(\cpx{X}, R)\cpx{\otimes}_R-\lra \cpx{\Hom}_R(\cpx{X},\,-)$$
defines a natural isomorphism of functors from $\C{R}$ to $\C{S}$.

In the following, we shall recall some definitions and basic facts
about derived functors defined on derived module categories. For
details and proofs, we refer to \cite{bn, Keller}.

Let $\K{R}_P$ (respectively, $\K{R}_I$) be the smallest full
triangulated subcategory of $\K{R}$ which

(i) contains all the bounded above (respectively, bounded below)
complexes of projective (respectively, injective) $R$-modules, and

(ii) is closed under arbitrary direct sums (respectively, direct
products).

\medskip
Note that $\K{R}_P$ is contained in $\K{\Pmodcat{R}}$, where
$\Pmodcat{R}$ is the full subcategory of $R$-Mod consisting of all
projective $R$-modules. Moreover, the composition functors
$$\K{R}_P\hookrightarrow\K{R}\to\D{R}\quad\mbox{and}\quad
\K{R}_I\hookrightarrow\K{R}\to\D{R}$$ are equivalences of
triangulated categories. This means that, for each complex $\cpx{X}$
in $\D{R}$, there exists a complex $_p\cpx{X}\in \K{R}_P$ together
with a quasi-isomorphism $_p\cpx{X}\to\cpx{X}$, as well as a complex
${}_i\cpx{X}\in \K{R}_I$ together with a quasi-isomorphism
$\cpx{X}\to{_i}\cpx{X}$. In this sense, we shall simply call
$_p\cpx{X}$ the projective resolution of $\cpx{X}$ in $\K{R}$. For
example, if $X$ is  an $R$-module, then we can choose $_pX$ to be a
deleted projective resolution of $_RX$.

Furthermore, if either $\cpx{X}\in\K{R}_P$ or $\cpx{Y}\in\K{R}_I$,
then
$\Hom_{\K{R}}(\cpx{X},\cpx{Y})\simeq\Hom_{\D{R}}(\cpx{X},\cpx{Y})$,
and this isomorphism is induced by the canonical localization
functor from $\K{R}$ to $\D{R}$.

For any triangle functor $H:\K{R}\to\K{S}$, there is a total
left-derived functor ${\mathbb L}H:\D{R}\to\D{S}$ defined by
$\cpx{X}\mapsto H(_p\cpx{X})$, a total right-derived functor
${\mathbb R}H:\D{R}\to\D{S}$ defined by $\cpx{X}\mapsto
H(_i\cpx{X})$. Observe that, if $H$ preserves acyclicity, that is,
$H(\cpx{X})$ is acyclic whenever $\cpx{X}$ is acyclic, then $H$
induces a triangle functor $D(H):\D{R}\to\D{S}$ defined by
$\cpx{X}\mapsto H(\cpx{X})$. In this case, we have ${\mathbb
L}H={\mathbb R}H=D(H)$ up to natural isomorphism, and $D(H)$ is then
called the derived functor of $H$.

Let $\cpx{M}$ be a complex of $R$-$S$-bimodules. Then the functors
$$\cpx{M}\cpx{\otimes}_S-:\K{S}\to\K{R}\quad\mbox{and}\quad\cpx{\Hom}_R(\cpx{M},-):\K{R}\to\K{S}$$
form a pair of adjoint triangle functors. Denote by
$\cpx{M}\otimesL_S-$ the total left-derived functor of
$\cpx{M}\cpx{\otimes}_S-$, and by ${\mathbb R}\Hom_R(\cpx{M},-)$ the
total right-derived functor of $\cpx{\Hom}_R(\cpx{M},-)$. It is
clear that $\big(\cpx{M}\otimesL_S-, {\mathbb
R}\Hom_R(\cpx{M},-)\big)$ is an adjoint pair of triangle functors.
Further, the corresponding counit adjunction
$$\varepsilon: \cpx{M}\otimesL_S\,{\mathbb
R}\Hom_R(\cpx{M},-)\lra Id_{\D{R}}$$ is given by the composite of
the following canonical morphisms in $\D{R}$:
$\cpx{M}\otimesL_S{\mathbb
R}\Hom_R(\cpx{M},\cpx{X})=\cpx{M}\otimesL_S\cpx{\Hom}_R(\cpx{M},\,{_i\cpx{X}})=
\cpx{M}\cpx{\otimes}_S\big(_p\cpx{\Hom}_R(\cpx{M},\,{_i\cpx{X}})\big)\lra
\cpx{M}\cpx{\otimes}_S\cpx{\Hom}_R(\cpx{M},\,{_i\cpx{X}})\lra{_i\cpx{X}}\lraf{\simeq}\cpx{X}.$
Similarly, we have a corresponding unit adjunction $\eta:
Id_{\D{S}}\lra {\mathbb R}\Hom_R(\cpx{M},\,\cpx{M}\otimesL_S-)$,
which is given by the following composites for $\cpx{Y}\in \D{S}$: $
\cpx{Y}\lraf{\simeq} {}_p\cpx{Y} \lra
\cpx{\Hom}_R(\cpx{M},\cpx{M}\cpx{\otimes}_S(_p\cpx{Y})) \lra
\cpx{\Hom}_R(\cpx{M},
{}_i(\cpx{M}\cpx{\otimes}_S(_p\cpx{Y})))=\mathbb{R}\Hom_R(\cpx{M},
\cpx{M}\cpx{\otimes}_S(_p\cpx{Y}))$ =
$\mathbb{R}\Hom_R(\cpx{M},\cpx{M}\otimes_S^{\mathbb{L}}\cpx{Y}).$

For $\cpx{X}\in \D{R}$ and $n\in\mathbb{Z}$, we have
$H^n(\mathbb{R}\Hom_R(\cpx{M},\cpx{X}))=H^n(\cpx{\Hom}_R(\cpx{M},
{}_i\cpx{X})) \simeq\Hom_{\K{R}}(\cpx{M}, {}_i\cpx{X}[n])$
$\simeq\Hom_{\D{R}}(\cpx{M}, {}_i\cpx{X}[n]) \simeq
\Hom_{\D{R}}(\cpx{M},\cpx{X}[n]).$

Let $T$ be another ring and $\cpx{N}$ a complex of
$S$-$T$-bimodules. If ${_S}\cpx{N}\in\K{S}_P$, then
$$\cpx{M}\otimesL_S(\cpx{N}\otimesL_T-)\lraf{\simeq}
(\cpx{M}\otimesL_S\cpx{N})\otimesL_T-\,=
(\cpx{M}\otimes_S\cpx{N})\otimesL_T-:\D{T}\lra\D{R}$$ In fact, since
${_S}\cpx{N}\in\K{S}_P$ by assumption, we have
$\cpx{N}\cpx{\otimes}_T(_p\cpx{W})\in\K{S}_P$ for $\cpx{W}\in\D{T}$.
It follows that $\cpx{M}\otimesL_S(\cpx{N}\otimesL_T\cpx{W})=
\cpx{M}\otimesL_S\big(\cpx{N}\cpx{\otimes}_T(_p\cpx{W})\big)
=\cpx{M}\cpx{\otimes}_S\big(\cpx{N}\cpx{\otimes}_T(_p\cpx{W})\big)\simeq
\big(\cpx{M}\cpx{\otimes}_S\cpx{N}\big)\cpx{\otimes}_T(_p\cpx{W})
=(\cpx{M}\otimes_S\cpx{N})\otimesL_T\cpx{W}=(\cpx{M}\otimesL_S\cpx{N})\otimesL_T\cpx{W}$.

\subsection{Homological ring epimorphisms and recollements\label{sect2.2}}

Let $\lambda: R\ra S$ be a homomorphism of rings.

We denote by $\lambda_*:S\Modcat\to R\Modcat$ the restriction
functor induced by $\lambda$, and by $D(\lambda_*):\D{S}\to\D{R}$
the derived functor of the exact functor $\lambda_*$. We say that
$\lambda$ is a \emph{ring epimorphism} if the restriction functor
$\lambda_*:S\Modcat\to R\Modcat$ is fully faithful. It is proved
that $\lambda$ is a ring epimorphism if and only if the
multiplication map $S\otimes_RS\ra S$ is an isomorphism as
$S$-$S$-bimodules if and only if, for any two homomorphisms
$f_1,f_2: S\ra T$ of rings, the equality $\lambda f_1=\lambda f_2$
implies that $f_1=f_2$. This means that, for a ring epimorphism, we
have $X\otimes_SY\simeq X\otimes_RY$ and $\Hom_S(Y,Z)\simeq
\Hom_R(Y,Z)$ for all right $S$-modules $X$, and for all $S$-modules
$Y$ and $Z$. Note that, for a ring epimorphism $\lambda: R\ra S$, if
$R$ is commutative, then so is $S$.

Following \cite{GL}, a ring epimorphism $\lambda: R\ra S$ is called
\emph{homological} if $\;\Tor^R_i(S, S)=0$ for all $i>0$. Note that
a ring epimorphism $\lambda$ is homological if and only if the
derived functor $D(\lambda_*):\D{S}\to\D{R}$ is fully faithful. This
is also equivalent to saying that $\lambda$ induces an isomorphism
$S\otimesL_RS\simeq S$ in $\D{S}$. Moreover, for a homological ring
epimorphism, we have $\Tor_i^R(X,Y)\simeq \Tor^S_i(X,Y)$ and
$\Ext_S^i(Y,Z)\simeq \Ext^i_R(Y,Z)$ for all $i\ge 0$ and all right
$S$-modules $X$, and for all $S$-modules $Y$ and $Z$ (see
\cite[Theorem 4.4]{GL}).

Clearly, if  $\lambda: R\ra S$ is a ring epimorphism such that
either $_RS$ or $S_R$ is flat, then $\lambda$ is homological. In
particular, if $R$ is commutative and $\Phi$ is a multiplicative
subset of $R$, then the canonical ring homomorphism $R\to
\Phi^{-1}R$ is homological, where $\Phi^{-1}R$ stands for the
(ordinary) localization of $R$ at $\Phi$.

As a generalization of localizations of commutative rings, noncommutative (``universal" in Cohen's terminology)
localizations of arbitrary rings were introduced in
\cite{cohenbook1} (see also \cite{Sch}) and provide a class of ring
epimorphisms with vanishing homology for the first degree.
Now we mention the following basic fact
on noncommutative localizations.

\begin{Lem}{\rm (see \cite{cohenbook1}, \cite{Sch})} Let $R$ be a ring and let
$\Sigma$ be a set of homomorphisms between finitely generated
projective $R$-modules. Then there is a ring $R_{\Sigma}$ and a
homomorphism  $\lambda_{\Sigma}: R\to R_{\Sigma}$ of rings such that

$(1)$ $\lambda_{\Sigma}$ is $\Sigma$-inverting, that is, if
$\alpha:P\to Q$ belongs to $\Sigma$, then
$R_{\Sigma}\otimes_R\alpha:R_{\Sigma}\otimes_RP\to
R_{\Sigma}\otimes_R Q$ is an isomorphism of $R_{\Sigma}$-modules,
and

$(2)$ $\lambda_{\Sigma}$ is universal $\Sigma$-inverting, that is,
if $S$ is a ring such that there exists a $\Sigma$-inverting
homomorphism $\varphi:R\to S$, then there exists a unique
homomorphism $\psi:R_{\Sigma}\to S$ of rings such that
$\varphi=\lambda\psi$.

$(3)$ $\lambda_{\Sigma}:R\to R_{\Sigma}$  is a ring epimorphism with
$\Tor^R_1(R_{\Sigma}, R_{\Sigma})=0.$ \label{lem2.3}
\end{Lem}

Following \cite{nr1}, the $\lambda_{\Sigma}: R\to R_{\Sigma}$ in Lemma \ref{lem2.3} is
called the \emph{noncommutative localization} of $R$ at $\Sigma$. One
should be aware that $R_{\Sigma}$ may not be flat as a right or left
$R$-module. Even worse, the map $\lambda_{\Sigma}$ in  general is
not homological (see \cite{nrs}). Thus it is a fundamental question
when $\lambda_{\Sigma}$ is homological.

\smallskip
Next, we recall the definition of coproducts of rings defined by
Cohn in \cite{cohn}, and point out that noncommutative localizations are
preserved by taking coproducts of rings.

Let $R_0$ be a ring. An \emph{$R_0$-ring} is a ring $R$ together
with a ring homomorphism $\lambda_R: R_0\ra R$. An
\emph{$R_0$-homomorphism} from an $R_0$-ring $R$ to another
$R_0$-ring $S$ is a ring homomorphism $f: R\ra S$ such that
$\lambda_S=\lambda_Rf$. Then we can form the category of $R_0$-rings
with $R_0$-rings as objects and with $R_0$-morphisms as morphisms.
Clearly, epimorphisms of this category are exactly ring epimorphisms
starting from $R_0$.

The \emph{coproduct} of a family $\{R_i\mid i\in I\}$ of $R_0$-rings
with $I$ an index set is defined to be an $R_0$-ring $R$ together
with a family $\{\rho_i: R_i\ra R\mid i\in I \}$ of
$R_0$-homomorphisms such that, for  any $R_0$-ring $S$ with a family
of $R_0$-homomorphisms $\{\tau_i: R_i\ra S\mid i\in I\}$, there
exists a unique $R_0$-homomorphism $\delta: R\ra S$ such that
$\tau_i=\rho_i\delta$ for all $i\in I$.

It is well known that the coproduct of a family $\{R_i\mid i\in I\}$
of $R_0$-rings always exists. We denote this coproduct by
$\sqcup_{R_0}R_{i}$. Note that if $I=\{1,2\}$, then $R_1\sqcup_{R_0}R_2$ is the push-out in the category of $R_0$-rings. This implies that if $\lambda_{R_1}:R_0\to R_1$ is a ring epimorphism, then so is the homomorphism $\rho_2: R_2\ra R_1\sqcup_{R_0}R_2$. Moreover, $R_0\sqcup_{R_0}R_1=R_1=R_1\sqcup_{R_0}R_0$ for every $R_0$-ring $R_1$, where $\lambda_{R_0}:R_0\to R_0$ is the identity.

In general, the coproduct of two $R_0$-algebras may not be
isomorphic to their tensor product over $R_0$. For example, given a
field $k$, the coproduct over $k$ of the polynomial rings $k[x]$ and
$k[y]$ is the free ring $k\langle x,y\rangle$ in two variables $x$
and $y$, while the tensor product over $k$ of $k[x]$ and $k[y]$ is
the polynomial ring $k[x,y]$.

The following result is taken from \cite[Lemma 6.2]{xc1} and will be used later.

\begin{Lem}\label{epi01}
Let $R_0$ be a ring, $\Sigma$ a set of homomorphisms between
finitely generated projective $R_0$-modules, and $\lambda_{\Sigma}:
R_0\to R_1:=(R_0)_{\Sigma}$ the noncommutative localization of $R_0$ at
$\Sigma$. Then, for any $R_0$-ring $R_2$, the coproduct
$R_1\sqcup_{R_0}R_2$ is isomorphic to the noncommutative localization
$(R_2)_{\Delta}$ of $R_2$ at the set
$\Delta:=\{R_2\otimes_{R_0}f\mid f\in\Sigma\}$.
\end{Lem}

Finally, we recall the notion of recollements of triangulated
categories, which was first defined 
in \cite{BBD} to study ``exact sequences" of derived categories of
coherent sheaves over geometric objects.

\begin{Def}\label{def01} \rm
Let  $\mathcal{D}$, $\mathcal{D'}$ and $\mathcal{D''}$ be
triangulated categories. We say that $\mathcal{D}$ is a
\emph{recollement} of $\mathcal{D'}$ and $\mathcal{D''}$ if there
are six triangle functors among the three categories:
$$\xymatrix{\mathcal{D''}\ar^-{i_*=i_!}[r]&\mathcal{D}\ar^-{j^!=j^*}[r]
\ar_-{i^!}@/^1.2pc/[l]\ar_-{i^*}@/_1.6pc/[l]
&\mathcal{D'}\ar_-{j_*}@/^1.2pc/[l]\ar_-{j_!}@/_1.6pc/[l]}$$ such
that

$(1)$ $(i^*,i_*),(i_!,i^!),(j_!,j^!)$ and $(j^*,j_*)$ are adjoint
pairs,

$(2)$ $i_*,j_*$ and $j_!$ are fully faithful functors,

$(3)$ $i^!j_*=0$ (and thus also $j^! i_!=0$ and $i^*j_!=0$), and

$(4)$ for each object $X\in\mathcal{D}$, there are two triangles in
$\mathcal D$:
$$
i_!i^!(X)\lra X\lra j_*j^*(X)\lra i_!i^!(X)[1],\quad j_!j^!(X)\lra
X\lra i_*i^*(X)\lra j_!j^!(X)[1].
$$

\end{Def}

Clearly, it follows from definition that, for any objects
$X\in\mathcal{D}'$ and $Y\in\mathcal{D''}$, we have
$$\Hom_{\mathcal{D}}\big(j_!(X),i_*(Y)\big)=0=\Hom_{\mathcal{D}}\big(i_*(Y),j_*(X)\big).$$

A typical example of recollements of derived module categories is
given by triangular matrix rings: Suppose that $A$ and $B$ are
rings, and that
$N$ is an $A$-$B$-bimodule. Let  $R=\left(\begin{array}{lc} A & N\\
0 & B\end{array}\right)$ be the triangular matrix ring associated
with $A,B$ and $N$. Then there is a recollement of derived module
categories:

$$\xymatrix@C=1.2cm{\D{A}\ar[r]&\D{R}\ar[r]
\ar@/^1.2pc/[l]\ar@/_1.2pc/[l]
&\D{B}\ar@/^1.2pc/[l]\ar@/_1.2pc/[l]}.\vspace{0.3cm}$$ In this case,
the six triangle functors in Definition \ref{def01} can be described
explicitly:

Let $e:=\left(\begin{array}{ll} 0 &0\\
0 & 1\end{array}\right)\in R$. Then we have $$j_!=Re\otimesL_B-,\,
j^!=eR\otimesL_R-,\, j_*={\mathbb R}\Hom_B(eR,-),\,
i^*=A\otimesL_R-,\, i_*=A\otimesL_A-,\, i^!={\mathbb
R}\Hom_R(A,-),$$ where $A$ is identified with $R/ReR$. Note that the
canonical surjection $R\to R/ReR$ is always a homological ring
epimorphism.

As a further generalization of the above situation, it was shown in
\cite[Section 4]{NS} that, for an arbitrary homological ring
epimorphism $\lambda:R\to S$, there is a recollement of triangulated
categories:

$$\xymatrix@C=1.2cm{\D{S}\ar[r]&\D{R}\ar[r]
\ar@/^1.2pc/[l]\ar@/_1.2pc/[l]
&{\rm{Tria}}(\cpx{Q})\ar@/^1.2pc/[l]\ar@/_1.2pc/[l]}\vspace{0.3cm}$$
where $\cpx{Q}$ is given by the distinguished triangle
$R\lraf{\lambda} S\lra\cpx{Q}\lra R[1]$ in $\D{R}$.  In this case,
the functor $j_!$ is the canonical embedding and
$$j^!=(\cpx{Q}[-1])\otimesL_R-,\,
i^*=S\otimesL_R-, \,i_*=S\otimesL_S-,\,i^!={\mathbb
R}\Hom_R(_RS,-).$$ Moreover, we have
$$\D{S}\simeq\Ker\big(\Hom_{\D{R}}({\rm{Tria}}(\cpx{Q}),-)\big):=\{\cpx{X}\in
\D{R}\mid \Hom_{\D{R}}(Y,\cpx{X})=0 \mbox{\, for  all\;}
Y\in\rm{Tria(\cpx{Q})}\}.$$ This clearly implies that
$\Hom_{\D{R}}(\cpx{Q},\cpx{X}[n])=0$ for all $\cpx{X}\in\C{S}$ and
$n\in\mathbb{Z}$.

\section{Definitions of rigid morphisms and exact contexts\label{sect3.0}}

In this section we introduce the notion of rigid morphisms in an
additive category, which occur almost everywhere in the
representation theory of algebras, and which will be used to
construct exact contexts.

Let $\mathcal C$ be an additive category. An object $\cpx{X}$ in
$\mathscr{C}(\mathcal{C})$ is \emph {rigid} if
$\Hom_{\mathscr{K}(\mathcal{C})}(\cpx{X},\cpx{X}[1])=0$.

\begin{Def}\label{def1}
A morphism $\cpx{f}: \cpx{Y}\ra \cpx{X}$ in
$\mathscr{C}(\mathcal{C})$ is said to be rigid if the object
$\cpx{Z}$ in a distinguished triangle
$\cpx{Y}\lraf{\cpx{f}}\cpx{X}\ra \cpx{Z}\ra\cpx{Y}[1]$ is rigid, or
equivalently, the mapping cone $\Cone(\cpx{f})$ of $\cpx{f}$ is
rigid in $\mathscr{C}(\mathcal{C})$.

A morphism $f: Y\ra X$ in $\mathcal{C}$ is said to be rigid if $f$,
considered as a morphism from the stalk complex $Y$ to the stalk
complex $X$, is rigid, or equivalently, the complex $\Cone(f): 0\ra Y\lraf{f}
X\ra 0$ is rigid in $\mathscr{C}(\mathcal{C})$.
\end{Def}

Note that the rigidity of a morphism $\cpx{f}$ does not depend on
the choice of the triangle which extends $\cpx{f}$.

If we consider a rigid morphism $f$ in $\mathcal C$ as a two-term complex over $\mathcal C$, then
$f$ is positively self-orthogonal in $\mathscr{K}(\mathcal{C})$, that is, $\Hom_{\mathscr{K}(\mathcal{C})}(f,f[n])=0$ for all $n> 0$.

\medskip
Clearly, a morphism $f:Y\ra X$ in
$\mathcal{C}$ is rigid if and only if $\Hom_{\mathcal
C}(Y,X)=\End_{\mathcal C}(Y)f+f\End_{\mathcal C}(X)$. Thus, the
zero map $Y\ra X$ is rigid if and only if $\Hom_{\mathcal
C}(Y,X)=0$, and any isomorphism $Y\ra X$ is always rigid.

Let us give some non-trivial examples of rigid morphisms, which show
that rigid morphisms exist in very general circumstances.

(i) For an additive category $\mathcal C$, if $f: Y\ra X$ is a
morphism in $\mathcal C$ such that the induced map $\Hom_\mathcal{C}(Y,f):\Hom_{\mathcal
C}(Y,Y)\ra \Hom_{\mathcal C}(Y,X)$ (respectively, $\Hom_\mathcal{C}(f,X):\Hom_{\mathcal C}(X,X)\ra
\Hom_{\mathcal C}(Y,X)$) is surjective, then  $\Hom_{\mathcal
C}(Y,X)=\End_{\mathcal C}(Y)f$ (respectively, $\Hom_{\mathcal
C}(Y,X)=f\End_{\mathcal C}(X)$), and therefore $f$ is rigid. Thus
all approximations in the sense of Auslander-Smalo (see \cite{as}) are
rigid morphisms.

This type of rigid morphisms includes the following three cases:

(a) Let $A$ be an Artin algebra, and let $0\ra Z\lraf{f} Y\lraf{g}
X\ra 0$ be an almost split sequence in $A$-mod. Then both $f$ and
$g$ are rigid since both $\Hom_A(Y,g)$ and $\Hom_A(f,Y)$ are
surjective. For the definition of almost split sequences, we refer
the reader to \cite{ars}.

(b) The covariant morphisms defined in \cite{xc4} are rigid. Recall
that a morphism $f: Y\ra X$ in an additive category $\mathcal C$ is
called covariant if the induced map $\Hom_{\mathcal C}(X,f):
\Hom_{\mathcal C}(X,Y)\ra \Hom_{\mathcal C}(X,X)$ is injective and
the induced map $\Hom_{\mathcal C}(Y,f): \Hom_{\mathcal C}(Y,Y)\ra
\Hom_{\mathcal C}(Y,X)$ is a split epimorphism of $\End_{\mathcal
C}(Y)$-modules.

(c) Let $S$ be a ring with identity. If $Y$ is a quasi-projective
$S$-module (that is, for any surjective homomorphism $Y\ra X$, the
induced map $\Hom_S(Y,Y)\ra \Hom_S(Y,X)$ is surjective), then we may
take a submodule $Z$ of $Y$ and consider the canonical map $f: Y\ra
X:=Y/Z$. Clearly, we have $\Hom_S(Y,X)=\End_S(Y)f$, and therefore
$f$ is rigid. Dually, if $X$ is a quasi-injective $S$-module, that
is, for any injective homomorphism $g: Y\ra X$, the induced map
$\Hom_S(X,X)\ra \Hom_S(Y,X)$ is surjective), then, for any submodule
$Y$ of $X$, we have $\Hom_S(Y,X)=\mu\End(X)$, where $\mu$ is the
inclusion of $Y$ into $X$. This means that $\mu$ is rigid. In
particular, every surjective homomorphism from a projective module
to a module is rigid, and every injective homomorphism from a module
to an injective module is rigid.

(ii) Let $R\subseteq S$ be an extension of rings, that is, $R$ is a
subring of the ring $S$ with the same identity. Then the canonical
map $\pi: S\ra S/R$ of $R$-modules is rigid.

In fact, for any $f\in\Hom_R(S, S/R)$, we choose an element
$s\in S$ such that $(s)\pi=(1)f$, and denote by $\cdot s: S\to S$
the right multiplication by $s$ map. Then the map $f-(\cdot s)\pi$ sends $1\in S$ to zero.
Thus there exists a unique homomorphism $g\in\End_R(S/R)$ such that $f=(\cdot s)\pi+ \pi\,g$.
This implies that $$\Hom_R(S,S/R)=\End_S(S)\pi +\pi \End_R(S/R).$$
Since $\End_S(S)\subseteq \End_R(S)$, we have $\Hom_R(S,S/R)=\End_R(S)\pi +\pi \End_R(S/R)$.
Thus the map $\pi$ is rigid.

We should observe that not every nonzero homomorphism is rigid. For
example, the right multiplication by $x$ map: $k[X]/(X^2)\ra
k[X]/(X^2)$ is not rigid in $k[X]/(X^2)$-Mod, where $k$ is a field
and  $x:=X+(X^2)$ is the coset of $X$ in $k[X]$. In general, an element $x$ in the radical of an Artin algebra $A$,
considered as the right multiplication by $x$ map from ${}_AA$ to itself, is never rigid.

\medskip
Motivated by the rigid morphisms, we introduce the notion of the
so-called hypercyclic bimodules.

Let $S$ and $T$ be two rings with identity, and let $M$ be an
$S$-$T$-bimodule. An element $m\in M$ is called a \emph{hypergenerator}
if $M=Sm+mT$. In this case, $M$ is said to be \emph{hypercyclic}.

Hypercyclic bimodules and rigid morphisms are intimately related in
the following way: If $f:Y\ra X$ is a rigid morphism in an additive
category $\mathcal{C}$, then $f$ is a hypergenerator of the
$\End_{\mathcal{C}}(Y)$-$\End_{\mathcal C }(X)$-bimodule
$\Hom_{\mathcal{C}}(Y,X)$,

If $M$ is hypercyclic with $m$ a hypergenerator, then we may define
a map
$$\zeta: S\oplus T\lra M \quad (s,t)\mapsto sm-mt \mbox{\;  for \; } s\in S \mbox{  and  } t\in T, $$ and get an
exact sequence of ableian groups
$$0\lra K\lra S\oplus T\lraf{\zeta} M\lra 0, $$
where $K:=\{(s,t)\in S\oplus T\mid sm=mt\} $ is a subring of the
ring $S\oplus T$.
Let $p$ and $q$ be the canonical projections from $K$ to $S$ and
$T$, respectively.  Then $S$ and $T$ can be considered as
$K$-$K$-bimodules, and therefore the above sequence is actually an
exact sequence of $K$-$K$-bimodules.

Thus, for each rigid morphism $f: Y\ra X$, there is an exact
sequence

$$0\lra R\lraf{(p,q)} \End_{\mathcal{C}}(Y)\oplus\End_{\mathcal C }(X)
\lraf{\left({\cdot\,f\,\atop{-f\,\cdot}}\right)} \Hom_{\mathcal
C}(Y,X)\lra 0$$ of $R$-$R$-bimodules, where $ R:=\{(s,t)\in
\End_\mathcal{C}(Y)\oplus \End_\mathcal{C}(X)\mid sf=ft\} $ is a
subring of the ring $\End_{\mathcal{C}}(Y)\oplus\End_{\mathcal C
}(X)$.

\smallskip
Now, we give the definition of exact contexts.

\begin{Def}\label{exact-context}
Let $R, S$ and $T$ be rings with identity, let $\lambda: R\to S$ and
$\mu: R\to T$ be ring homomorphisms, and let $M$ be an
$S$-$T$-bimodule with $m\in M$. The quadruple $(\lambda, \mu, M, m)$
is called an \emph{exact context} if
$$(\ast)\quad\;
0\lra R\lraf{(\lambda,\,\mu)}S\oplus T
\lraf{\left({\cdot\,m\,\atop{-m\,\cdot}}\right)}M\lra 0
$$
is an exact sequence of abelian groups, where $\cdot m$ and $m \cdot$ stand for the right and left multiplication by $m$ maps, respectively.
In this case, we also say that$(M,m)$ is an exact complement of $(\lambda,\mu)$.

If $(\lambda, \mu, S\otimes_RT, 1\otimes 1)$ is an exact context, then we
say simply that $(\lambda,\mu)$ is an \emph{exact pair}.
\end{Def}

Note that the sequence $(\ast)$ is exact in the category of abelian groups if and only if

$(E_1)$ the $S$-$T$-bimodule $M$ is hypercyclic with $m$ as a hypergenerator, and

$(E_2)$ the ring homomorphism $R\lraf{(\lambda,\,\mu)}S\oplus T$ induces a ring isomorphism from $R$ to $K$.

\smallskip
Note that the quadruple $(\lambda, \mu, M, m)$ is an exact context if and only if the following
diagram $$(\sharp)\quad
\xymatrix{
R\ar[d]_-{\mu} \ar[r]^-{\lambda} & S \ar[d]^-{\cdot m}\\
T \ar[r]^-{m \cdot} &  M }
$$
is both a push-out and a pull-back in the category of $R$-$R$-bimodules.

Let $(\lambda,\mu, M, m)$ be an exact context. Then, from $(\sharp)$ we see
that, for an $S$-$T$-bimodule $N$ with an element $n\in N$, the pair
$(N,n)$ is an exact complement of $(\lambda, \mu)$ if and only if
there exists a unique isomorphism $\omega:M\to N$ of
$R$-$R$-bimodules such that $(sm)\omega=sn$ and $(mt)\omega=nt$ for
all $s\in S$ and $t\in T$. Clearly, $\omega$ preserves
hypergenerators, that is $(m)\omega=n$. In general, $\omega$ has not
to be an isomorphism of $S$-$T$-bimodules, that is, $M$ and $N$ may
not be isomorphic as $S$-$T$-bimodules (see the examples in Subsection \ref{sect3.2.1}).

\smallskip
Next, we mention several examples of exact contexts.

$(1)$ Let $M$ be a hypercyclic $S$-$T$-bimodule with $m$ a
hypergenerator. Then the pair $(p,q)$ of ring homomorphisms $p: K\to
S$ and $q: K\to T$  together with $(M, m)$ forms an exact context.
So rigid morphisms always provide us with a class of exact contexts.
Conversely, every exact context appears in this form. In fact, for a
given exact context $(\lambda, \mu, M,m)$, we may define $B =
\left(\begin{array}{lc} S& M\\ 0 & T\end{array}\right)$ and consider
the canonical map $\varphi$ from the first column to the second
column of $B$ defined by $\cdot m$. It is easy to see that this
$\varphi$ is rigid and the induced exact context is precisely the
given one. So, rigid morphisms describe exact contexts.

$(2)$ Suppose that $R\subseteq S$ is an extension of rings. Let $\lambda:R\to S$ be the inclusion   with $\pi:S\to S/R$ the canonical surjection. Define $S':=\End_R(S/R)$ and
$$\lambda':\,R \lra S':\;\;r\mapsto (x\mapsto xr)\;\,\mbox{for}\;r\in R\;\mbox{and}\; x\in S/R.$$
Then $\Hom_R(S,S/R)$ is an $S$-$S'$-bimodule, and the quadruple
$\big(\lambda, \lambda', \Hom_R(S, S/R),\, \pi \big)$ is an exact
context since the following diagram
$$
\xymatrix{ 0\ar[r] &R \ar[d]_-{\lambda'} \ar[r]^-{\lambda}
                & S \ar[d]_{\cdot\,\pi}\ar[r]^-{\pi}
                & S/R\ar[d]_-{\simeq}\ar[r]& 0\\
0\ar[r] & S'\ar[r]^-{\pi\,\cdot} & \Hom_R(S,
S/R)\ar[r]^-{\lambda\,\cdot}& \;\Hom_R(R,S/R)\ar[r]& 0 }
$$
is commutative and the sequence of $R$-$R$-bimodules
$$
\xymatrix{0\ar[r] & R\ar[r]^-{(\lambda,\,\lambda')\,}&  S\oplus S'
\ar[r]^-{\left({\cdot \pi\,\atop{-\pi\cdot}}\right)}&
\Hom_R(S,S/R)\ar[r]& 0.}
$$
is exact. In general, the exact context presented here is different
from the one induced  from the rigid morphism $\pi$, and the pair
$(\lambda, \lambda')$ may not be exact, because either $S\simeq
\End_R(S)$ as rings or $S\otimes_RS'\simeq \Hom_R(S,S/R)$ as
$S$-$S'$-bimodules may fail.

A more general construction of exact contexts from a (not
necessarily injective) ring homomorphism will be discussed in Lemma
\ref{general case}.

$(3)$ Milnor squares, defined by Milnor in \cite[Sections 2 and
3]{Milnor}, also provide a class of  exact contexts. Recall that a
Milnor square is a commutative diagram of ring
homomorphisms
$$
\xymatrix{
\Lambda \ar[d]_-{i_2} \ar[r]^-{i_1} & \Lambda_1 \ar[d]^-{j_1}\\
\Lambda_2\ar[r]^-{j_2} &  \Lambda'}
$$
satisfying the following two conditions:

$(M1)$ The ring $\Lambda$ is the pull-back of $\Lambda_1$ and $\Lambda_2$ over $\Lambda'$, that is, given
a pair $(\lambda_1, \lambda_2)\in\Lambda_1\oplus \Lambda_2$ with $(\lambda_1)j_1=(\lambda_2)j_2\in\Lambda'$,
there is one and only one element $\lambda\in\Lambda$ such that $(\lambda)i_1=\lambda_1$ and
$(\lambda)i_2=\lambda_2$.

$(M2)$ At least one of the two homomorphisms $j_1$ and $j_2$ is surjective.

\smallskip
Clearly, $\Lambda'$ can be regarded as an $\Lambda_1$-$\Lambda_2$-bimodule via the ring homomorphisms $j_1$ and $j_2$.
Let $1$ be the identity of $\Lambda'$. Then $j_1$  and $j_2$ are exactly the multiplication maps $\cdot\,1$ and $1\,\cdot$, respectively.

Now, we claim that the pair $(i_1, i_2)$  together with $(\Lambda', 1)$ forms an exact context.
Indeed, it follows from the condition $(M2)$ that $\Lambda'$ is hypercyclic with $1$ as a hypergenerator.
With the help of the condition $(M1)$, the following sequence
$$
0\lra \Lambda\lraf{(i_1,\,i_2)}\Lambda_1\oplus \Lambda_2
\lraf{\left({j_1\,\atop{-j_2}}\right)}\Lambda'\lra 0
$$
is an exact sequence of $\Lambda$-$\Lambda$-bimodules. This verifies the claim.

Even more, the pair $(i_1, i_2)$ is exact. Without loss of generality, assume that
$j_2$ is surjective. Then, by $(M1)$, the map $i_1$ is also surjective and $i_2$
induces an isomorphism $\Ker(i_1)\simeq \Ker(j_2)$ of $\Lambda$-$\Lambda$-bimodules.
Now, we can check that the map $\Lambda_1\otimes_\Lambda\Lambda_2\to \Lambda'$, defined by
$\lambda_1\otimes\lambda_2\mapsto (\lambda_1)j_1 (\lambda_2)j_2$ for $\lambda_1\in\Lambda_1$ and
$\lambda_2\in\Lambda_2$, is an isomorphism of
$\Lambda_1$-$\Lambda_2$-bimodules. Actually, this follows from the
following isomorphisms:
$$\Lambda_1\otimes_{\Lambda}\Lambda_2\simeq\big(\Lambda/\Ker(i_1)\big)\otimes_\Lambda\Lambda_2
\simeq\Lambda_2/\big(\Ker(i_1)\Lambda_2\big)=\Lambda_2/(\Ker(j_2)\Lambda_2)\simeq\Lambda'.$$

Similarly, we can check that the pair $(i_2, i_1)$ is also exact
with $\Lambda_2\otimes_\Lambda\Lambda_1\simeq \Lambda'$ as $\Lambda_2$-$\Lambda_1$-bimodules.

\section{Noncommutative tensor products of exact contexts\label{sect3}}

In this section, we shall define a new ring for each exact context.
This is the so-called noncommutative tensor product which includes
the notion of coproducts of rings, usual tensor products and so on.
These noncommutative tensor products can be constructed from both
Morita context rings and strictly pure extensions, and will play a
crucial role in construction of recollements of derived module
categories in the next section.

\subsection{Definition of noncommutative tensor products\label{sect3.1}}
From now on, let $\lambda:R\to S$ and $\mu:R\to T$ be two arbitrary but fixed ring homomorphisms.
Unless stated otherwise, we always assume that $(\lambda, \mu, M, m)$ is  an exact context.

First, we characterize when the pair $(\lambda, \mu)$ in the exact context is exact.
Recall that we have the following exact sequence of $R$-$R$-bimodules:
$$(*)\quad
\xymatrix{0\ar[r] & R\ar[r]^-{(\lambda,\,\mu)\,}& S\oplus
\;T\ar[r]^-{\left({\cdot m\,\atop{-m \cdot}}\right)}&
M\ar[r]& 0.}
$$
According to $(\ast)$, there exist two unique homomorphisms
$$
\alpha:M\lra S\otimes_RT,\; x\mapsto s_x\otimes 1+1\otimes t_x\quad \mbox{and}\quad
\beta: M\lra T\otimes_RS,\; x\mapsto 1\otimes s_x+t_x\otimes 1,
$$
where  $x\in M$ and $(s_x, t_x)\in S\oplus T$  with $x=s_xm+mt_x$, such that the following two diagrams
$$(\dag)\quad
\xymatrix{R\ar@{=}[d]\ar[r]^-{(\lambda,\,\mu)\,} &S\oplus T \ar@{=}[d]\ar[r]^-{\left({\cdot m\,\atop{-m \cdot}}\right)}
                & M \ar[d]^-{\alpha}\\
 R\ar[r]^-{(\lambda,\,\mu)\,} & S\oplus T\ar[r]^-{\left({\mu\,'\,\atop{-\lambda\,'}}\right)}& S\otimes_RT}
\quad \mbox{and}\;\quad(\ddag)\quad
\xymatrix{R\ar@{=}[d]\ar[r]^-{(\lambda,\,\mu)\,} &S\oplus T \ar@{=}[d]\ar[r]^-{\left({\cdot m\,\atop{-m \cdot}}\right)}
                & M \ar[d]^-{\beta}\\
 R\ar[r]^-{(\lambda,\,\mu)\,} & S\oplus T\ar[r]^-{\left({\rho\,\atop{-\phi}}\right)}& T\otimes_RS.}$$
are commutative, where
$$\lambda'=\lambda\otimes_RT: \; T\lra S\otimes_RT, \; t\mapsto 1\otimes t \quad \mbox{and}\quad \mu'=S\otimes_R\mu: \; S\lra S\otimes_RT, \; s\mapsto s\otimes 1, $$
$$\rho=\mu\otimes S:\; S\lra T\otimes_RS,\; s\mapsto 1\otimes s \quad \mbox{and}\quad \phi=T\otimes\lambda:\; T\lra T\otimes_RS,\; t\mapsto t\otimes 1$$
for $s\in S$ and $t\in T$. Note that $(x)\alpha$ and $(x)\beta$ are independent of different choices of
$(s_x,t_x)$ in $S\oplus T$.

Further, let
$$\gamma:S\otimes_RT\lra M,\; s\otimes t\mapsto smt.$$
Clearly, $\alpha$ and $\beta$ are homomorphisms of $R$-$R$-bimodules,
$\gamma$ is a homomorphism of $S$-$T$-bimodules and $\alpha\,\gamma=Id_M$. In particular, $\alpha$ is injective and $\gamma$ is surjective.

\begin{Lem}\label{exact pair}
The following statements are equivalent:

$(1)$ The pair $(\lambda, \mu)$ is an exact pair.

$(2)$ The map $\gamma$ is an isomorphism.

$(3)$ $\Coker(\lambda)\otimes_R\Coker(\mu)=0$.

$(4)$ $(M/mT)\otimes_R(M/Sm)=0$.
\end{Lem}

{\it Proof.} Note that $\gamma$ is a homomorphism of $S$-$T$-bimodules and $$\big(s(1\otimes 1)\big)\gamma=(s\otimes 1)\gamma=sm\quad\mbox{and}\quad \big((1\otimes 1)t\big)\gamma=(1\otimes t)\gamma=mt$$ for $s\in S$ and $t\in T$.
This implies that the following diagram
$$
\xymatrix{& R\ar@{=}[d]\ar[r]^-{(\lambda,\,\mu)\,} &S\oplus T\ar[r]^-{\left({\mu\,'\atop{-\lambda\,'}}\right)} \ar@{=}[d] &S\otimes_RT\ar[d]^-{\gamma}&\\
0\ar[r] &R\ar[r]^-{(\lambda,\,\mu)\,} & S\oplus T\ar[r]^-{\left({\cdot m\,\atop{-m \cdot}}\right)}& M\ar[r]& 0}
$$
is commutative, where the second row is assumed to be exact. Consequently, $(1)$ and $(2)$ are equivalent.

According to $(\sharp)$, we know that $\Coker(\lambda)\simeq M/mT$
and $\Coker(\mu)\simeq M/Sm$ as $R$-$R$-bimodules. Thus $(3)$ and $(4)$ are equivalent.

Now, we verify the equivalences of $(2)$ and $(3)$.

In fact, since $\alpha\,\gamma=Id_M$, the map $\gamma$ is an
isomorphism if and only if $\alpha$ is a surjection, while the
latter is equivalent to that the map
$$\xi:=\left({\mu\,'\atop{-\lambda\,'}}\right):S\oplus T\lra S\otimes_RT$$ is a surjection
by $(\dag)$. Therefore, it is enough to show that $\xi$ is
surjective if and only if $\Coker(\lambda)\otimes_R\Coker(\mu)=0$.
To check this condition, we consider the following two complexes
$$\Cone(\lambda): 0\to R\lraf{\lambda} S\to 0\quad \mbox{and}\quad \Cone(\mu): 0\to R\lraf{\mu} T\to 0$$
of $R$-$R$-bimodules, where both $S$ and $T$ are of degree $0$, and calculate the tensor complex of them over $R$:
$$\Cone(\lambda)\cpx{\otimes}_R\Cone(\mu):
\xymatrix{0\ar[r] & R\otimes_RR\ar[rr]^-{\big(\lambda\otimes R,\; -R\otimes\mu\big)}&&S\otimes_RR\oplus R\otimes_RT
\ar[rr]^-{\left({S\otimes\mu\;\atop{\lambda\otimes T }}\right)}&&
S\otimes_RT\ar[r]& 0}$$
where $R\otimes_RR$ is of degree $-2$. If we identify $R\otimes_RR$, $S\otimes_RR$ and $R\otimes_RT$
with $R$, $S$ and $T$, respectively, then $\Cone(\lambda)\cpx{\otimes}_R\Cone(\mu)$
is precisely the complex:
$$\xymatrix{0\ar[r] & R\ar[r]^-{\big(\lambda,\; -\mu\big)}& S\oplus
T\ar[r]^-{\left({\mu\,'\;\atop{\lambda\,'}}\right)}&
S\otimes_RT\ar[r]& 0}$$
which is isomorphic to the following complex
$$\xymatrix{0\ar[r] & R\ar[r]^-{\big(\lambda,\; \mu\big)}& S\oplus
T\ar[r]^-{\xi}&
S\otimes_RT\ar[r]& 0.}$$
It follows that $\xi$ is surjective if and only if $H^0\big(\Cone(\lambda)\cpx{\otimes}_R\Cone(\mu)\big)=0$.
Since $$H^0\big(\Cone(\lambda)\cpx{\otimes}_R\Cone(\mu)\big)\simeq H^0(\Cone(\lambda))\otimes_R H^0(\Cone(\mu))
\simeq \Coker(\lambda)\otimes_R\Coker(\mu),$$
the map $\xi$ is surjective
if and only if $\Coker(\lambda)\otimes_R\Coker(\mu)=0$. Thus
$\gamma$ is an isomorphism if and only if $\Coker(\lambda)\otimes_R\Coker(\mu)=0$.
This shows the equivalences of $(2)$ and $(3)$.  $\square$

\begin{Rem}
By the equivalences of $(1)$ and $(2)$ in Lemma \ref{exact pair}, if
the  pair $(\lambda, \mu)$ is exact, then it admits a unique
complement $(S\otimes_RT, 1\otimes 1)$ up to isomorphism (preserving
hypergenerators) of $S$-$T$-bimodules.
\end{Rem}

A sufficient condition to guarantee the isomorphism of $\gamma$ is
the following result.

\begin{Koro}\label{app1}
If either $\lambda:R\to S$ or $\mu:R\to T$ is a ring epimorphism, then
$\gamma:S\otimes_RT\to M, s\otimes t\mapsto smt$ is an isomorphism of $S$-$T$-bimodules.
\end{Koro}

{\it Proof.} Suppose that $\lambda$ is a ring epimorphism. Then, for any
$S$-module $X$, the map $\lambda\otimes X:R\otimes_RX\to S\otimes_R X$ is an
isomorphism. This implies that $\Coker(\lambda)\otimes_RX=0$. Since $\Coker(\mu)
\simeq M/Sm$ as $R$-modules by $(\sharp)$ and since $M/Sm$ is an $S$-module,
we have $\Coker(\lambda)\otimes_R\Coker(\mu)\simeq \Coker(\lambda)\otimes_R(M/Sm)=0$.
By Lemma \ref{exact pair}, the map $\gamma$ is an isomorphism.

Similarly, if $\mu$ is a ring epimorphism, then $\gamma$ is an isomorphism. $\square$

\medskip
As examples of exact pairs, we see from Corollary \ref{app1} that the rigid morphisms from an almost split sequence always provide us with exact pairs.

\medskip
Next, we shall introduce the so-called noncommutative tensor products $T\boxtimes_RS$ of the exact context $(\lambda, \mu, M, m)$. That is, we endow $T\otimes_RS$ with an associative multiplication
$\circ: (T\otimes_RS)\times (T\otimes_RS)\to T\otimes_RS$, under which it becomes an associative ring with the identity $1\otimes 1$.

Let
$$\delta:=\gamma\,\beta:\;S\otimes_RT\lra T\otimes_RS,\; s\otimes t\mapsto 1\otimes s_{smt}+ t_{smt}\otimes 1$$
for $s\in S$ and $t\in T$, where the pair $(s_{smt}, t_{smt})\in S\oplus T$ is chosen such that $smt=s_{smt}m+m t_{smt}$. Then $\delta$ is a homomorphism of $R$-$R$-bimodules such that $(s\otimes 1)\delta=1\otimes s$ and $(1\otimes t)\delta=t\otimes 1$.

The multiplication $\circ$ is induced by the following homomorphisms:
{\small $$
(T\otimes_RS)\otimes_R(T\otimes_RS)\lraf{\simeq}T\otimes_R(S\otimes_RT)\otimes_RS\lraf{\,T\otimes \delta\otimes S\,}
T\otimes_R(T\otimes_RS)\otimes_R S\lraf{\simeq} (T\otimes_RT)\otimes_R(S\otimes_R S)\lraf{\,\mu_T\otimes \mu_S\,} T\otimes_RS
$$}\noindent where $\mu_T: T\otimes_RT\to T$ and $\mu_S:S\otimes_RS\to
S$ are the multiplication maps. More precisely, for $(t_i, s_i)\in
T\otimes_RS$ with $i=1,2$, we have
$$(t_1\otimes s_1)\circ (t_2\otimes s_2):=t_1(s_1\otimes t_2)\delta\,s_2
=t_1\big( 1\otimes s_{s_1mt_2}+ t_{s_1mt_2}\otimes 1\big)s_2.$$

\medskip
The following lemma reveals a crucial property of this multiplication.

\begin{Lem}\label{ring}
The following statements are true.

$(1)$ With the multiplication
$\circ$, the abelian group $T\otimes_RS$ becomes an associative ring with the identity $1\otimes 1$.

$(2)$ The maps $\rho: S\to T\otimes_RS$ and $\phi: T\to T\otimes_RS$ are ring homomorphisms. In particular,
$T\otimes_RS$ can be regarded as an $S$-$T$-bimodule via $\rho$ and $\phi$.

$(3)$ The map $\beta:M\to T\otimes_RS$ is a homomorphism of $S$-$T$-bimodules such that $(m)\beta=1\otimes 1$.
\end{Lem}

{\it Proof.} $(1)$ It suffices to show that the multiplication $\circ$ is associative and that $1\otimes 1$
is the identity of $T\otimes_RS$.

To check the associativity of $\circ$, we take elements $t_i\in T$ and $s_i\in S$ for $1\leq i\leq 3$, and choose two pairs $(x,y)$ and $(u,v)$ in $S\times T$ such that
$$s_1mt_2=xm+my\quad \mbox{and}\quad s_2mt_3=um+mv.$$

On the one hand,
$\big((t_1\otimes s_1)\circ (t_2\otimes s_2)\big)\circ(t_3\otimes s_3)=\big(t_1(s_1\otimes t_2)\delta\,s_2\big)\circ (s_3\otimes t_3)=\big(t_1\big(1\otimes x+ y\otimes 1\big)s_2\big)\circ (t_3\otimes s_3)=(t_1\otimes xs_2)\circ(t_3\otimes s_3)+ (t_1y\otimes s_2)\circ(t_3\otimes s_3)=t_1\big((xs_2\otimes t_3)\delta +y(s_2\otimes t_3)\delta\big)s_3.$

\smallskip
On the other hand, $(t_1\otimes s_1)\circ\big((t_2\otimes s_2)\circ(t_3\otimes s_3)\big)=(t_1\otimes s_1)\circ \big(t_2(s_2\otimes t_3)\delta\,s_3\big)=(t_1\otimes s_1)\circ \big(t_2(1\otimes u+ v\otimes 1)s_3\big)=(t_1\otimes s_1)\circ(t_2\otimes us_3)+(t_1\otimes s_1)\circ(t_2v\otimes s_3)=t_1\big((s_1\otimes t_2)\delta\,u +(s_1\otimes t_2v)\delta\big)s_3.$

So, to prove that $$\big((t_1\otimes s_1)\circ (t_2\otimes s_2)\big)\circ(t_3\otimes s_3)=(t_1\otimes s_1)\circ\big((t_2\otimes s_2)\circ(t_3\otimes s_3)\big),$$ it is enough to verify that $$(xs_2\otimes t_3)\delta +y(s_2\otimes t_3)\delta=(s_1\otimes t_2)\delta\,u +(s_1\otimes t_2v)\delta.$$
In fact, since $xs_2mt_3=x(um+mv)=xum+xmv$ and $xu\in S$, we have $(xs_2\otimes t_3)\delta=(xum+xmv)\beta=1\otimes xu+(xmv)\beta$. Similarly, $(s_1\otimes t_2v)\delta=yv\otimes 1+(xmv)\beta$.
It follows that
$$(xs_2\otimes t_3)\delta +y(s_2\otimes t_3)\delta=1\otimes xu+(xmv)\beta+y(1\otimes u+v\otimes 1)
=1\otimes xu+(xmv)\beta+y\otimes u+yv\otimes 1$$
$$
\quad\;=1\otimes xu+ y\otimes u+yv\otimes 1+(xmv)\beta=(1\otimes x+ y\otimes 1)u+yv\otimes 1+(xmv)\beta=(s_1\otimes t_2)\delta\,u +(s_1\otimes t_2v)\delta.
$$
This shows that the multiplication $\circ$ is associative.

Note that $(t_1\otimes s_1)\circ (1\otimes 1)=t_1(s_1\otimes 1)\delta=t_1(1\otimes s_1)=t_1\otimes s_1$ and
$(1\otimes 1)\circ (t_1\otimes s_1)=(1\otimes t_1)\delta\,s_1=(t_1\otimes 1)s_1=t_1\otimes s_1$. Thus
$(T\otimes_RS, \circ)$ is an associative ring with the identity $1\otimes 1$.

$(2)$ Since $(s_1)\rho \circ (s_2)\rho=(1\otimes s_1)\circ (1\otimes s_2)=(s_1\otimes 1)\delta\,s_2=(1\otimes s_1)s_2=1\otimes s_1s_2=(s_1s_2)\rho$, the map $\rho:S\to T\otimes_RS$ is a ring homomorphism. Similarly, we can show that $\phi:T\to T\otimes_RS$ is also a ring homomorphism.

$(3)$ Clearly, by the definition of $\beta$, we have
$(m)\beta=1\otimes 1$. It remains to check that $\beta$ is a
homomorphism of $S$-$T$-bimodules, or equivalently, that
$$(sat)\beta=(s)\rho \circ (a)\beta \circ(t)\phi$$ for $s\in S$,  $a\in M$ and $t\in T$.

To check this, we pick up $s_a\in S$ and $t_a\in T$ such that $a=s_am+mt_a$. Then $(sat)\beta=(ss_amt+smt_at)\beta
=(ss_amt)\beta+(smt_at)\beta=(ss_a\otimes t)\delta+(s\otimes t_at)\delta=(1\otimes ss_a)\circ (t\otimes 1)+
(1\otimes s)\circ (t_at\otimes 1)=(1\otimes s)\circ (1\otimes s_a)\circ (t\otimes 1)+(1\otimes s)\circ
(t_a\otimes 1)\circ (t\otimes 1)=(1\otimes s)\circ (1\otimes s_a+t_a\otimes 1)\circ (t\otimes 1)
=(s)\rho \circ (a)\beta \circ(t)\phi$. $\square$

\medskip
Thanks to Lemma \ref{ring}, the ring$(T\otimes_RS,\circ)$ will be
called the \textbf{noncommutative tensor product} of the exact context $(\lambda,\mu,
M, m)$, denoted simply by $T\boxtimes_RS$ if the exact context
$(\lambda,\mu, M,m)$ is clear.

We should note that the ring $T\boxtimes_RS$ is not the usual tensor
product of two $R$-algebras: First, the ring $R$ is not necessarily
commutative, this means that the usual tensor product of
$R$-algebras on the abelian group $T\otimes_RS$ does not make sense.
Second, even if the ring $R$ is commutative, we cannot ensure that
the product has to coincide with the usual tensor product because
the image of $\lambda: R\ra S$ does not have to be in the center of
$S$. This means that $S$ is not necessarily an $R$-algebra.
Nevertheless, the ring $T\boxtimes_RS$ does generalize the usual
tensor product of $R$-algebras in the following sense:

Let $R$ be a commutative ring. Suppose that $S$ and $T$ are $R$-algebras via $\lambda$ and $\mu$, respectively, that is,
the images of $\lambda$ and $\mu$ are contained in the centers of $S$ and $T$, respectively. If $(\lambda,\mu)$ is an exact pair, then the noncommutative tensor product $T\boxtimes_RS$ coincides with the usual tensor product $T\otimes_RS$ of $R$-algebras $T$ and $S$.

In fact, by our notation, we have $M=S\otimes_RT$, $\gamma=Id_{S\otimes_RT}$ and $\delta=\beta:S\otimes_RT\to T\otimes_RS$,
where $\beta$ is determined uniquely by the diagram:
$$ (\ddag)\qquad
\xymatrix{R\ar@{=}[d]\ar[r]^-{(\lambda,\,\mu)\,} &S\oplus T \ar@{=}[d]\ar[r]^-{\left({\mu\,'\atop{-\lambda\,'}}\right)}
                & S\otimes_RT \ar[d]^-{\beta}\\
R\ar[r]^-{(\lambda,\,\mu)\,} & S\oplus T\ar[r]^-{\left({\rho\,\atop{-\phi}}\right)}& T\otimes_RS}$$
However, since $(\lambda,\mu)$ is exact, we can check that the switch map $\omega: S\otimes_RT\to T\otimes_RS$, defined by $s\otimes t\mapsto t\otimes s$ for $s\in S$ and $t\in T$,
also makes the above diagram commutative, that is, $\left({\mu\,'\atop{-\lambda\,'}}\right) \omega=\left({\rho\,\atop{-\phi}}\right)$. This implies
that $\beta=\omega$. Thus the multiplication $\circ:(T\boxtimes_RS)\times (T\boxtimes_RS)\to T\boxtimes_RS$ coincides with the usual tensor product of $R$-algebras $T$ and $S$ over $R$.

\subsection{Examples of noncommutative tensor products \label{sect4.2}}

In this section, we present two general receipts for constructing noncommutative tensor products, which show that noncommutative tensor products cover a large variety of interesting algebras.

\subsubsection{From Morita context rings \label{sect3.2.1}}

Let $(A,C,X,Y,f,g)$ be an arbitrary but fixed Morita context, that is, $A$ and $C$ are rings with identity, $X$ is an $A$-$C$-bimodule, $Y$ is a $C$-$A$-bimodule, $f:X\otimes_CY\to A$ is a homomorphism of $A$-$A$-bimodules and $g:Y\otimes_AX\to C$ is a homomorphism of $C$-$C$-bimodules, such that
$$(x_1\otimes y_1)f x_2=x_1(y_1\otimes x_2)g\;\;\mbox{and}\;\;(y_1\otimes x_1)g y_2=y_1(x_1\otimes y_2)f$$ for $x_i\in X$ and $y_i\in Y$ with $i=1,2$. For simplicity, we denote by $x_1y_1$ and $y_1x_1$ the elements $(x_1\otimes y_1)f$ and $(y_1\otimes x_1)g$, respectively.

Given a Morita context $(A,C,X,Y,f,g)$, we can define the Morita context ring $\Gamma:=\left(\begin{array}{lc} A &X \\
Y & C\end{array}\right)$, where the multiplication is given by
$$\left(\begin{array}{lc} a_1 &x_1 \\
y_1& c_1\end{array}\right)\left(\begin{array}{lc} a_2 &x_2 \\
y_2& c_2\end{array}\right)=\left(\begin{array}{lc} a_1a_2+x_1y_2 &a_1x_2+x_1c_2 \\
y_1a_2+c_1y_2& c_1c_2+y_1x_2\end{array}\right)$$
for $a_i\in A$, $c_i\in C$, $x_i\in X$ and $y_i\in Y$.

Let $$R:=\left(\begin{array}{lc} A &0 \\
0& C\end{array}\right),\; S:=\left(\begin{array}{lc} A &X \\
0& C\end{array}\right),\; T:=\left(\begin{array}{lc} A &0 \\
Y& C\end{array}\right),\; M:=\left(\begin{array}{lc} A &X \\
Y & C\end{array}\right),\; m:=\left(\begin{array}{lc} 1 &0 \\
0& 1\end{array}\right)\in M,$$
and let $\lambda:R\to S$ and $\mu:R\to T$ be the canonical inclusions. Note that the $S$-$T$-bimodule structure on $M$ is induced from the ring structure of the Morita context ring $\Gamma$.
Since $R=S\cap T$ and $M =S + T$, the quadruple $(\lambda, \mu, M, m)$ is an exact context. So we can consider the noncommutative tensor product $T\boxtimes_RS$ of this exact context. In fact, the multiplication in $T\boxtimes_RS$ can be described explicitly as follows:

We identify $R\Modcat$ with the product $A\Modcat \times C\Modcat$. In this sense, $_RS=(A\oplus X)\times C$ and $T_R=(A\oplus Y)\times C$. It follows that the following homomorphism $$T\otimes_RS\lra \left(\begin{array}{lc} A &X \\
Y & C\oplus (Y\otimes_AX)\end{array}\right)=:\Lambda,$$ defined by
$$\left(\begin{array}{lc} a_1 &0 \\
y_1& c_1\end{array}\right)\otimes\left(\begin{array}{lc} a_2 &x_2 \\
0& c_2\end{array}\right)\mapsto \left(\begin{array}{lc} a_1a_2 &a_1x_2 \\
y_1a_2& (c_1c_2,y_1\otimes x_2)\end{array}\right),$$ is an
isomorphism of abelian groups. Via this isomorphism, we identify
$T\otimes_RS$ with $\Lambda$ and translate the multiplication of
$T\boxtimes_RS$ into the one of $\Lambda$. By calculation, this
multiplication on $\Lambda$ is exactly given by the following
formula:
$$\left(\begin{array}{lc} a_1 & x_1\\
y_1& (c_1, y\otimes x)\end{array}\right)\circ \left(\begin{array}{lc} a_2 &x_2 \\
y_2& (c_2, y'\otimes x')\end{array}\right)$$ $$ = \left(\begin{array}{lc} a_1a_2+x_1y_2 &a_1x_2+x_1c_2+x_1(y'x') \\
y_1a_2+c_1y_2+(yx)y_2& \big(c_1c_2,\,y_1\otimes x_2+(c_1y')\otimes x'+y\otimes (xc_2)+ y\otimes (xy')x'\big)\end{array}\right),$$
where $x,x'\in X$ and $y,y'\in Y$. Thus $T\boxtimes_RS=\Lambda$. In this sense, the associated homomorphisms
$\rho:S\to T\boxtimes_RS$, $\phi:T\to T\boxtimes_RS$ and $\beta: M\to T\otimes_RS$ are given by
{\small $$
\left(\begin{array}{lc} a_1 & x_1\\
0& c_1\end{array}\right)\stackrel{\rho}{\mapsto} \left(\begin{array}{lc} a_1 & x_1\\
0& (c_1, 0)\end{array}\right),\;\left(\begin{array}{lc} a_1 & 0\\
y_1& c_1\end{array}\right)\stackrel{\phi}{\mapsto} \left(\begin{array}{lc} a_1 & 0\\
y_1& (c_1, 0)\end{array}\right),\;\left(\begin{array}{lc} a_1 & x_1\\
y_1& c_1\end{array}\right)\stackrel{\beta}{\mapsto} \left(\begin{array}{lc} a_1 & x_1\\
y_1& (c_1, 0)\end{array}\right),
$$}respectively. Note that both $\rho$ and $\phi$ are ring homomorphisms. However,
$\beta$ does not have to be a ring homomorphism in general.
Actually, it is a ring homomorphism if and only if $Y\otimes_A X=0$. Moreover, it follows from the multiplication of $\Lambda$ that the map $$\pi:\Lambda\lra \Gamma,\;\left(\begin{array}{lc} a_1 & x_1\\
y_1& (c_1, y\otimes x)\end{array}\right)\mapsto \left(\begin{array}{lc} a_1 & x_1\\
y_1& c_1+ yx\end{array}\right)$$
is a surjective ring homomorphism such that $\beta\pi=Id_M$. Further, let
$e:=\left(\begin{array}{lc} 1 & 0\\
0& (0,0)\end{array}\right)\in\Lambda$. Then $e^2=e$, $e\Lambda e=A$, $\Lambda e=A\oplus Y$,
$e\Lambda=A\oplus X$, $\Lambda e\Lambda=\left(\begin{array}{lc} A &X \\
Y & Y\otimes_AX\end{array}\right)$ and $\Lambda/(\Lambda e\Lambda)=C$.
This also implies that the canonical multiplication map
$\Lambda e\otimes_A e\Lambda\to \Lambda e\Lambda$ is an isomorphism of $\Lambda$-$\Lambda$-bimodules. So, if $\Tor_i^A(Y,X)=0$ for all $i>0$, then the canonical surjective map $\Lambda\ra \Lambda/\Lambda e\Lambda$ is homological.

For each $i\geq 1$, we have $\Tor_i^R(T,S)\simeq \Tor_i^A(Y,X)$.
Thus $\Tor_i^R(T,S)=0$ if and only if $\Tor_i^A(Y,X)=0$.

\medskip
Let us give some examples to illustrate how the choices of structure
maps in the Morita contexts influence the noncommutative tensor
products of exact contexts in the above construction.

Let $k$ be a field, and let $A=C=X=Y=k$. Now we take two different
kinds of structure maps $f:X\otimes_RY\to A$ and $g:Y\otimes_RX\to
C$ as follows:

(i) Let $f$ and $g$ be the canonical isomorphism
$k\otimes_kk\lraf{\simeq} k$. Then the Morita context ring is the
matrix ring $M_2(k)$ of $2\times 2$ matrices over $k$. In this case,
the noncommutative tensor product of the corresponding exact context
is
$$\Lambda:=T\otimes_RS=\left(\begin{array}{lc} k &k \\
k & k\oplus k\end{array}\right)$$
with the multiplication given by
{\small $$\left(\begin{array}{lc} a_1 & x_1\\
y_1& (c_1, x)\end{array}\right)\circ \left(\begin{array}{lc} a_2 &x_2 \\
y_2& (c_2, x')\end{array}\right)\mapsto \left(\begin{array}{lc} a_1a_2+x_1y_2 &a_1x_2+x_1c_2+x_1x' \\
y_1a_2+c_1y_2+xy_2& \big(c_1c_2,\,y_1x_2+c_1x'+xc_2+ xx'\big)\end{array}\right),$$}
\smallskip
\noindent where $x,x', a_i, c_i, x_i,y_i\in k$ for $i=1,2$. Actually, this ring is Morita equivalent to $k\times k$ since $e_2 := \left(\begin{array}{lc} 0 & 0\\
0& (1, 0)\end{array}\right)$ = $\left(\begin{array}{lc} 0 & 0\\
0& (1, -1)\end{array}\right)+\left(\begin{array}{lc} 0 & 0\\
0& (0, 1)\end{array}\right)=e_2'+e_2''$ and $\Lambda e\simeq \Lambda e_2''$.

\medskip
(ii) Let $f$ and $g$ be the zero homomorphism $k\otimes_kk\to  k$.
Then the Morita context ring, denoted by $M_2(k)_0$, has the vector
space $M_2(k)$ and admits a new multiplication:
{\small $$\left(\begin{array}{lc} a_1 & x_1\\
y_1& c_1\end{array}\right)\left(\begin{array}{lc} a_2 &x_2 \\
y_2& c_2\end{array}\right)= \left(\begin{array}{lc} a_1a_2 &a_1x_2+x_1c_2 \\
y_1a_2+c_1y_2& c_1c_2\end{array}\right).$$}
Note that $M_2(k)_0$ can be identified with the following quiver algebra with
relations
$$\xymatrix{1\,\bullet\ar^-{\bf{\alpha}}@/^0.8pc/[r]
&\bullet\, 2\ar^-{\bf{\beta}}@/^0.8pc/[l]},
\quad\alpha\beta=\beta\alpha=0.$$
In this case, the noncommutative tensor product $T\boxtimes_RS$ can be calculated analogously and turns out to be isomorphic to the quiver algebra of the same quiver as the above, but with only one zero relation: $\alpha\beta=0$. Clearly, this noncommutative tensor product $T\boxtimes_RS$ is a quasi-hereditary algebra and has $M_2(k)_0$ as its quotient algebra, as the foregoing general fact indicated.

Note that the noncommutative tensor products in both (i) and (ii) are not derived equivalent to the coproduct $S\sqcup_RT$ of $\lambda$ and $\mu$. In fact, $\lambda$ and $\mu$ are independent of the choices of structure maps $f$ and $g$, and moreover, $S\sqcup_RT$ is given by the following quiver algebra
$$\xymatrix{1\,\bullet\ar^-{\bf{\alpha}}@/^0.8pc/[r]
&\bullet\, 2\ar^-{\bf{\beta}}@/^0.8pc/[l]}$$ which is
infinite-dimensional and hereditary. Note that if a $k$-algebra is
derived equivalent to another finite-dimension $k$-algebra, then the
algebra itself must be finite-dimensional. Since the noncommutative
tensor products in both (i) and (ii) are finite-dimensional, they
are not derived equivalent to $S\sqcup_RT$.

\subsubsection{From strictly pure extensions}
An extension $D\subseteq C$ of rings is said to be \emph{strictly
pure} if $C$ has an ideal $X$ such that there exists a splitting
$C=D\oplus X$ of $D$-$D$-bimodules. Such a kind of extensions was
used by Waldhausen to compute the algebraic $K$-theory of
generalized free products in \cite{wald1}.

Now, let $\lambda:R\to S$ and $\mu:R\to T$ be two arbitrary strictly
pure extensions. We shall construct an exact context $(\lambda,\mu,
M,m)$ from the pair $(\lambda, \mu)$. First of all, we fix two split
decompositions of $R$-$R$-bimodules:
$$
S=R\oplus X\;\;\mbox{and}\;\; T=R\oplus Y
$$
where $X$ and $Y$ are ideals of $S$ and $T$, respectively, and
define $M:=R\oplus X\oplus Y$, the direct sum of abelian groups.
Next, we endow $M$ with a ring structure such that $S$ and $T$ are
subrings of $M$. Here, we define a multiplication on $M$ as follows:
$$
(r_1+x_1+y_1)(r_2+x_2+y_2):=r_1r_2+(r_1x_2+x_1r_2+x_1x_2)+(r_1y_2+y_1r_2+y_1y_2)
$$
for $r_i\in R$, $x_i\in X$ and $y_i\in Y$ with $i=1,2$. In
particular, we have $x_1y_1=0=y_1x_1$ in $M$.  One can check that,
under this multiplication, $M$ is a ring with identity $1$, and
contains both $S$ and $T$ as subrings. Since the intersection of $S$
and $T$ in $M$ is equal to $R$ and since $M=S+T$, we see that the
quadruple $(\lambda, \mu, M, 1)$ is an exact context. Clearly, $\Tor_j^R(T,S)=0$ if and only if $\Tor_j^R(Y_R, {}_RX)=0$.

Now, we calculate the noncommutative tensor product $T\boxtimes_RS$ of the exact content $(\lambda, \mu, M, 1)$.

Actually, as $R$-$R$-bimodules, we have
$$T\otimes_RS=R\oplus X\oplus Y\oplus Y\otimes_RX.$$ In this case,
the map $\gamma:S\otimes_RT\to M$ is given by $s\otimes t\mapsto st$
for $s\in S$ and $t\in T$, and the map $\beta:M\to T\otimes_RS$ is
exactly the canonical inclusion. It follows that
$\delta:S\otimes_RT\to T\otimes_RS$ is defined as follows:
$$
(r+x)\otimes(r'+y)\mapsto rr'+ry+xr'
$$
for $r, r'\in R$, $x\in X$ and $y\in Y$. In particular, we have
$(x\otimes y)\delta=0$. Now, we can check that the multiplication
$\circ:(T\otimes_RS)\times (T\otimes_RS)\to T\otimes_RS$ is actually
given by
$$
\big(r_1+x_1+y_1+y_3\otimes x_3\big)\circ \big(r_2+x_2+y_2+y_4\otimes x_4\big)$$
$$
=r_1r_2+(r_1x_2+x_1r_2+x_1x_2)+(r_1y_2+y_1r_2+y_1y_2)+\big(y_1\otimes x_2+y_3\otimes
(x_3r_2)+(r_1y_4)\otimes x_4+(y_1y_4)\otimes x_4+y_3\otimes (x_3x_2)\big).
$$
where $r_1,r_2\in R$, $x_i\in X$ and $y_i\in Y$ for $1\leq i\leq 4$.
Here, we have $x_1\circ y_2=0$ and $y_1\circ x_2=y_1\otimes x_2$.
Moreover, the following map
$$\pi:T\otimes_RS\lra M,\;r_1+x_1+y_1+y_3\otimes x_3\mapsto
r_1+x_1+y_1$$ is a surjective ring homomorphism with
$\beta\pi=Id_M$. Note that $\beta$ may not be a ring homomorphism in
general.

In the following, we show that noncommutative tensor products
induced from strictly pure extensions cover the trivially twisted
extensions in \cite{x6}.

Let $A$ be an Artin algebra, and let $A_0$, $A_1$ and $A_2$ be three Artin subalgebras of $A$ with the same identity. We say that $A$ decomposes as a \emph{twisted tensor product} of $A_1$ and $A_2$ over $A_0$ (see \cite{xx}) if the following three conditions hold:

$(1)$ $A_0$ is a semisimple $k$-algebra such that $A_1\cap A_2=A_0$ and $A=A_0\oplus \rad(A)$ as a direct sum of $A_0$-$A_0$-bimodules, where $\rad(A)$ denotes the Jacobson radical of $A$.

$(2)$ The multiplication map $\sigma: A_2\otimes_{A_0}A_1\to A$ is an isomorphism of $A_2$-$A_1$-bimodules.

$(3)$ $\rad(A_1)\rad(A_2)\subseteq \rad(A_2)\rad(A_1)$.

\smallskip
Now, we assume that $A$ decomposes as a twisted tensor product of $A_1$ and $A_2$ over $A_0$. Then we always have the following decompositions of $A_0$-$A_0$-bimodules:
$$A_1=A_0\oplus \rad(A_1)\;\;\mbox{and}\;\;A_2=A_0\oplus \rad(A_2),$$
where $A_0$ is a common semisimple subalgebra of $A$, $A_1$ and
$A_2$. If $\rad(A_1)\rad(A_2)=0$, then $A$ is called \emph {the
trivially twisted tensor product} of $A_1$ and $A_2$ over $A_0$.

Let $A$ be the trivially twisted tensor product of $A_1$ and $A_2$ over $A_0$. Then we may take
$$R:=A_0, \;\;S:=A_1,\;\;T:=A_2,\;\;X:=\rad(A_1),\;\;Y:=\rad(A_2),$$
and let $\lambda:R\to S$ and $\mu:R\to T$ be the inclusions. Clearly, both $\lambda$ and $\mu$ are strictly pure.
By the foregoing discussion, $M:=R\oplus X\oplus Y$ is a ring and $(\lambda, \mu, M, 1)$ is an exact context. So
the noncommutative tensor product $T\boxtimes_RS$ of this exact context can be defined.
Since $XY=\rad(A_1)\rad(A_2)=0$ in $A$,  the multiplication of the noncommutative tensor product $T\boxtimes_RS$ implies that the map $\sigma:T\boxtimes_RS\to A$ is actually an isomorphism of rings. Thus $A\simeq T\boxtimes_RS$ as rings.

We do not know whether all twisted tensor products of Artin algebras can be realized as the noncommutative tensor products of some exact contexts.

\section{Recollements arising from exact contexts\label{sect5}}

In this section, we shall give a procedure to construct recollements of derived module categories of rings from exact contexts.

Throughout this section, we assume that $(\lambda:R\ra S, \mu: R\ra T, M, m)$ is an exact context.

\subsection{Proof of Theorem \ref{th1}\label{sect5.1}}
In the following, we shall first show that noncommutative tensor
products $T\boxtimes_RS$ can be used to describe noncommutative
localizations.

Let $$B:=\left(\begin{array}{lc} S & M\\
0 & T\end{array}\right),\quad  C:=\left(\begin{array}{lc} T\boxtimes_RS & T\boxtimes_RS\\
T\boxtimes_RS& T\boxtimes_RS\end{array}\right).$$ We define a ring homomorphism
$$\theta:=\left(\begin{array}{cc} \rho & \beta\\
0 & \phi\end{array}\right): \; B \lra C. $$
See Section \ref{sect3.1} for notation.

Furthermore, let
$$e_1:= \left(\begin{array}{ll} 1 &0\\
0 & 0\end{array}\right),\; e_2:=\left(\begin{array}{ll} 0 &0\\
0 & 1\end{array}\right)\in B\quad \mbox{and}\quad\varphi: Be_1\lra Be_2,\;\left(\begin{array}{l} s \\
0 \end{array}\right)\mapsto \left(\begin{array}{c} sm\\
0 \end{array}\right)\;\;\mbox{for}\;\; s\in S.$$ Then $\varphi$ is a
homomorphism of finitely generated projective $B$-modules. If we
identify $\Hom_B(Be_1, Be_2)$ with $M$, then $\varphi$ corresponds
to the element $m\in M$. Let us now look at the noncommutative
localization $\lambda_\varphi:B\to B_\varphi$ of $B$ at $\varphi$.

\smallskip
\begin{Lem}\label{universal}
Up to isomorphism, the map $\theta:B\to C$ is the noncommutative localization of $B$ at $\varphi$.
 \end{Lem}

{\it Proof.}  We first recall a characterization of $B_\varphi$ in terms of generators and relations in \cite{Sheiham}.

Let $\Lambda$ be the ring defined by the following generators and relations:

Generators: $a_x$ for any $x\in M$;

Relations:
$$(1)\;a_m=1; \quad (2)\;a_x+a_{y}=a_{x+y}\;\,\mbox{for}\;y\in M;\quad  (3)\;a_{sm}a_{x}=a_{sx}\;\,\mbox{for}\; s\in S; \quad (4)\;a_xa_{mt}=a_{xt}\;\,\mbox{for}\; t\in T.$$

Define
$\rho_S:S\ra \Lambda,\, s\mapsto a_{sm}, \quad \rho_T: T\ra \Lambda, \,t\mapsto a_{mt} \quad\mbox{and}\quad
\rho_M: M\ra \Lambda,\, x\mapsto a_x$ for $s\in S$, $t\in T$ and
$x\in M$. Then $\rho_S$ and $\rho_T$ are ring homomorphisms.
Moreover, by \cite[Theorem 2,4]{Sheiham}, the noncommutative localization
$\lambda_\varphi:B\to B_\varphi$ is (isomorphic to) the following
map
$${\left(\begin{array}{cc} \rho_S &\rho_M\\
0 & \rho_T\end{array}\right)}: \; \left(\begin{array}{lc} S & M\\
0 & T\end{array}\right)\lra  \left(\begin{array}{lc} \Lambda & \Lambda\\
\Lambda& \Lambda\end{array}\right).$$

Let $$\omega: T\boxtimes_RS \lra \Lambda, \;\,t\otimes s\mapsto (t)\rho_T (s)\rho_S=a_{mt}a_{sm}$$ for $t\in T$ and
$s\in S$. In the following, we shall show that $\omega$ is a ring isomorphism such that $\rho\omega=\rho_S$, $\beta\omega=\rho_M$ and $\phi\omega=\rho_T$. Thus, up to isomorphism, the map $\theta$ can be regarded as the noncommutative localization of $B$ at $\varphi$.
This also means that the noncommutative tensor product of an exact context can be described by generators and relations.

Now, we show that $\omega$ is a ring homomorphism. Clearly, $(1\otimes
1)\omega=a_ma_m=a_m$ by the relation $(3)$. To show that $\omega$
preserves multiplications, that is,
$$\big((t_1\otimes s_1)\circ (t_2\otimes s_2)\big)\omega=(t_1\otimes s_1)\omega (t_2\otimes s_2)\omega.$$
for $s_i\in S$ and $t_i\in T$ for $i=1,2$, we pick up $u\in S$ and
$v\in T$ such that $s_1mt_2=um+mv$. Then $\big((t_1\otimes s_1)\circ
(t_2\otimes s_2)\big)\omega=\big(t_1(1\otimes u+v\otimes
1)s_2\big)\omega=(t_1\otimes us_2+t_1v\otimes
s_2)\omega=(t_1)\rho_T(us_2)\rho_S+(t_1v)\rho_T(s_2)\rho_S=(t_1)\rho_T(u)\rho_S(s_2)\rho_S+(t_1)\rho_T(v)\rho_T(s_2)\rho_S=
(t_1)\rho_T\big( (u)\rho_S+(v)\rho_T \big)(s_2)\rho_S.$ Note that
$(t_1\otimes s_1)\omega (t_2\otimes s_2)\omega=(t_1)\rho_T
(s_1)\rho_S (t_2)\rho_T (s_2)\rho_S$. So it is sufficient to prove
that $(u)\rho_S+(v)\rho_T=(s_1)\rho_S (t_2)\rho_T$, or equivalently,
that $a_{um}+a_{mv}=a_{s_1m}a_{mt_2}$. Actually, due to the
relations $(2)$ and $(3)$, we obtain
$$a_{um}+a_{mv}=a_{um+mv}=a_{s_1mt_2}=a_{s_1m}a_{mt_2}.$$ Thus
$\omega$ is a ring homomorphism.

Next, we show that $\omega$ is a bijection. In fact, the element
$1\otimes 1$ is the identity of $T\boxtimes_RS$ and
$(m)\beta=1\otimes 1$ by Lemma \ref{ring} (3), Moreover, for any
$s\in S$, $t\in T$ and $x\in M$, we have
$(sm)\beta\circ(x)\beta=(s)\rho\circ
(m)\beta\circ(x)\beta=(s)\rho\circ (x)\beta=(sx)\beta$ and
$(x)\beta\circ (mt)\beta=(x)\beta\circ (m)\beta\circ
(t)\phi=(x)\beta\circ (t)\phi=(xt)\beta$. This implies that there
exists a unique ring homomorphism $\psi:\Lambda\to T\boxtimes_RS$
sending $a_x$ to $(x)\beta$. Now, we check that
$\omega\psi=Id_{T\otimes_RS}$ and $\psi\omega=Id_M$. Indeed, the
former follows from
$$
(t\otimes s)\omega\psi=(a_{mt}a_{sm})\psi=(mt)\beta\circ (sm)\beta=(m)\beta\circ(t)\phi\circ (s)\rho\circ (m)\beta=(t)\phi\circ (s)\rho=(t\otimes 1)\circ (1\otimes s)=t\otimes s,
$$
while the latter follows from
$$
(\diamondsuit)\quad
(a_x)\psi\omega=\big((x)\beta\big)\omega=(t_x\otimes 1+1\otimes s_x)\omega=a_{mt_x}a_m+a_ma_{s_xm}=a_{mt_x}+a_{s_xm}=a_x
$$
where $s_x\in S$ and $t_x\in T$ such that $x=s_xm+mt_x$. Thus $\omega$ is a ring isomorphism.

Note that $\beta\omega=\rho_M$ by $(\diamondsuit)$. Since $(s)\rho\omega=(1\otimes s)\omega=a_ma_{sm}=a_{sm}=(s)\rho_S$ and $(t)\phi\omega=(t\otimes 1)\omega
=a_{mt}a_m=a_{mt}=(t)\rho_T$, we see that $\rho\omega=\rho_S$ and $\phi\omega=\rho_T$. $\square$

\begin{Rem}\label{coproduct}
$(1)$
If $(\lambda,\mu)$ is an exact pair, then it follows from Lemma \ref{universal} and \cite[Theorem 4.10, p. 59]{Sch} that the noncommutative tensor product $T\boxtimes_RS$, together with the ring homomorphisms $\rho:S\to T\boxtimes_RS$ and $\phi:T\to T\boxtimes_RS$, is the coproduct $S\sqcup_RT$ of the $R$-rings $S$ and $T$ over $R$ (via the ring homomorphisms $\lambda:R\to S$ and $\mu:R\to T$), that is the \emph{push-out} in the category of $R$-rings. In this case,
the map $\theta: B\to C$ is actually given by the following:
$$\left(\begin{array}{lc} S & S\otimes_RT\\
0 & T\end{array}\right)\lra \left(\begin{array}{lc} T\boxtimes_RS & T\boxtimes_RS\\
T\boxtimes_RS& T\boxtimes_RS\end{array}\right), $$
$$
\quad \left(\begin{array}{cc} s_1 & s_2\otimes t_2\\
0 & t_1\end{array}\right)\mapsto\left(\begin{array}{cc} (s_1)\rho
&(s_2)\rho(t_2)\phi\\
0 & (t_1)\phi \end{array}\right)
$$
for  $s_i\in S$ and $t_i\in T$ with $i=1,2$.

In fact, since $(\lambda,\mu)$ is an exact pair, we have $M=S\otimes_RT$,  $\alpha=Id_M$ and $\delta=\beta:S\otimes_RT\to T\boxtimes_RS$ (see Section \ref{sect3.1} for notation). Further,
$\delta$ is equal to the following map
$$
S\otimes_RT\lra T\boxtimes_RS,\; s_2\otimes t_2\mapsto (s_2)\rho
\circ (t_2)\phi.
$$
In general, for an exact context, its noncommutative tensor product
may not be isomorphic to the coproduct of the $R$-rings $S$ and $T$.

$(2)$ If $\lambda$ is a ring epimorphism, then $T\boxtimes_RS\simeq\End_T(T\otimes_RS)$ as rings.

Actually, in this case, the pair $(\lambda,\mu)$ is an exact pair by
Corollary \ref{app1} and Lemma \ref{exact pair}. It follows from
$(1)$ that $S\sqcup_RT=T\boxtimes_RS$. Further, the ring
homomorphism $\phi:T\to T\boxtimes_RS$ is a ring epimorphism. Thus
$T\boxtimes_RS\simeq \End_{T\boxtimes_RS}(T\boxtimes_RS)\simeq
\End_T(T\boxtimes_RS)= \End_T(T\otimes_RS)$ as rings.

\end{Rem}

From now on, let $\cpx{P}$ be the complex $$0\lra Be_1\lraf{\varphi}
Be_2\lra 0$$ in $\C B$ with $Be_1$ and $Be_2$ in degrees $-1$ and
$0$, respectively, that is $\cpx{P}=\Cone(\varphi)$. Further, let
$\cpx{P}{^*}:=\Hom_B(\cpx{P},B)$ which is isomorphic to the complex
$$0\lra e_2B\lraf{\varphi_*} e_1B\lra 0$$ in $\C {B^{op}}$ with $e_2B$
and $e_1B$ in degrees $0$ and $1$, respectively.

Note that $Be_1$ and $Be_2$ are also right $R$-modules via $\lambda:R\to S$
and $\mu:R\to T$, respectively, and that the map $\cdot m: S\to M$ is a homomorphism
of $S$-$R$-bimodules. Thus $\varphi$ is actually a
homomorphism of $B$-$R$-bimodules. This implies that $\cpx{P}$ is a bounded complex over $B\otimes_\mathbb{Z}R\opp$, and that there is a distinguished triangle in $\K{B\otimes_\mathbb{Z}R\opp}$:
$$Be_1\lraf{\varphi}Be_2\lra\cpx{P}\lra Be_1[1].$$
By Lemma \ref{universal}, the ring homomorphism $\theta: B\to C$ is
a ring epimorphism, and therefore the restriction functor
$\theta_*:C\Modcat\to B\Modcat$ is fully faithful. Now, we define a
full subcategory of $\D{B}$:
$$\mathscr{D}(B)_{C{\tiny\mbox{-Mod}}}:=\{\cpx{X}\in\D{B}~|~H^n(\cpx{X})\in
C\Modcat \mbox{ for\, all \,}n\in\mathbb{Z}\}.
$$
Clearly, we have $X[n]\in \D{B}_{C{\tiny\mbox{-Mod}}}$ for all $X\in
C\Modcat$ and all $n\in \mathbb{Z}$. Also, by \cite[Proposition 3.3
(3)]{xc1}, we have
$$\mathscr{D}(B)_{C{\tiny\mbox{-Mod}}}=
\small\small{\Ker\big(\Hom_{\D{B}}({\rm{Tria}}(\cpx{P}),-)\big)=
\{\cpx{X}\in\D{B}\mid \Hom_{\D B}(\cpx{P},\cpx{X}[n])=0 \mbox{ for\,
all \,}n\in\mathbb{Z}\}},$$ or equivalently,
$$\mathscr{D}(B)_{C{\tiny\mbox{-Mod}}}=\{\cpx{X}\in\D{B}~|~H^n\big(\cpx{\Hom}_{B}(\cpx{P},\cpx{X})\big)=0
\mbox{ for\, all \,}n\in\mathbb{Z}\}.$$

\smallskip
The following result is taken from \cite[Proposition 3.6 (a) and (b)
(4-5)]{xc1}. See also \cite[Theorem 0.7 and Proposition 5.6]{nr1}.

\begin{Lem}\label{good}
Let ${i_*}$ be the canonical embedding of $\D{B}_{C{\tiny\mbox{-{\rm
Mod}}}}$ into $\mathscr{D}(B)$. Then there is a recollement
$$\xymatrix@C=1.2cm{\D{B}_{C{\tiny\mbox{-{\rm Mod}}}}\ar[r]^-{{i_*}}
&\D{B}\ar[r]\ar@/^1.2pc/[l]\ar_-{i^*}@/_1.2pc/[l]
&{\rm{Tria}}{(\cpx{P})} \ar@/^1.2pc/[l]\ar@/_1.2pc/[l]}$$

\medskip
\noindent such that $i^*$ is the left adjoint of $i_*$. Moreover,
the map $\theta: B\to C$ is homological if and only if
$H^n\big({i_*i^*}(B)\big)=0$ for all $n\neq 0$. In this case, the
derived functor $D(\theta_*): \D{C}\to\D{B}_{C{\tiny{\rm -Mod}}}$ is
an equivalence of triangulated categories.
\end{Lem}

To realize ${\rm{Tria}}{(\cpx{P})}$ in Lemma \ref{good} by the
derived module category of a ring, we first show that $\cpx{P}$ is a
self-orthogonal complex in $\D B$. Recall that a complex $\cpx{X}$
in $\D{B}$ is called \emph{self-orthogonal} if
$\Hom_{\D{B}}(\cpx{X},\cpx{X}[n])=0$ for any $n\ne 0$.

\begin{Lem} \label{sum}
The following statements are true:

$(1)$ $\End_{\D{B}}(\cpx{P})\simeq R$ as rings.

$(2)$ The complex $\cpx{P}$ is self-orthogonal in $\D{B}$, that is $\Hom_{\D B}\big(\cpx{P},\cpx{P}[n]\big)=0$ for any $n\neq 0$.

$(3)$ There exists a recollement of triangulated categories:
$$(\star)\quad
\xymatrix@C=1.2cm{\D{B}_{C{\tiny\mbox{-{\rm Mod}}}}\ar[r]^-{{i_*}}
&\D{B}\ar_-{i^!}@/^1.2pc/[l]\ar[r]^-{\;j^!}\ar@/^1.2pc/[l]\ar_-{i^*}@/_1.2pc/[l]
&\D{R} \ar_-{j_*}@/^1.2pc/[l]\ar@/^1.2pc/[l]\ar@/_1.2pc/[l]_{j_!\;}}$$
where $i_*$ is the canonical embedding and
$$
j_!:={_B}\cpx{P}\otimesL_R-,\;\;j^!:=\cpx{\Hom}_B(\cpx{P},-)\simeq {_R}\cpx{P}{^*}\cpx{\otimes}_B-,\;\;
j_*:={\mathbb{R}}\Hom_R(\cpx{P}{^*},-).
$$

\end{Lem}

{\it Proof.}
$(1)$ Note that $\cpx{P}$ is a bounded complex over $B$
consisting of finitely generated projective $B$-modules. It follows
that $\End_{\D B}(\cpx{P})\simeq\End_{\K B}(\cpx{P})$ as rings.
Since $\Hom_B(Be_2,Be_1)=0$, we clearly have
$\End_{\K B}(\cpx{P})=\End_{\C B}(\cpx{P})$. Moreover, if $\End_B(Be_1)$
and $\End_B(Be_2)$ are identified with $S$ and $T$, respectively,
then we can identify $\End_{\C B}(\cpx{P})$ with $K:=\{(s,t)\in S\oplus T\mid sm=mt\}$
which is a subring of $S\oplus T$. Since $(\lambda,\mu, M, m)$ is an exact context, we see that $R\simeq K$ as rings. Thus $\End_{\D B}(\cpx{P})\simeq R$ as rings.

$(2)$ It is clear that $\Hom_{\D B}\big(\cpx{P},\cpx{P}[n]\big)
\simeq \Hom_{\K B}\big(\cpx{P},\cpx{P}[n]\big)=0$ for all
$n\in\mathbb{Z}$ with $|n|\ge 2$. Since $\Hom_B(Be_2,Be_1)$ = $0$,
we get $\Hom_{\D B}\big(\cpx{P},\cpx{P}[-1]\big)=0$. Observe that
$\Hom_{\K B}\big(\cpx{P},\cpx{P}[1]\big)=0$ if and only if
$\Hom_B(Be_1,Be_2)=\varphi\End_B(Be_2)+\End_B(Be_1)\varphi$. If we
identify $\Hom_B(Be_1,Be_2)$, $\End_B(Be_1)$ and $\End_B(Be_2)$ with
$M$, $S$ and $T$, respectively, then the latter
condition is equivalent to that the map $$\left(\begin{array}{l} \;\;\cdot m \\
-m\cdot \end{array}\right): \; S\oplus T\lra M, \quad (s,t)\mapsto
sm-mt, \;\;\mbox{for}\;\; s\in S, \; t\in T,$$ is surjective. Clearly,
this is guaranteed by the definition of exact contexts. Thus $(2)$
holds.

(3) The idea of our proof is motivated by \cite{Keller}.
Since $\cpx{P}$ is a complex of $B$-$R$-bimodules, the total
left-derived functor $\cpx{P}\otimesL_R-:\D{R}\to\D{B}$ and the
total right-derived functor
${\mathbb{R}}\Hom_B(\cpx{P},-):\D{B}\to\D{R}$ are well defined.
Moreover, since $\cpx{P}$ is a bounded complex of finitely generated
projective $B$-modules, the functor
$\cpx{\Hom}_B(\cpx{P},-):\K{B}\to\K{R}$ preserves acyclicity, that
is, $\cpx{\Hom}_B(\cpx{P},\cpx{W})$ is acyclic whenever
$\cpx{W}\in\C{B}$ is acyclic. This automatically induces a derived
functor $\D{B}\ra\D{R}$, which is defined by $\cpx{W}\mapsto
\cpx{\Hom}_B(\cpx{P},\cpx{W})$.  Therefore, we can replace
${\mathbb{R}}\Hom_B(\cpx{P},-)$  with the Hom-functor
$\cpx{\Hom}_B(\cpx{P},-)$ up to natural isomorphism.

Now, we claim that the functor $\cpx{P}\otimesL_R-$ is fully
faithful and induces a triangle equivalence from $\D{R}$ to
${\rm{Tria}}{(\cpx{P})}$.

To prove this claim, we first show that the functor
$\cpx{P}\otimesL_R-: \D{R}\lra \D{B}$ is fully faithful.

Let $$\mathscr{Y}:=\{\cpx{Y}\in\D{R} \mid
\cpx{P}\otimesL_R-:\Hom_{\D R}(R,\cpx{Y}[n])\lraf{\simeq} \Hom_{\D
B}(\cpx{P}\otimesL_RR, \cpx{P}\otimesL_R\cpx{Y}[n])\mbox{ for\, all
\,}n\in\mathbb{Z}\}.$$ Clearly, $\mathscr{Y}$ is a full triangulated
subcategory of $\D{R}$. Since $\cpx{P}\otimesL_R-$ commutates with
arbitrary direct sums and since $\cpx{P}$ is compact in $\D{B}$, we
see that $\mathscr{Y}$ is closed under arbitrary direct sums in $\D{R}$.

In the following, we shall show that $\mathscr{Y}$ contains $R$. It
is sufficient to prove that

$(a)\,\cpx{P}\otimesL_R-$ induces an isomorphism of rings from
$\End_{\D{R}}(R)$ to $\End_{\D{R}}(\cpx{P}\otimesL_RR)$, and

$(b)\,\Hom_{\D B}(\cpx{P}\otimesL_RR, \cpx{P}\otimesL_RR[n])=0$ for
any $n\neq 0$.

\smallskip
Since $\cpx{P}\otimesL_RR\simeq \cpx{P}$ in $\D{B}$, we know that
$(a)$ is equivalent to that the right multiplication map
$R\to\End_{\D{R}}(\cpx{P})$ is an isomorphism of rings, and that
$(b)$ is equivalent to $\Hom_{\D B}(\cpx{P}, \cpx{P}[n])=0$ for any
$n\neq 0$. Actually, $(a)$ and $(b)$ follow directly from (1) and
(2), respectively. This shows $R\in\mathscr{Y}$.

Thus we have $\mathscr{Y}=\D{R}$ since $\D{R}={\rm{Tria}}{(R)}$.
Consequently, for any $\cpx{Y}\in\D{R}$, there is the following
isomorphism:
$$\cpx{P}\otimesL_R-:\Hom_{\D
R}(R,\cpx{Y}[n])\lraf{\simeq} \Hom_{\D B}(\cpx{P}\otimesL_RR,
\cpx{P}\otimesL_R\cpx{Y}[n])\mbox{ for\, all \,}n\in\mathbb{Z}.$$
Now, fix $\cpx{N}\in\D{R}$ and consider
$$\mathscr{X}_{\cpx{N}}:=\{\cpx{X}\in\D{R} \mid \cpx{P}\otimesL_R-:\Hom_{\D
R}(\cpx{X},\cpx{N}[n])\lraf{\simeq} \Hom_{\D
B}(\cpx{P}\otimesL_R\cpx{X}, \cpx{P}\otimesL_R\cpx{N}[n])\mbox{
for\, all \,}n\in\mathbb{Z}\}.$$ Then, one can check that
$\mathscr{X}_{\cpx{N}}$ is a full triangulated subcategory of
$\D{R}$, which is closed under arbitrary direct sums in $\D{R}$.
Since $R\in\mathscr{X}_{\cpx{N}}$ and $\D{R}={\rm{Tria}}{(R)}$, we
get $\mathscr{X}_{\cpx{N}}=\D{R}$. Consequently, for any
$\cpx{M}\in\D{R}$, we have the following isomorphism:
$$\cpx{P}\otimesL_R-: \Hom_{\D R}\big(\cpx{M},\cpx{N}[n]\big)\lraf{\simeq}
\Hom_{\D B}\big(\cpx{P}\otimesL_R\cpx{M},
\cpx{P}\otimesL_R\cpx{N}[n]\big)$$ for all $n\in\mathbb{Z}$. This
means that $\cpx{P}\otimesL_R-:\D{R}\to\D{B}$ is fully faithful.

Recall that $\Tria(\cpx{P})$ is the smallest full triangulated
subcategory of $\D{B}$, which contains $\cpx{P}$ and is closed under
arbitrary direct sums in $\D{B}$. It follows that the image of $\D
R$ under $\cpx{P}\otimesL_R-$ is $\Tria(\cpx{P})$ (see the property
(2) in Section \ref{sect2.1}) and that $\cpx{P}\otimesL_R-$ induces
a triangle equivalence from $\D{R}$ to ${\rm{Tria}}{(\cpx{P})}$.

Note that $\cpx{\Hom}_B(\cpx{P},-)$ is a right adjoint of
$\cpx{P}\otimesL_R-$. This means that the restriction of the functor
$\cpx{\Hom}_B(\cpx{P},-)$ to ${\rm{Tria}}{(\cpx{P})}$ is the
quasi-inverse of the functor
$\cpx{P}\otimesL_R-:\D{R}\to{\rm{Tria}}{(\cpx{P})}.$ In particular,
$\cpx{\Hom}_B(\cpx{P},-)$ induces an equivalence of triangulated
categories:
$$\begin{CD}{\rm{Tria}}{(\cpx{P})}@>{\simeq}>>\D{R}.\end{CD}$$ Furthermore, it follows from
\cite[Proposition 3.3 (3)]{xc1} that
$$\mathscr{D}(B)_{C{\tiny{\rm -Mod}}}=\{\cpx{X}\in\D{B}\mid \Hom_{\D B}(\cpx{P},\cpx{X}[n])=0
\mbox{ for\, all
\,}n\in\mathbb{Z}\}=\Ker\big(\cpx{\Hom}_B(\cpx{P},-)\big).$$
Therefore, we can choose $j_!=\cpx{P}\otimesL_R-$ and
$j^!=\cpx{\Hom}_B(\cpx{P},-)$.

Since $\cpx{P}$ is a bounded complex of $B$-$R$-bimodules with all
of its terms being finitely generated and projective as $B$-modules,
there exists a natural isomorphism of functors (see Section
\ref{sect2.1}):
$$\cpx{P}{^*}\cpx{\otimes}_B-\lraf{\simeq}\cpx{\Hom}_B(\cpx{P},-):\C{B}\lra\C{R}.$$
This implies that the former functor preserves acyclicity, since the
latter always admits this property. It follows that the functors
$\cpx{P}{^*}\otimesL_B-$ and
$\cpx{P}{^*}\cpx{\otimes}_B-:\D{B}\to\D{R}$ are naturally
isomorphic, and therefore $j^!\simeq \cpx{P}{^*}\otimesL_B-$.
Clearly, the functor $\cpx{P}{^*}\otimesL_B-$ has a right adjoint
${\mathbb{R}}\Hom_R(\cpx{P}{^*},-)$. This means that the functor
$j^!$ can also have ${\mathbb{R}}\Hom_R(\cpx{P}{^*},-)$ as a right
adjoint functor (up to natural isomorphism). However, by the
uniqueness of adjoint functors in a recollement, we see that $j_*$
is naturally isomorphic to ${\mathbb{R}}\Hom_R(\cpx{P}{^*},-)$.
Thus, we can choose $j_*={\mathbb{R}}\Hom_R(\cpx{P}{^*},-)$. This
finishes the proof of $(3)$. $\square$

\begin{Lem} \label{calculation}
The following statements hold true:

$(1)$ $i_*i^*(Be_1)\simeq i_*i^*(Be_2)$ in $\D{B}$.

$(2)$ $H^n\big(i_*i^*(Be_1)\big)\simeq
\left\{\begin{array}{ll} 0 & \mbox{if}\; n>0,\\
\Tor^R_{-n}\,(T,S)\oplus \Tor^R_{-n}\,(T,S) & \mbox{if }\; n\leq 0.\end{array} \right.$

\smallskip
$(3)$ If $\lambda:R\to S$ is homological, then $i_*i^*(Be_1)\simeq Be_2\otimesL_RS$ in $\D{B}$.

\end{Lem}

{\it Proof.} We keep the notation introduced in Lemma \ref{sum}.

$(1)$ Applying the triangle functor $i_*i^*:\D{B}\to \D{B}$ to the
distinguished triangle:
$$\cpx{P}[-1]\lra Be_1\lraf{\varphi}Be_2\lra\cpx{P}$$
in $\D{B}$, we obtain another distinguished triangle in $\D{B}$:
$$i_*i^*(\cpx{P})[-1]\lra i_*i^*(Be_1)\lraf{i_*i^*(\varphi)} i_*i^*(Be_2)\lra
i_*i^*(\cpx{P}).$$ Since the composition functor
$i^*j_!:\D{R}\to\D{B}_{C{\tiny\mbox{-{\rm Mod}}}}$ is zero in the
recollement $(\star)$, we clearly have $i^*(\cpx{P})\simeq
i^*j_!(R)=0$. Thus $i_*i^*(\varphi): i_*i^*(Be_1)\to i_*i^*(Be_2)$
is an isomorphism.

$(2)$ First, we show that if $n>0$ or $n<-1$, then $$H^n(i_*i^*(Be_1))\simeq \Tor^R_{-n}\,(T,S)\oplus \Tor^R_{-n}\,(T,S)$$
where $\Tor^R_{-n}\,(T,S):=0$ for $n>0$.

In fact, let $\varepsilon: j_!j^!\to Id_{\D{B}}$ and $\eta: Id_{\D{B}}\to i_*i^* $
be the counit and unit adjunctions with
respect to the adjoint pairs $(j_!,j^!)$ and $(i^*,i_*)$ in the recollement $(\star)$, respectively.
Then, for any $\cpx{X}\in\D{B}$, there is a canonical triangle in $\D{B}$:
$$j_!j^!(\cpx{X})\lraf{\varepsilon_{\cpx{X}}}\cpx{X}\lraf{\eta_{\cpx{X}}} i_*i^*(\cpx{X})\lra
j_!j^!(\cpx{X})[1].$$
In particular, we have the following triangle in $\D{B}$:
$$j_!j^!(Be_1)\lraf{\varepsilon_{Be_1}}Be_1\lraf{\eta_{Be_1}} i_*i^*(Be_1)\lra
j_!j^!(Be_1)[1].$$

Note that $j^!(Be_1)=\Hom_B(\cpx{P},Be_1)\simeq S[-1]$ as complexes of $R$-modules¡£ In the following, we always identify $\Hom_B(\cpx{P},Be_1)$ with $S[-1]$. Under this identification, we obtain the following triangle in $\D{B}$:
$$\cpx{P}\otimesL_RS[-1]\lraf{\varepsilon_{Be_1}}Be_1\lraf{\eta_{Be_1}} i_*i^*(Be_1)\lra
\cpx{P}\otimesL_RS.$$ Now, for each $n\in\mathbb{Z}$, we apply the
$n$-th cohomology functor $H^n:\D{B}\to B\Modcat$ to this triangle,
and conclude that if $n>0$ or $n<-1$, then $H^n(i_*i^*(Be_1))\simeq
H^n(\cpx{P}\otimesL_RS)$. Moreover, we have $$\cpx{P}=T\oplus (0\ra
S\lraf{\cdot m} M\ra 0)=T\oplus \Cone(\cdot m)\in\D{R\opp}.$$ Sincer
$(\lambda,\mu, M,m)$ is an exact context, it follows from the
diagram $(\sharp)$ that the chain map $(\lambda, m\cdot):
\Cone(\mu)\to \Cone(\cdot m)$ is a quasi-isomorphism. This implies
that $$\Cone(\cdot m)\simeq \Cone(\mu) \; \;\mbox{in}\;\;\D{R\opp}.$$ Thus $\cpx{P}\simeq T\oplus
\Cone(\mu)$ in $\D{R\opp}$ and $\cpx{P}\otimesL_RS\simeq
(T\otimesL_RS)\oplus (\Cone(\mu)\otimesL_RS)$ in $\D{\mathbb Z}$. In
particular, we have
$$H^n(\cpx{P}\otimesL_RS)\simeq H^n(T\otimesL_RS)\oplus
H^n(\Cone(\mu)\otimesL_RS)$$ for all $n\in\mathbb{Z}$. Applying the
functor $-\otimesL_RS$ to the canonical triangle $$R\lraf{\mu}
T\lra\Cone(\mu)\lra R[1]$$ in $\D{R\opp}$, we obtain another
triangle $S\to T\otimesL_RS\to \Cone(\mu)\otimesL_RS\to S[1]$ in
$\D{\mathbb{Z}}$. This implies that if $n>0$ or $n<-1$, then
$H^n(T\otimesL_RS)\simeq H^n(\Cone(\mu)\otimesL_RS)$, and therefore
$$H^n(i_*i^*(Be_1))\simeq H^n(\cpx{P}\otimesL_RS)\simeq H^n(T\otimesL_RS)\oplus
H^n(T\otimesL_RS)\simeq \Tor^R_{-n}\,(T,S)\oplus \Tor^R_{-n}\,(T,S).$$

Next, we shall show that $H^{-1}(i_*i^*(Be_1))\simeq \Tor^R_{1}\,(T,S)\oplus \Tor^R_{1}\,(T,S)$.

Indeed, we have the following two homomorphisms:
$$\sigma: S\otimes_RS\lra S, \quad s_1\otimes{s_2}\mapsto
s_1s_2,\quad\quad\quad\varphi_1:\;S\otimes_RS\lra M\otimes_RS, \quad
s_1\otimes{s_2}\mapsto s_1m\otimes s_2$$ for $s_1,s_2\in S$, and
can identify $Be_1\otimes_RS$ and $Be_2\otimes_RS$ with
$\left(S\otimes_RS\,\atop{0}\right)$ and $\left(M\otimes_RS
\atop{T\otimes_RS}\right)$ as $B$-modules, respectively. Then there
is a chain map in $\C{B}$:
$$\xymatrix{
\cpx{P}\otimes_R(S[-1]):\ar[d]^-{\cpx{g}}&0\ar[r]&\left(S\otimes_RS\,\atop{0}\right)
\ar^-{\left({\sigma\,\atop{0}}\right)}[d]\ar[r]^-{\left({\varphi_1\,\atop{0}}\right)}&
\left(M\otimes_RS \atop{T\otimes_RS}\right)
\ar^-{0}[d]\ar[r]&0\\
Be_1:&0\ar[r]&\left(S\atop{0}\right) \ar[r]^-{0}&0\ar[r]&0}$$ Let
$_pS$ be a deleted projective resolution of the module $_RS$ with
$\tau:{_p}S\ra S$ a quasi-isomorphism. Recall that
$j^!(Be_1)=\Hom_B(\cpx{P},Be_1)=S[-1]$. Then the counit
$\varepsilon_{Be_1}: j_!j^!(Be_1)\lra Be_1$ is just the composite of
the following homomorphisms:
$$j_!j^!(Be_1)=\cpx{P}\otimesL_R\cpx{\Hom}_B(\cpx{P},Be_1)=
\cpx{P}\cpx{\otimes}_R({_p}S)[-1]\lraf{\small{1\otimes\tau[-1]}}
\cpx{P}\cpx{\otimes}_RS[-1]\lraf{\cpx{g}} Be_1.$$
Further, let $\cpx{h}$ be the following chain map:
$$
\xymatrix{
\cpx{P}:\ar[d]^-{\cpx{h}}&0\ar[r]&\left(S\,\atop{0}\right)
\ar@{=}[d]\ar[r]^-{\left(\cdot m\,\atop{0}\right)}&
\left(M \atop{T}\right)
\ar^-{0}[d]\ar[r]&0\\
Be_1[1]:&0\ar[r]&\left(S\atop{0}\right) \ar[r]^-{0}&0 \ar[r]&0}
$$
Then we have a commutative diagram:
$$
\xymatrix{
\cpx{P}\otimes_R({_p}S)[-1]\ar[r]^-{1\otimes\tau[-1]}\ar[d]_-{\cpx{h}\otimes 1}
& \cpx{P}\otimes_RS[-1]\ar[r]^-{\cpx{g}}\ar[d]^-{\cpx{h}\otimes 1} & Be_1\ar@{=}[d]\\
Be_1[1]\otimes_R({_p}S)[-1]\ar[r]^-{1\otimes \tau[-1]} &
Be_1[1]\otimes_RS[-1]\ar[r]^-{\left({\sigma\,\atop{0}}\right)}& Be_1
}$$ This implies that the following diagram
$$(\ast\ast)\quad
\xymatrix{
\cpx{P}\otimesL_R\Hom_B(\cpx{P},Be_1)\ar[r]^-{\cpx{h}\otimesL 1}\ar@{=}[d]&Be_1[1]\otimesL_R\Hom_B(\cpx{P},Be_1)
\ar[d]_-{\big(1\otimes\tau[-1]\big)\,\left({\sigma\,\atop{0}}\right)}\\
\cpx{P}\otimesL_R\Hom_B(\cpx{P},Be_1)\ar[r]^-{\varepsilon_{Be_1}}
&Be_1=\left(S\atop{0}\right)}$$ is commutative in $\D{B}$. Since we
have the following distinguished triangle
$$Be_2\lra\cpx{P}\lraf{\cpx{h}} Be_1[1]\lraf{\varphi[1]} Be_2[1]$$
in $\K{B\otimes_{\mathbb{Z}}R^{\opp}}$, there is a homomorphism $$\xi:Be_2\otimesL_R\Hom_B(\cpx{P},Be_1)[1]\lra i_*i^*(Be_1)$$
in $\D{B}$ and a complex $W\in\D{B}$ such that $(\ast\ast)$ is completed to the following commutative diagram:
{\footnotesize $$\xymatrix{
& W\ar@{=}[r]\ar[d]^-{\zeta} &W\ar[d]^-{\zeta(\varphi[1]\otimesL 1)}&  \\
\cpx{P}\otimesL_R\Hom_B(\cpx{P},Be_1)\ar[r]^-{\cpx{h}\otimesL 1}\ar@{=}[d]&Be_1[1]\otimesL_R\Hom_B(\cpx{P},Be_1)
\ar[r]^-{\varphi[1]\otimesL 1}\ar[d]_-{(1\otimes\tau[-1])\,\left({\sigma\,\atop{0}}\right)}
&Be_2[1]\otimesL_R\Hom_B(\cpx{P},Be_1)\ar@{-->}[d]_-{\xi} \ar[r]&\cpx{P}[1]\otimesL_R\Hom_B(\cpx{P},Be_1)\ar@{=}[d]\\
\cpx{P}\otimesL_R\Hom_B(\cpx{P},Be_1)\ar[r]^-{\varepsilon_{Be_1}}
&Be_1=\left(S\atop{0}\right)\ar[r]^-{\eta_{Be_1}}\ar[d]&i_*i^*(Be_1)\ar[r]\ar[d]& \cpx{P}[1]\otimesL_R\Hom_B(\cpx{P},Be_1)\\
& W[1]\ar@{=}[r]&W[1]& }$$}\noindent with rows and columns being
distinguished triangles in $\D{B}$. Note that such a homomorphism
$\xi$ is unique. In fact, this follows from
$$\Hom_{\D{B}}(\cpx{P}[1]\otimesL_R\Hom_B(\cpx{P},Be_1),
i_*i^*(Be_1))=\Hom_{\D{B}}(j_!(S), i_*i^*(Be_1))=0.$$Now, we obtain
the following triangle in $\D{B}$:
$$W \lraf{\psi} Be_2\otimesL_RS\lraf{\xi} i_*i^*(Be_1)\lra W[1]$$
where $\psi:=\zeta(\varphi[1]\otimesL 1)$. This yields a long exact sequence of abelian groups:
$$
H^{-1}(W)\lraf{H^{-1}(\psi)} H^{-1}(Be_2\otimesL_RS)\lraf{H^{-1}(\xi)} H^{-1}(i_*i^*(Be_1))\lra
H^0(W)\lraf{H^{0}(\psi)} H^0(Be_2\otimesL_RS)
$$

In the sequel, we show that the map $H^{0}(\psi):H^0(W)\to
H^0(Be_2\otimesL_RS)$ is always injective.

Note that $H^0(\psi)$ is the composite of $H^0(\zeta):H^0(W)\to H^0(Be_1\otimesL_RS)$ with
$$H^0(\varphi[1]\otimesL 1):H^0(Be_1\otimesL_RS)\lra H^0(Be_2\otimesL_RS),$$ and that
$H^0(Be_1\otimesL_RS)=Be_1\otimes_RS=S\otimes_RS$ and
$H^0(Be_2\otimesL_RS)= Be_2\otimes_RS$. On the one hand, applying
the functor $H^0$ to the triangle
$$W \lraf{\zeta} Be_1\otimesL_RS\lraf{(1\otimes\tau[-1])\,\left({\sigma\,\atop{0}}\right)} Be_1\lra W[1],$$
we obtain a short exact sequence $$0\lra H^0(W)\lraf{H^0(\zeta)}
S\otimes_RS\lraf{\sigma} S\lra 0,$$ where $Be_1$ is identified with
$S$. This implies that $H^0(\zeta):
H^0(W)\lraf{\simeq}\Ker(\sigma)$. On the other hand, we can identify
$H^0(\varphi[1]\otimesL 1)$ with the map $\varphi\otimes_RS:
Be_1\otimes_RS\to Be_2\otimes_RS$ induced from $\varphi:Be_1\to
Be_2$. Consequently, $H^0(\psi)$ is the composite of
$H^0(\zeta):H^0(W)\to S\otimes_RS=Be_1\otimes_RS$ with
$\varphi\otimes_RS:Be_1\otimes_RS\to Be_2\otimes_RS.$ Thus
$H^0(\psi)$ is injective if and only if so is the restriction of
$\varphi\otimes_RS$ to $\Ker(\sigma)$, while the latter is also
equivalent to saying that the restriction of the map
$\varphi_1=(\cdot m)\otimes_RS: S\otimes_RS\to M\otimes_RS$ to
$\Ker(\sigma)$ is injective. Hence, we need to show that
$\Ker(\varphi_1)\cap \Ker(\sigma)=0$.

In fact, for the ring homomorphism $\lambda:R\to S$, the sequence  $ 0\ra \Ker(\sigma) \ra S\otimes_RS\lraf{\sigma} S\ra 0$ always splits in the the category of $R$-$S$-bimodules since the composite of $\lambda\otimes_RS:R\otimes_RS\to S\otimes_RS$ with $\sigma$ is an isomorphism of $R$-$S$-bimodules.
It follows that $\lambda\otimes_RS$ is injective, $\Img(\lambda\otimes_RS)\cap \Ker(\sigma)=0$ and $S\otimes_RS=\Ker(\sigma)\oplus \Img(\lambda\otimes_RS).$  Now, we apply the tensor functor $-\otimes_RS$ to the diagram $(\sharp)$, which is a push-out and  pull-back diagram in the category of $R$-$R$-bimodules, and obtain another diagram
$$
\xymatrix{
R\otimes_RS\ar[d]_-{\mu\otimes_RS} \ar[r]^-{\lambda\otimes_RS} & S\otimes_RS\ar[d]^-{\varphi_1}\\
T\otimes_RS \ar[r]^-{(m \cdot)\otimes_RS} &  M\otimes_RS}
$$
which is a push-out and pull-back diagram in the category of $R$-$S$-bimodules. This implies that
the map $\lambda\otimes_RS$ induces an isomorphism from $\Ker(\mu\otimes_RS)$ to $\Ker(\varphi_1)$.
In particular, we have $\Ker(\varphi_1)\subseteq \Img(\lambda\otimes_RS)$. It follows from
$\Img(\lambda\otimes_RS)\cap \Ker(\sigma)=0$ that $\Ker(\varphi_1)\cap\Ker(\sigma)=0$.

Thus $H^{0}(\psi):H^0(W)\to H^0(Be_2\otimesL_RS)$ is injective. Consequently, the map $H^{-1}(\xi)$ is surjective and  $H^{-1}(i_*i^*(Be_1))\simeq \Coker(H^{-1}(\psi))$. Observe that $H^{-1}(\psi)$ is the composite of the isomorphism $H^{-1}(\zeta):H^{-1}(W)\lraf{\simeq} H^{-1}(Be_1\otimesL_RS)$ with the map
$$H^{-1}(\varphi[1]\otimesL 1):H^{-1}(Be_1\otimesL_RS)\lra H^{-1}(Be_2\otimesL_RS).$$
Therefore, we have $$H^{-1}(i_*i^*(Be_1))\simeq \Coker(H^{-1}(\psi))\simeq\Coker\big(H^{-1}(\varphi[1]\otimesL 1)\big).$$ So, to show that $H^{-1}(i_*i^*(Be_1))\simeq \Tor^R_{1}\,(T,S)\oplus \Tor^R_{1}\,(T,S)$, it suffices to prove that $$\Coker\big(H^{-1}(\varphi[1]\otimesL 1)\big)\simeq \Tor^R_{1}\,(T,S)\oplus \Tor^R_{1}\,(T,S).$$
Recall that $Be_1=S$, $Be_2=M\oplus T$ and $\varphi=(\cdot m, 0): S\to M\oplus T$ in $R\opp\Modcat$. Moreover, we have $H^{-1}(Be_1\otimesL_RS)=\Tor^R_1(S,S)$ and $H^{-1}(Be_2\otimesL_RS)=\Tor^R_1(M\oplus T, S)$.
In this sense, the map $H^{-1}(\varphi[1]\otimesL 1)$ is actually given by
$$\big(\Tor^R_1(\cdot m, S), 0\big): \Tor^R_1(S,S)\lra \Tor^R_1(M, S)\oplus \Tor^R_1(T,S)$$
Thus  $\Coker\big(H^{-1}(\varphi[1]\otimesL 1)\big)\simeq \Coker(\Tor^R_1(\cdot m,\, S))\oplus \Tor^R_{1}\,(T,S)$.

Finally, we show that $\Coker(\Tor^R_1(\cdot m, S))\simeq \Tor^R_1(T,S)$.

In fact, since the quadruple $(\lambda, \mu, M, m)$ is an exact context, we have the following exact sequence of $R$-$R$-bimodules:
$$(*)\quad
\xymatrix{0\ar[r] & R\ar[r]^-{(\lambda,\,\mu)\,}& S\oplus
\;T\ar[r]^-{\left({\cdot m\,\atop{-m \cdot}}\right)}&
M\ar[r]& 0.}
$$
Applying $\Tor^R_i(-,S)$ for $i=0, 1$ to this sequence, we obtain a long exact sequence of abelian groups:
$$
0=\Tor^R_1(R,S)\lra \Tor^R_1(S,S)\oplus \Tor^R_1(T,S)\lraf{\left({\Tor^R_1(\cdot m,\,S)\,\atop{-\Tor^R_1(m \cdot,\, S)}}\right)}\Tor^R_1(M,S)\lra R\otimes_RS\lraf{\big(\lambda\otimes_RS,\,\mu\otimes_RS\big)} S\otimes_RS\oplus S\otimes_RT.
$$
Since $\lambda\otimes_RS: R\otimes_RS\to S\otimes_RS$ is injective, the map
$$\left({\Tor^R_1(\cdot m,\,S)\,\atop{-\Tor^R_1(m \cdot,\, S)}}\right): \Tor^R_1(S,S)\oplus\Tor^R_1(T,S)\lra\Tor^R_1(M,S)$$ is an isomorphism, which gives rise to
$\Coker(\Tor^R_1(\cdot m, S))\simeq \Tor^R_1(T,S)$. It follows that $$H^{-1}(i_*i^*(Be_1))\simeq\Coker\big(H^{-1}(\varphi[1]\otimesL 1)\big)\simeq \Coker(\Tor^R_1(\cdot m,\, S))\oplus \Tor^R_{1}\,(T,S)\simeq \Tor^R_{1}\,(T,S)\oplus \Tor^R_{1}\,(T,S).$$

Hence, we have shown that $H^n\big(i_*i^*(Be_1)\big)\simeq
\Tor^R_{-n}\,(T,S)\oplus \Tor^R_{-n}\,(T,S)$ for any $n\in\mathbb{Z}$. This finishes the proof of $(2)$.

$(3)$ Note that if $\lambda$ is homological, then both $1\otimes \tau[-1]: Be_1[1]\otimesL_R\Hom_B(\cpx{P},Be_1)\to Be_1[1]\otimes_RS[-1]$ and $\sigma: S\otimes_RS\to S$ are isomorphisms. This implies that the morphism $$(1\otimes\tau[-1])\,\left({\sigma\,\atop{0}}\right):Be_1[1]\otimesL_R\Hom_B(\cpx{P},Be_1)
\lra Be_1$$
is an isomorphism. Thus $W\simeq 0$ and $i_*i^*(Be_1)\simeq Be_2\otimesL_RS$ in $\D{B}$.
This shows $(3)$. $\square$

\medskip
For exact pairs, we establish the following result which will be
used in the proof of Theorem \ref{th1}.

\begin{Lem}\label{gen}
Suppose that $(\lambda,\mu)$ is an exact pair and that $\lambda$ is homological. Then:

$(1)$ $\Tor^R_i(S,T)=0$ for all $i>0$.

$(2)$ Given a commutative diagram of ring homomorphisms:
$$\xymatrix{
  R \ar[d]_-{\mu} \ar[r]^-{\lambda}
                & S \ar[d]_{f}  \\
  T \ar[r]^-{g}
                & \,\; \Gamma,
 }$$
the following statements are equivalent:

$(a)$ The ring homomorphism $g: T\ra \Gamma$ is homological.

$(b)$ The ring homomorphism $$\theta_{f,g}:B\lra M_2(\Gamma),\quad\left(\begin{array}{cc} s_1 & s_2\otimes t_2\\
0 & t_1\end{array}\right)\mapsto\left(\begin{array}{cc} (s_1)f
&(s_2)f(t_2)g\\
0 & (t_1)g \end{array}\right), \, s_i\in S, t_i\in T, i=1,2$$ is
homological.
\end{Lem}

{\it  Proof.} $(1)$ Let $\cpx{Q}$ be the mapping cone of $\lambda$. Then
there is a distinguished triangle in $\D{R}$:
$$R\lraf{\lambda} S\ra \cpx{Q}\ra R[1].$$
Since $\lambda$ is homological, it follows from \cite[Theorem
4.4]{GL} that $\lambda$ induces the following isomorphisms
$$S\lraf{\simeq} S\otimesL_RR\lraf{S\otimesL_R\lambda}
S\otimesL_RS$$ in $\D{S}$. This implies that $S\otimesL_R\cpx{Q}=0$
in $\D{S}$, and therefore $S\otimesL_R\cpx{Q}=0$ in $\D{R}$. Let
$\cpx{\mu}:=(\mu^i)_{i\in\mathbb{Z}}$ be the chain map defined by
$\mu^{-1}:=\mu$, $\mu^0:=\mu\,'$ and $\mu^{\,i}=0$ for $i\ne -1,0.$
Since $(\lambda,\mu)$ is an exact pair, we see that
$\cpx{\mu}:\cpx{Q}\to\cpx{Q}\otimes_RT$ is an isomorphism in
$\D{R}$. It follows that
$S\otimesL_R\big(\cpx{Q}\otimes_RT\big)\simeq S\otimesL_R\cpx{Q}=0$
in $\D S$. Now, applying $S\otimesL_R-$ to the triangle $
T\lraf{\lambda'} S\otimes_RT\lra\cpx{Q}\otimes_RT\lra T[1]$, we
obtain $S\otimesL_RT\simeq S\otimesL_R\big(S\otimes_RT\big)$ in
$\D{S}$ (and also in $\D{R}$). This yields that $\Tor_i^R(S,T)\simeq
\Tor_i^R(S,S\otimes_RT)$ for all $i\ge 0$. As $S\otimes_RT$ is a
left $S$-module and $\lambda$ is homological, it follows that
$\Tor^R_i(S,S\otimes_RT)=\Tor^S_i(S,S\otimes_RT)=0$ for all $i>0$,
and therefore $\Tor^R_i(S,T)=0$. This shows $(1)$.

$(2)$ Set $\Lambda:=M_2(\Gamma)$. Let $e_1:= \left(\begin{array}{ll} 1 &0\\
0 & 0\end{array}\right)$ and $e_2:=\left(\begin{array}{ll} 0 &0\\
0 & 1\end{array}\right)\in B$, and let
$e:=(e_2)\theta_{f,g}\in\Lambda$. Then we have
$e=e^2,\;\End_{\Lambda}(\Lambda e)\simeq \Gamma$ and
$\End_{B}(Be_2)\simeq T$. Observe that $\Lambda e$ is a projective
generator for $\Lambda$-Mod. Then, by Morita theory, the tensor
functor $ e\Lambda\otimes_{\Lambda}-:\Lambda\Modcat\lra
\Gamma\Modcat$ is an equivalence of module categories, which can be
canonically extended to a triangle equivalence
$D(e\Lambda\otimes_{\Lambda}-):\D{\Lambda}\to\D{\Gamma}$.

It is clear that $e_2B\otimes{_B}\Lambda\simeq
e_2\cdot\Lambda=e\Lambda$ as $T$-$\Lambda$-bimodules, where the left
$T$-module structure of $e\Lambda$ is induced by $g:T\to \Gamma$.
Thus the following diagram of functors between module categories
$$\xymatrix{
\Lambda\Modcat\ar^-{e\Lambda\otimes_{\Lambda}-}[rr]\ar_-{\big(\theta_{f,g}\big)_*}[d]
&&\Gamma\Modcat\ar^{g_{*}}[d] \\
B\Modcat\ar^{e_2B\otimes_{B}-}[rr] &&T\Modcat}
$$
is commutative, where  $\big(\theta_{f,g}\big)_*$ and $g_{*}$ stand
for the restriction functors induced by the ring homomorphisms
$\theta_{f,g}$ and $g$, respectively. Since all of the functors
appearing in the diagram are exact, we can pass to derived module
categories and get the following commutative diagram of functors
between derived module categories:
$$(\dag)\quad\quad
\xymatrix{
\D{\Lambda}\ar^-{D\big(e\Lambda\otimes_{\Lambda}-\big)}[rr]\ar_-{D\big((\theta_{f,g})_*\big)}[d]
&&\D{\Gamma}\ar^{D(g_{*})}[d] \\
\D{B}\ar^{D\big(e_2B\otimes_{B}-\big)}[rr] &&\D{T}}
$$
where the functor $D\big(e\Lambda\otimes_{\Lambda}-\big)$ in the
upper row is a triangle equivalence.

Note that $\theta_{f,g}: B\to \Lambda$ (respectively,
$g:T\to\Gamma$\,) is homological if and only if the functor
$D\big((\theta_{f,g})_*\big)$ $\big($respectively, $D(g_{*})\big )$ is fully
faithful. This means that, to prove that $(a)$ and $(b)$ are
equivalent, it is necessary to establish some further connection
between $D\big((\theta_{f,g})_*\big)$ and $D(g_{*})$ in the diagram
$(\dag)$.

Actually, the triangle functor $D(e_2B\otimes_B-)$ induces a
triangle equivalence from ${\rm{Tria}}(Be_2)$ to $\D{T}$. This can
be obtained from the following classical recollement of derived
module categories:
$$\xymatrix{\D{S}\ar[rr]^-{S\otimesL_S-}&&\D{B}\ar[rr]^{D(e_2B\otimes_B-)}\ar@/^1.8pc/[ll]\ar@/_1.8pc/[ll]
&&\D{T}\ar@/^1.8pc/[ll]\ar_-{Be_2\otimesL_{T}-}@/_2.0pc/[ll].}\vspace{0.4cm}$$
which arises form the triangular structure of the ring $B$.

Suppose that the image $\Img\big(D\big((\theta_{f,g})_*\big)\big)$
of the functor $D\big((\theta_{f,g})_*\big)$ belongs to
${\rm{Tria}}(Be_2)$. Then we can strengthen the diagram $(\dag)$ by
the following commutative diagram of functors between triangulated
categories:
$$\xymatrix{
&\ar[ld]_-{D\big((\theta_{f,g})_*\big)}\D{\Lambda}\ar^-{D(e\Lambda\otimes_{\Lambda}-)}_-{\simeq}[rr]
\ar^-{D\big((\theta_{f,g})_*\big)}[d]
&&\D{\Gamma}\ar^{D (g_{*})}[d] \\
\D{B}&\ar@{_{(}->}[l]\;{\rm{Tria}}(Be_2)\ar^-{D(e_2B\otimes_{B}-)}_-{\simeq}[rr]
&&\D{T}}
$$
This implies that $D\big((\theta_{f,g})_*\big)$ is fully faithful if
and only if so is $D (g_{*})$, and therefore $\theta_{f,g}$ is
homological if and only if $g$ is homological.

So, to finish the proof of Lemma \ref{gen} (2), it suffices to prove
that $\Img\big(D\big((\theta_{f,g})_*\big)\big)\subseteq
{\rm{Tria}}(Be_2).$ In the following, we shall concentrate on
proving this inclusion.

In fact, it is known that $\D{\Lambda}={\rm Tria}(\Lambda e)$ and
$D\big((\theta_{f,g})_*\big)$ commutes with small coproducts since
it admits a right adjoint. Therefore, according to the property (2)
in Section \ref{sect2.1}, in order to check the above inclusion, it
is enough to prove that $\Lambda e\in {\rm Tria}(Be_2)$ when considered
as a $B$-module via $\theta_{f,g}$. If we identify $e_2B\otimes_B-$
with the left multiplication functor by $e_2$, then $\Lambda e\in
{\rm Tria}(Be_2)$ if and only if $Be_2\otimesL_{T}e_2\cdot(\Lambda
e)\lraf{\simeq}\Lambda e$ in $\D{B}$. Clearly, the latter is
equivalent to that $\Tor^{T}_n(Be_2,\,e_2\cdot(\Lambda e))=0$ for
any $n>0$ and the canonical multiplication map
$Be_2\otimes_Te_2\cdot(\Lambda e)\ra\Lambda e$ is an isomorphism.

Set $M:=S\otimes_RT$ and write $B$-modules in the form of triples
$(X,Y,h)$ with $X\in T\Modcat,\,Y\in S\Modcat$ and
$h:M\otimes_{T}X\to Y$ a homomorphism of $S$-modules. The morphisms
between two modules $(X,Y,h)$ and $(X',Y',h')$ are pairs of
morphisms $(\alpha, \beta)$, where $\alpha:X\to X'$ and $\beta:Y\to
Y'$ are homomorphisms in $T$-Mod and $S$-Mod, respectively, such
that $h\beta=(M\otimes_{T}\alpha )h'$.

With these interpretations, we rewrite $\Lambda e=(\Gamma, \,\Gamma,
\delta_{\Gamma})\in B\Modcat$, where
$\delta_{\Gamma}:M\otimes_{T}\Gamma\to \Gamma$ is defined by
$(s\otimes t)\otimes\gamma\mapsto (s)f (t)g\gamma\;$ for $s\in S,
t\in T$ and $\gamma\in\Gamma$. Then $e_2\cdot(\Lambda e)=e\Lambda
e\simeq \Gamma$ as left $T$-modules, and $Be_2\simeq M\oplus T$ as
right $T$-modules. Consequently, we have
$$Be_2\otimes_{T}e_2\cdot(\Lambda e)\simeq
Be_2\otimes_{T}\Gamma\simeq (\Gamma,\,M\otimes_{T}\Gamma,\,1)\;\,
\mbox{and}\,\;\Tor^{T}_n(Be_2,\,e_2\cdot(\Lambda
e))\simeq\Tor^T_n(M\oplus T,\Gamma)\simeq\Tor^T_n(M,\Gamma)$$ for
any $n>0$. This implies that the multiplication map
$Be_2\otimes_Te_2\cdot(\Lambda e)\ra\Lambda e$ is an isomorphism if
and only if so is the map $\delta_{\Gamma}$. It follows that
$Be_2\otimesL_{T}e_2\cdot(\Lambda e)\simeq \Lambda e$ in $\D{B}$ if
and only if $\delta_{\Gamma}$ is an isomorphism of $S$-modules and
$\Tor^T_n(M,\Gamma)=0$ for any $n>0$.

In order to verify the latter conditions just mentioned, we shall
prove the following general result:

For any $\Gamma$-module $W$, if we regard $W$ as a left $T$-module
via $g$ and an $S$-module via  $f$, then the map
$\delta_{W}:M\otimes_{T}W\to W$,  defined by $(s\otimes t)\otimes
w\mapsto (s)f (t)g\,w\;$ for $s\in S, t\in T$ and $w\in W$, is an
isomorphism of $S$-modules, and $\Tor^{T}_i(M, W)=0$ for any $i>0$.

To prove this general result, we fix a projective resolution
$\cpx{V}$ of $S_R$:
$$\cdots
\lra V^n\lra V^{n-1}\lra\cdots \lra V^1\lra V^0 \lra S_R\lra 0$$
with $V^i$ a projective right $R$-module for each $i$. By $(1)$,
we have $\Tor^R_j(S,T)=0$ for any $j>0$. It follows
that the complex $\cpx{V}\otimes_RT$ is a projective resolution of
the right $T$-module $M$. Thus the following isomorphisms of
complexes of abelian groups:
$$\big(\cpx{V}\otimes_RT\big)\otimes_{T}W\simeq
\cpx{V}\otimes_R\big(T\otimes_{T}W\big)\simeq\cpx{V}\otimes_RW$$
imply that $\Tor^{T}_i(M, W)\simeq\Tor^{R}_i(S, W)$ for any $i>0$.
Recall that $W$ admits an $S$-module structure via the map $f$.
Moreover, it follows from $\lambda f=\mu\, g$ that the $R$-module
structure of $W$ endowed via the ring homomorphism $\mu\, g$ is the
same as the one endowed via the ring homomorphism $\lambda f$. Then,
it follows from $\lambda$ being a homological ring epimorphism that
the multiplication map $S\otimes_RW\ra W$  is an isomorphism of
$S$-modules and that $\Tor^R_i(S, W)=0$ for all $i>0$ (see
\cite[Theorem 4.4]{GL}). Therefore, for any $i>0$, we have
$\Tor^{T}_i(M,W)\simeq \Tor^{R}_i(S, W)=0$.  Note that
$$M\otimes_{T}W=(S\otimes_RT)\otimes_T W\simeq
S\otimes_R(T\otimes_TW)\simeq S\otimes_RW\simeq W$$ as $S$-modules.
Thus the map $\delta_{W}$ is an isomorphism of $S$-modules. So the
above-mentioned general result follows.

Now, by applying the above general result to the ring $\Gamma$, we
can show that $\delta_{\Gamma}$ is an isomorphism and
$\Tor^T_n(M,\Gamma)=0$ for any $n>0$. This completes the proof of
Lemma \ref{gen} (2). $\square$

\medskip
{\bf Proof of Theorem \ref{th1} (1)}

Let $(\lambda, \mu, M, m)$ be a given exact context,
where $\lambda: R\ra S$ and $\mu:R\ra T$ are ring homomorphisms. By Lemma \ref{good},
the map $\theta$ is homological if and only if
$H^n\big({i_*i^*}(B)\big)=0$ for all $n\neq 0$. However, by Lemma \ref{calculation},
we see that $H^n\big({i_*i^*}(B)\big)\simeq H^n\big({i_*i^*}(Be_1)\big)\oplus H^n\big({i_*i^*}(Be_1)\big)\simeq \bigoplus_{i=1}^4\Tor^R_{-n}(T,S)$ for each $n\in\mathbb{Z}$. Thus $(a)$
and $(b)$ are equivalent. This shows the first part of Theorem \ref{th1} (1).

Assume that $(\lambda,\mu)$ is an exact pair such that $\lambda$ is
homological. Let $\Lambda:=T\boxtimes_RS$ be the noncommutative
tensor product of $(\lambda, \mu, M, m)$ (see Lemma \ref{ring}), and
$C:=M_2(\Lambda)$. Note that we have the following
commutative diagram of ring homomorphisms:
$$\xymatrix{
  R \ar[d]_-{\mu} \ar[r]^-{\lambda}
                & S \ar[d]_{\rho}  \\
  T \ar[r]^-{\phi}
                & \,\; \Lambda,
 }$$
and that the map $\theta_{\rho,\phi}$ defined in Lemma \ref{gen} (b) is equal to $\theta:B\to C$ by Remark \ref{coproduct} (1). It follows from Lemma \ref{gen} (2) that the statements $(a)$ and $(c)$ in Theorem \ref{th1} are equivalent. This finishes the proof of Theorem \ref{th1} (1). $\square$

%

\medskip
Combining Theorem \ref{th1} (1) with  Lemmas \ref{good} and \ref{sum} (3), we have the following result.

\begin{Koro} \label{key lemma}
If one of the assertions in Theorem \ref{th1} (1) holds, then there is a
recollement of derived module categories:
$$
\xymatrix@C=1.2cm{\D{C}\ar[r]^-{D(\theta_*)}
&\D{B}\ar[r]^-{\;j^!}\ar@/^1.2pc/[l]\ar_-{C\otimesL_B-}@/_1.2pc/[l]
&\D{R} \ar@/^1.2pc/[l]\ar@/_1.2pc/[l]_{j_!\;}}$$

\medskip
\noindent where $D(\theta_*)$ is the restriction functor induced by
$\theta:B\to C$, and where
$$j_!={_B}\cpx{P}\otimesL_R-\;\; \mbox{and}\;\; j^!=\cpx{\Hom}_B(\cpx{P},-)\simeq {_R}\cpx{P}{^*}\cpx{\otimes}_B-.$$
\end{Koro}

\smallskip
To prove Theorem \ref{th1} (2), we first establish the following
result which describes relationships among projective dimensions of
special modules over different rings. For an $R$-module $X$, we
denote the projective dimension of $X$ by $\pd(_RX)$.

\begin{Koro}\label{proj.dimension}
Assume that one of the assertions in Theorem \ref{th1} (1) holds. Then we have the following:

$(1)$ $\pd(_RS)\leq \max{\{1, \pd(_BC)\}}$ and $\pd(_BC)\leq \max{\{2, \pd(_RS)+1\}}$. In particular, $\pd(_RS)<\infty$ if and only if $\pd(_BC)<\infty$.

\smallskip
$(2)$ $\pd(T_R)\leq\max{\{1, \pd(C_B)\}}$ and $\pd(C_B)\leq \max{\{2, \pd(T_R)+1\}}$. In particular, $\pd(T_R)<\infty$ if and only if $\pd(C_B)<\infty$.
\end{Koro}

{\it Proof.} Note that the ring homomorphisms $\mu\opp: R\opp\to
T\opp$ and $\lambda\opp: R\opp\to S\opp$, together with $(M,m)$ form
an exact context. So it is sufficient to show $(1)$ because $(2)$
can be shown similarly.

We first show that $\pd(_RS)\leq \max{\{1, \pd(_BC)\}}$.

To see this inequality, we use the recollement given in Corollary
\ref{key lemma}. Clearly, there is a triangle in $\D{B}$:
$$
\cpx{P}\otimesL_R\cpx{P}{^*}\lra B\lraf{\theta} C\lra
\cpx{P}\otimesL_R\cpx{P}{^*}[1].
$$
This implies that $\Cone(\theta)$ is isomorphic in $\D{B}$ to the complex
$\cpx{P}\otimesL_R\cpx{P}{^*}[1]$. Since $\Cone(\lambda)\simeq\Cone(m\cdot)$ in $\D{R}$ by the diagram $(\sharp)$ in Section \ref{sect3.0}, we see that $\cpx{P}{^*}[1]\simeq S\oplus\Cone(m\cdot) \simeq S\oplus
\Cone(\lambda)$ in $\D{R}$, and therefore
$$\Cone(\theta)\simeq\cpx{P}\otimesL_RS\oplus P\otimesL_R\Cone(\lambda)\;\;\mbox{in}\;\;\D{B}.$$
As $\cpx{P}\otimesL_R-:\D{R}\to \D{B}$ is fully
faithful, we have $$\Hom_{\D R}(S, Y[n])\simeq \Hom_{\D
B}(\cpx{P}\otimesL_RS,\cpx{P}\otimesL_RY[n])$$ for every $Y\in
R\Modcat$ and $n\in\mathbb{N}$.

Suppose that $\pd(_BC)<\infty$, and let $s:=\max{\{1, \pd(_BC)\}}.$
We claim that $\Hom_{\D R}(S, Y[n])=0$ for any $n>s$, and therefore
$\pd(_RS)\leq s$. Since $\cpx{P}\otimesL_RS$ is a direct summand of
$\Cone(\theta)$ in $\D{B}$, it is enough to show that
$\Hom_{\D{B}}(\Cone(\theta),\cpx{P}\otimesL_RY[n])=0$ for any $n>s$.

Recall that $\Cone(\theta)$ is the complex $0\to B\lraf{\theta} C\to
0$ with $B$ and $C$ in degrees $-1$ and $0$, respectively. Then
$\Cone(\theta)$ is isomorphic in $\D{B}$ to a bounded complex
$$
\cpx{X}: 0\lra X^{-s} \lra X^{1-s}\lra \cdots \lra X^{-1} \lra
X^0\lra 0
$$
such that $X^i$ are projective $B$-modules for all $0\leq i\leq s$.
Let $_pY$ be a deleted projective resolution of $Y$ in $R\Modcat$.
Then $\Hom_{\D B}(\Cone(\theta), \cpx{P}\otimesL_RY[n]) \simeq
\Hom_{\D B}(\cpx{X}, \cpx{P}\otimesL_RY[n])=\Hom_{\D B}(\cpx{X},
\cpx{P}\cpx{\otimes}_R(_pY)[n])\simeq\Hom_{\K B}(\cpx{X},
\cpx{P}\cpx{\otimes}_R({_pY})[n])=0$ for any $n>s$, where the last
equality is due to the observation that all positive terms of the
complex $\cpx{P}\cpx{\otimes}_R({_pY})$ are zero. This verifies the claim
and shows that $\pd(_RS)\leq s$.

Next, we show that $\pd(_BC)\leq \max{\{2, \pd(_RS)+1\}}$.

Suppose that $\pd(_RS)=m<\infty$, and let $$\cpx{M}: 0\lra
M^{-m} \lra M^{1-m}\lra \cdots \lra M^{-1} \lra M^0\lra 0$$ be a
deleted projective resolution of $_RS$, where $M^i$ are projective
$R$-modules for all $-m\leq i\leq 0$. Then
$\cpx{P}\otimesL_RS=\cpx{P}\cpx{\otimes}_R\cpx{M}$ in
$\D{B}$. Note that $\Cone(\lambda)$ is isomorphic in $\D{R}$ to a
complex of the form:
$$\widetilde{\cpx{M}}: 0\lra M^{-m}
\lra M^{1-m}\lra \cdots \lra M^{-1}\oplus R \lra M^0\lra 0.$$ In
particular, $\widetilde{M^i}=0$ for $i>0$ or $i<-\max{\{1, m\}}$.
Then $\cpx{P}\otimesL_R\Cone(\lambda)=\cpx{P}\cpx{\otimes}_R\widetilde{\cpx{M}}$ in $\D{B}$. Thus $$\Cone(\theta)\simeq
(\cpx{P}\cpx{\otimes}_R\cpx{M})\oplus(\cpx{P}\cpx{\otimes}_R\widetilde{\cpx{M}})\;\;\mbox{in}\;\;\D{B}.$$
Recall that $\cpx{P}$ is the two-term complex $0\ra Be_1\lraf{\varphi} Be_2\ra 0$
over $B$ with $Be_1$ and $Be_2$ in degrees $-1$ and $0$, respectively. This implies that $\Cone(\theta)$
is isomorphic in $\D{B}$ to a complex of the following form:
$$\cpx{N}: 0\lra N^{-t}
\lra N^{1-t}\lra \cdots \lra N^{-1} \lra N^0\lra 0,$$ where $t:=
\max{\{2, m+1\}}$ and $N^i$ are projective $B$-modules for all
$-t\leq i\leq 0$. Now, let $X\in B\Modcat$. Then
$$
\Hom_{\D{B}}(\Cone(\theta), X[n])\simeq
\Hom_{\D{B}}(\cpx{N}, X[n])\simeq \Hom_{\K{B}}(\cpx{N},
X[n])=0
$$
for any $n>t$. Applying $\Hom_{\D{B}}(-, X[n])$ to the canonical distinguished triangle
$$
\Cone(\theta)[-1]\lra B\lraf{\theta} C\lra \Cone(\theta)
$$
in $\D{B}$, we have $\Ext^n_{B}(C, X)\simeq \Hom_{\D{B}}(C, X[n])=0$ for any
$n>t$. This shows $\pd(_BC)\leq t$. $\square$

\medskip
{\bf Proof of Theorem \ref{th1} (2)}

Let $\Lambda:=T\boxtimes_RS$ and $C:=M_2(\Lambda)$. Then the $\Lambda$-$C$-bimodule  $(\Lambda, \Lambda)$ induces an equivalence of module categories: $$(\Lambda,\Lambda)\otimes_C-:C\Modcat\lra \Lambda\Modcat. $$
In view of derived module categories, we obtain a triangle equivalence $$(\Lambda,\Lambda)\otimes_C-:\D{C}\lraf{\simeq} \D{\Lambda}.$$

Now, assume that one of the assertions in Theorem \ref{th1} (1) holds. By the above equivalence, we know from Corollary \ref{key lemma} that there exists a recollement of derived module categories:
$$
\xymatrix@C=1.2cm{\D{\Lambda}\ar[r]
&\D{B}\ar[r]\ar@/^1.2pc/[l]\ar_-{G}@/_1.2pc/[l]
&\D{R} \ar@/^1.2pc/[l]\ar@/_1.2pc/[l]_{j_!\;}}\vspace{0.4cm}$$
where $G:=(\Lambda,\Lambda)\otimesL_B-$ and $j_!={_B}\cpx{P}\otimesL_R-$.
This shows the first part of Theorem \ref{th1} (2).

By \cite[Theorem 3]{NS2}, the recollement in Corollary \ref{key lemma} can be
restricted to a recollement at $\mathscr{D}^-$-level:
$$\xymatrix{
\Df{C}\quad\ar^-{D(\theta_*)}[r] &\Df{B}\ar[r]^-(.6){j^!}
\ar@/^1.2pc/[l]\ar_-{C\otimesL_B-}@/_1.4pc/[l] &
\Df{R}\ar_-{j_!}@/_1.2pc/[l]\ar@/^1.2pc/[l]}\vspace{0.2cm}$$ if and only if the image of the object $C\in\D{C}$ under the functor $D(\theta_*)$ is isomorphic to a bounded complex of projective $B$-modules, that is $\pd(_BC)<\infty$. Furthermore, this $\mathscr{D}^-$-level
recollement can be restricted to $\mathscr{D}^b$-level
$$\xymatrix{
\Db{C}\quad\ar^-{D(\theta_*)}[r] &\Db{B}\ar[r]^-(.6){j^!}
\ar@/^1.2pc/[l]\ar_-{C\otimesL_B-}@/_1.4pc/[l] &
\Db{R}\ar_-{j_!}@/_1.2pc/[l]\ar@/^1.2pc/[l]}$$

\medskip \noindent
provided that $\pd(C_B)<\infty$. However, by Corollary \ref{proj.dimension}, we see that $\pd(_BC)<\infty$ if and only if $\pd(_RS)<\infty$, and that $\pd(C_B)<\infty$ if and only if $\pd(T_R)<\infty$. Identifying $\Db{C}$ with $\Db{\Lambda}$ up to equivalence, we finish the proof of the second part of Theorem \ref{th1} (2). $\square$

\subsection{Proofs of Corollaries\label{sect4}}
In this section, we shall prove all corollaries of Theorem
\ref{th1}, which were mentioned in the introduction.

All notation introduced in the previous sections will be kept.
As in Section
1, we fix a ring homomorphism $\lambda: R\to S$, and let
$$(\ast\ast)\quad R\lraf{\lambda} S\lraf{\pi}\cpx{Q}\lraf{\nu} R[1]$$
be the distinguished triangle in the homotopy category $\K{R}$ of
$R$, where the complex $\cpx{Q}$ stands for the mapping cone of
$\lambda$. Now, we set $S':=\End_{\D {R}}(\cpx{Q})$ and define $\lambda':R\to
S'$ by $r\mapsto \cpx{f}$ for $r\in R$, where $\cpx{f}$ is the chain
map with $f^{-1}:=\cdot r$, \,$f^0:=\cdot(r)\lambda$ and $f^{\,i}=0$
for $i\neq 0, -1.$ Here, $\cdot\,r$ and $\cdot\,(r)\lambda$ stand
for the right multiplication maps by $r$ and $(r)\lambda$,
respectively. These data can be recorded in the following commutative diagram:
$$ \xymatrix{
R\ar[d]^-{\cdot\,r}\ar[r]^-{\lambda}&S\ar[d]^-{\cdot\,(r)\lambda}\ar[r]^-{\pi}&\cpx{Q}
\ar[d]^-{\cpx{f}}\ar[r]^-{\nu}&R[1]\ar[d]^-{(\cdot\,r)[1]}\\
R\ar[r]^-{\lambda}&S\ar[r]^-{\pi}&\cpx{Q}\ar[r]^-{\nu}&R[1]}
$$
The map $\lambda'$ is called the ring homomorphism \emph{associated
to} $\lambda$. If $\lambda$ is injective, then we shall identify $\cpx{Q}$
with $S/R$ in $\D{R}$, and further, identify $\lambda'$ with the induced map
$R\to\End_R(S/R)$ by the right multiplication map.

Recall that $\Lambda$ denotes the ring $\End_{\D R}\big(S\oplus
\cpx{Q}\big)$ and that $\pi^*$ is the induced map
$$\Hom_{\D{R}}(S\oplus \cpx{Q},\, \pi):\Hom_{\D{R}}(S\oplus
\cpx{Q},\,S)\lra\Hom_{\D{R}}(S\oplus \cpx{Q},\,\cpx{Q}).$$ Let
$\lambda_{\pi^*}: \Lambda\to \Lambda_{\pi^*}$ stand for the
noncommutative localization  of $\Lambda$ at $\pi^*$.

Note that $\Hom_{\D R}(S,\cpx{Q})$ is an $S$-$S'$-bimodule containing $\pi$. Now we define a homomorphism of $S$-$S'$-bimodules:
$$\gamma:\;\;S\otimes_RS'\lra \Hom_{\D R}(S,\cpx{Q}),\;\; s\otimes f\mapsto
(\cdot s)(\pi\,f)$$
for $s\in S$ and $f\in S'$. This induces the following ring homomorphism:
$$\tau:=\left(\begin{array}{cc} \sigma & \gamma\\
0 & 1\end{array}\right):\;
\left(\begin{array}{lc}S& S\otimes_RS'\\
0& S'\end{array}\right)\lra\left(\begin{array}{lc} \;\;\quad \End_R(S)& \Hom_{\D R}(S,\cpx{Q})\\
\Hom_{\D R}(\cpx{Q},S)& S'\end{array}\right)=\Lambda$$
where $\sigma:S\to \End_R(S)$ is the inclusion under the identification of $S$ with $\End_S(S)$.

It is natural to ask when the quadruple $\big(\lambda,
\lambda',\Hom_{\D R}(S,\cpx{Q}), \pi\big)$ is an exact context.
Actually, in \cite[Lemma 6.5 (3)]{xc1}, we proved that if $\lambda$
is an injective ring epimorphism with $\Tor_1^R(S,S)=0$, then the
pair $(\lambda,\lambda')$ is exact. As a generalization of this
result, we show the following statement in which $\lambda$ is not
necessarily injective and $\cpx{Q}$ may have non-zero
cohomology in two degrees.

\begin{Lem}\label{general case}
Assume that $\Hom_R\big(S, \Ker(\lambda)\big)=0$. Then the quadruple $\big(\lambda, \lambda',\Hom_{\D R}(S,\cpx{Q}), \pi\big)$ is an exact context. If $\lambda$ is additionally a ring epimorphism, then $(\lambda,\lambda')$ is an exact pair. In this case, both $\gamma$ and $\tau$ are isomorphisms.
\end{Lem}
{\it Proof.} Applying $\Hom_{\D{R}}(-,\cpx{Q})$ to the triangle $(**)$, we have the following long exact sequence: {\footnotesize $$\Hom_{\D R}(S[1],\cpx{Q})\lraf{(\lambda[1])_*}\Hom_{\D R}(R[1],\cpx{Q})\lraf{\nu_*}\Hom_{\D{R}}(\cpx{Q},\cpx{Q})\lraf{\pi_*}
\Hom_{\D{R}}(S,\cpx{Q})\lraf{\lambda_*}\Hom_{\D R}(R,\cpx{Q}).$$}Since $\Hom_R\big(S, \Ker(\lambda)\big)=0$, the map $\Hom_R(S,\lambda)$ is injective. As $\Hom_{\D R}(S[1], S)\simeq\Hom_{\D R}(S, S[-1])\simeq \Ext^{-1}_R(S,S)=0$, we obtain $\Hom_{\D R}(S[1],\cpx{Q})=0$ by applying $\Hom_{\D R}(S[1],-)$ to the triangle ($**$).
Thus the above map $\nu_*$ is injective.

Next, we show that $\lambda_*=\Hom_{\D R}(\lambda, \cpx{Q}): \Hom_{\D{R}}(S, \cpx{Q})\lra \Hom_{\D{R}}(R,\cpx{Q})$ is surjective. In fact, the following diagram: $$
\begin{CD}\Hom_{\D{R}}(S,\cpx{Q})@>{\lambda_*}>>\Hom_{\D{R}}(R,\cpx{Q})\\
@AAA @AA{\simeq}A\\
\Hom_{\K{R}}(S,\cpx{Q})@>{\Hom_{\K R}(\lambda, \cpx{Q})}>>\Hom_{\K{R}}(R,\cpx{Q})
\end{CD}$$
is commutative, where the vertical maps are the canonical localization maps from $\K{R}$ to $\D{R}$. Since $\Hom_R(\lambda, S):\Hom_R(S,S)\to\Hom_R(R,S)$ is surjective, we see that $\Hom_{\K R}(\lambda, \cpx{Q})$ is surjective. This implies that $\lambda_*$ is surjective, as desired.

Let $\cdot \pi: S\to \Hom_{\D R}(S, \cpx{Q})$ stand for the right multiplication by $\pi$ map. Then we have the following exact commutative diagram
$$\xymatrix{
0\ar[r]&\Ker(\lambda)\ar[r]\ar[d]_{\simeq}&R
\ar[d]_-{\lambda'}\ar[r]^-{\lambda}
& S\ar[d]_{\cdot\pi}\ar[r]&\Coker(\lambda)\ar[d]_{\simeq}\ar[r]&0\\
0\ar[r] &\Hom_{\D R}(R[1],\cpx{Q})\ar[r]^-{\nu_*}
&\Hom_{\D{R}}(\cpx{Q},\cpx{Q})\ar[r]^-{\pi_*} &
\Hom_{\D{R}}(S,\cpx{Q})\ar[r]^-{\lambda_*}&\Hom_{\D R}(R,\cpx{Q})\ar[r] & 0}$$
where the isomorphisms follow from the fact that $\Hom_{\D R}(R[n], \cpx{Q})\simeq H^{-n}(\cpx{Q})$ for $n\in\mathbb{Z}$. Consequently, the chain map
$$\big(\lambda', \,\cdot\pi\big):\Cone(\lambda)\lra \Cone(\pi_*)$$ is a quasi-isomorphism. This implies that all of their cohomologies are isomorphic.
Note that $\pi_*$ is exactly the left multiplication by $\pi$ map. Thus the square in the middle of the above diagram is a pull-back and push-out diagram, and therefore the quadruple $\big(\lambda, \lambda',\Hom_{\D R}(S,\cpx{Q}), \pi\big)$ is an exact context.

Assume that $\lambda$ is a ring epimorphism. By Corollary
\ref{app1}, the map $\gamma:S\otimes_RS'\to \Hom_{\D R}(S,\cpx{Q})$
is an isomorphism of $S$-$S'$-bimodules. This implies that
$(\lambda,\lambda')$ is an exact pair by Lemma \ref{exact pair}. It
remains to show that the ring homomorphism $\tau$ is an isomorphism.
In fact, since $\lambda$ is a ring epimorphism, we have $S\simeq
\End_S(S)=\End_R(S)$ as rings. Thus $\sigma: S\ra \End_R(S)$ is an
isomorphism. Note that the composite of $\sigma$ with the valuation
map
$$\Hom_R(\lambda, S):\End_R(S)\lra \Hom_R(R,S)=S: \;f\mapsto (1)f
\quad \mbox{for}\;  f\in \End_R(S)$$ coincides with the identity map
of $S$. This implies that $\Hom_R(\lambda, S)$ is an isomorphism,
and therefore $\Hom_{\D R}(\cpx{Q}, S)=0$. Consequently, the map
$\tau$ is an isomorphism.  $\square$

\bigskip
{\bf Proof of Corollary \ref{new}}

By Lemma \ref{general case}, there is an isomorphism of rings:
$$\tau^{-1}:\quad \Lambda:=\End_{\D R}\big(S\oplus \cpx{Q}\big)
\lraf{\simeq} B:=\left(\begin{array}{lc} S & S\otimes_RS'\\
0 & S'\end{array}\right)$$ which sends
$\left(\begin{array}{cc} 0 & \pi\\
0 & 0\end{array}\right)$
to $\left(\begin{array}{cc} 0 & 1\otimes 1\\
0 & 0\end{array}\right)$. Set $e_1:= \left(\begin{array}{ll} 1 &0\\
0 & 0\end{array}\right)$ and $e_2:=\left(\begin{array}{ll} 0 &0\\
0 & 1\end{array}\right)\in B$. Let
$\varphi: Be_1\to Be_2$ be the map sending $\left(\begin{array}{l} s \\
0 \end{array}\right)$ to $\left(\begin{array}{c} s\otimes 1\\
0 \end{array}\right)$ for $s\in S$. Then $\pi^*$ corresponds to
$\varphi$ under the isomorphism $\tau^{-1}$. Thus $\Lambda_{\pi_*}=B_\varphi$ and $\lambda_{\pi^*}=\tau^{-1}\lambda_\varphi$.
Moreover, by Corollary \ref{general case}, the pair $(\lambda,\lambda')$ is exact. It follows from Lemma \ref{universal} that $\lambda_\varphi$ coincides with $\theta:B\to C:=M_2(S\boxtimes_R S')$. This means that $\lambda_{\pi^*}$ is homological if and only if $\theta$ is homological.  Note that $S\boxtimes_R S'\simeq S\sqcup_RS'\simeq \End_{S'}(S'\otimes_RS)$ as rings by Remark \ref{coproduct} (2). Since $\lambda$ is homological, Corollary \ref{new} follows immediately from Theorem \ref{th1}. $\square$

\begin{Rem}\label{stronger}
The equivalences of $(1)$ and $(3)$ in Corollary \ref{new} can be
obtained under a weaker assumption, instead of the 'homological'
assumption on $\lambda$. Precisely, we have the following result:

If $\Hom_R\big(S,\Ker(\lambda)\big)=0$ and $\End_S(S)=\End_R(S)$, then the map $\lambda_{\pi^*}:\Lambda\to \Lambda_{\pi^*}$ is homological if and only if $\Tor^R_i(S',S)=0$ for each $i\geq 1$.

\smallskip
{\it Proof.} If $\Hom_R\big(S,\Ker(\lambda)\big)=0$, then the
quadruple $\big(\lambda, \lambda',\Hom_{\D R}(S,\cpx{Q}), \pi\big)$
is an exact context by Lemma \ref{general case}. Moreover, if
$\End_S(S)=\End_R(S)$, then the map $\Hom_R(\lambda,
S):\Hom_R(S,S)\to\Hom_R(R,S)$
is an isomorphism, which leads to $$\Hom_{\D R}(\cpx{Q}, S)=0\;\;\mbox{and}\;\;\Lambda\lraf{\simeq}\left(\begin{array}{lc} S& \Hom_{\D R}(S,\cpx{Q})\\
0 & S'\end{array}\right).$$ Now, the above-mentioned result follows
immediately from Lemma \ref{universal} and Theorem \ref{th1} (1).
$\square$
\end{Rem}

\medskip
Combining Corollary \ref{new} with Lemma \ref{epi01}, we get the
following criterion for $\lambda_{\pi^*}$ to be homological.

\begin{Koro}\label{new1}
Let $\Sigma$ be a set of homomorphisms between finitely generated
projective $R$-modules. Suppose that $\lambda_{\Sigma}:R\to
R_{\Sigma}$ is homological such that
$\Hom_R\big(R_\Sigma,\Ker(\lambda_{\Sigma})\,\big)=0$. Set
$S:=R_{\Sigma}$, $\lambda:=\lambda_{\Sigma}$ and
$\Phi:=\{S'\otimes_{R}f\mid f\in\Sigma\}$ . Then the noncommutative
localization $\lambda_{\pi^*}:\Lambda\to \Lambda_{\pi^*}$ of
$\Lambda$ at $\pi^*$ is homological if and only if the noncommutative
localization $\lambda_{\Phi}:S'\to S'_{\Phi}$ of $S'$ at $\Phi$ is
homological. In particular, if one of the above equivalent
conditions holds , then there is a recollement of derived module
categories:

$$\xymatrix@C=1.2cm{\D{S'_{\Phi}}\ar[r]&\D{\Lambda}\ar[r]
\ar@/^1.2pc/[l]\ar@/_1.2pc/[l]
&\D{R}\ar@/^1.2pc/[l]\ar@/_1.2pc/[l]}.$$
\end{Koro}

\bigskip
As a consequence of Corollary \ref{new1}, we obtain the following
result which can be used to adjudge whether a noncommutative
localizations of the form $\lambda_{\pi^*}: \Lambda\to
\Lambda_{\pi^*}$ is homological or not.

\begin{Koro}\label{new2}
Let $F\subseteq D$ be an arbitrary extension of rings. Let $\omega:
D\ra D/F$ be the canonical surjection of $F$-modules.
Set $R:=\left(\begin{array}{cc} D & D\\
0 & F\end{array}\right)$ and $S:=M_2(D)$. Let $\lambda:R\to S$ be
the canonical inclusion, and let $\pi: S\ra S/R$ be the canonical
surjection. Then the noncommutative localization $\lambda_{\pi^*}:\Lambda\to \Lambda_{\pi^*}$ of
$\Lambda$ at $\pi^*$ is homological if and only if the noncommutative
localization $\lambda_{\omega^*}:E\to E_{\omega^*}$ of $E$ at
${\omega^*}$ is homological, where $E:=\End_F(D\oplus D/F)$, and
$\omega^*: \Hom_F(D\oplus D/F, \,D)\to\Hom_F(D\oplus D/F,\, D/F)$ is the
homomorphism of $E$-modules induced by $\omega$.
\end{Koro}

{\it Proof.} Since $\cpx{Q}$ can be identified with $S/R$ in
$\D{R}$, we have $S'=\End_R(S/R)$. Thus the map $\lambda': R\ra S'$
is given by the right multiplication.
Set $e_1:= \left(\begin{array}{ll} 1 &0\\
0 & 0\end{array}\right)$, $e_2:=\left(\begin{array}{ll} 0 &0\\
0 & 1\end{array}\right)$ and $e_{12}:=\left(\begin{array}{ll} 0 &1\\
0 & 0\end{array}\right)\in R$. Furthermore, let $\varphi: Re_1\ra
Re_2$ and $\varphi':S'(e_1)\lambda'\to S'(e_2)\lambda'$ be the right
multiplication maps of $e_{12}$ and $(e_{12})\lambda'$,
respectively.

It follows from Lemma \ref{universal} (see also \cite[Theorem 4.10]{Sch})
and $D\sqcup_{F}F=D$ that $\lambda:R\to S$ is  the noncommutative localization of $R$ at $\varphi$.
In particular, $\lambda$ is a ring epimorphism. Since $S\simeq
e_1R\oplus e_1R$ as right $R$-modules, the embedding $\lambda$ is
even homological. Note that $S'\otimes_R\varphi$ can be identified
with $\varphi'$. By Corollary \ref{new1}, the map $\lambda_{\pi^*}:
\Lambda\to \Lambda_{\pi^*}$ is homological if and only if the map
$\lambda_{\varphi'}:S'\to S'_{\varphi'}$ is homological.

Clearly, $R/Re_1R\simeq F$ as rings. So, every $F$-module can be
regarded as an $R$-module. In particular, the $F$-module $D\oplus D/F$ can be considered
as an $R$-module. Further, one can check that the map
$$\alpha:\;  D\oplus D/F\to S/R, \quad
(d,t+F)\mapsto \left(\begin{array}{lc} 0 & 0\\
d & t\end{array}\right)+R$$ for $d,t\in D$, is an isomorphism of
$R$-modules. Thus $S'\simeq E$ as rings. Under this isomorphism,
$\varphi'$ corresponds to $\omega^*$, and therefore
$S'_{\varphi'}\simeq E_{\omega^*}$ as rings. It follows that $\lambda_{\varphi'}:S'\to
S'_{\varphi'}$ is homological if and only if so is $\lambda_{\omega^*}:E\to E_{\omega^*}$. This finishes the proof. $\square$

\medskip
Before starting with the proof of Corollary \ref{tilting}, we
introduce a couple of more definitions and notation.

Recall from \cite{Rk} that a complex $\cpx{U}$ in $\D{R}$ is called a
\emph{tilting} complex if $\cpx{U}$ is self-orthogonal, isomorphic
in $\D{R}$ to a bounded complex of finitely generated projective
$R$-modules, and ${\rm{Tria}}(\cpx{U})=\D{R}$. It is well known that
if $\cpx{U}$ is a tilting complex over $R$, then $\D{R}$ is
equivalent to $\D{\End_{\D{R}}(\cpx{U})}$ as triangulated categories
(see \cite[Theorem 6.4]{Rk}). In this case,  $R$ and
$\End_{\D{R}}(\cpx{U})$ are called \emph{derived equivalent}. We refer the reader to \cite{hx2} for some new advances in constructions of derived equivalences.

If $I$ is an index set, we denote by $\cpx{U}\,^{(I)}$ the direct
sum of $I$ copies of $\cpx{U}$ in $\D R$.

The following result generalizes some known results in the
literature. See, for example, \cite[Theorem 4.14]{GL}, \cite[Theorem
3.5 (5)]{HJ} and \cite[Lemma 3.1 (3)]{x}, where the ring homomorphism
$\lambda:R\to S$ is required to be injective. We shall use this
generalization to prove Corollary \ref{tilting}.

\begin{Lem}\label{der}
Let $\lambda:R\to S$ be a ring homomorphism, and let $I$ be an
arbitrary nonempty set. Define $\cpx{U}:=S\oplus \cpx{Q}$. Then
$\Hom_{\D R}(\cpx{U},\cpx{U}\,^{(I)}[n])=0$ for any $0\neq n\in
\mathbb{Z}$ if and only if the following conditions hold:

$(1)$ $\Hom_R\big(S,\Ker(\lambda))=0$ and

$(2)$ $\Ext^i_R(S,S^{(I)})=0=\Ext^{i+1}_R(S,R^{(I)})$ for any $i\geq
1$.

In particular, the complex $\cpx{U}$ is a tilting complex in $\D{R}$
if and only if $\Hom_R\big(S,\Ker(\lambda))=0$, $\Ext^1_R(S,S)=0$
and there is an exact sequence: $ 0 \ra P_1\ra P_0\ra {}_RS \ra 0$
of $R$-modules, such that $P_i$ is finitely generated and projective
for $i=0,1$.
\end{Lem}
{\it Proof.} Recall that we have a distinguished triangle
$(\ast\ast)\quad R\lraf{\lambda} S\lraf{\pi}\cpx{Q}\lraf{\nu}
R[1]$ in $\K{R}$.

First of all, we mention two general facts: Let $I$ be an arbitrary
nonempty set.

$(a)$ By applying $\Hom_{\D{R}}(-,S^{(I)})$ to $(\ast\ast)$, one can
prove that $$\Hom_{\D
R}(\cpx{Q},S^{(I)}[i])\simeq\Hom_{\D{R}}(S,S^{(I)}[i])\;\;\mbox{for}\;\;i\in\mathbb{Z}\setminus\{0\}\;\;\mbox{and}\;\;
\Hom_{\D{R}}(\cpx{Q},S^{(I)})\simeq\Ker{\big(\Hom_R(\lambda,S^{(I)})\big)}.$$

$(b)$ By applying $\Hom_{\D{R}}(-,R^{(I)})$ to $(\ast\ast)$, one can
show that $$\Hom_{\D
R}(\cpx{Q},R^{(I)}[j])\simeq\Hom_{\D{R}}(S,R^{(I)}[j])\;\;\mbox{for}\;\;j\in\mathbb{Z}\setminus\{0,1\}.$$

Next, we show the necessity of the first part of Lemma \ref{der}.

Suppose that $\Hom_{\D R}(\cpx{U},\cpx{U}\,^{(I)}[n])=0$ for any
$n\neq 0$. Then
$\Ext^i_R(S,S^{(I)})\simeq\Hom_{\D{R}}(\cpx{Q},S^{(I)}[i])=0$ for
any $i\geq 1$, and $\Hom_{\D{R}}(S,\cpx{Q}[-1])=0$. Consequently,
the map $\Hom_R(S,\lambda):\Hom_R(S,R)\ra\Hom_R(S,S)$ is injective.
This means that the condition $(1)$ holds. Further, applying
$\Hom_{\D{R}}(S,-)$ to the triangle $ R^{(I)}\lraf{\lambda^{(I)}}
S^{(I)}\lraf{\pi^{(I)}} \cpx{Q}{^{(I)}}\lra R^{(I)}[1]$, we get
$\Ext^{i+1}_R(S,R^{(I)})\simeq\Hom_{\D{R}}\big(S,\cpx{Q}{^{(I)}}[i]\big)=0$.
Thus, the conditions $(1)$ and $(2)$ in Lemma \ref{der} are
satisfied.

In the following, we shall show the sufficiency of the first part of
Lemma \ref{der}.

Assume that the conditions $(1)$ and $(2)$ in Lemma \ref{der} hold
true. Then, it follows from $(a)$ and $(b)$ that
$$\Hom_{\D{R}}(\cpx{Q},S^{(I)}[n])=0=
\Hom_{\D{R}}(\cpx{Q},R^{(I)}[m+1])$$  for
$n\in\mathbb{Z}\setminus\{0\}$ and
$m\in\mathbb{Z}\setminus\{-1,0\}$. Applying
$\Hom_{\D{R}}(\cpx{Q}\,-)$ to the triangle
$R^{(I)}\lraf{\lambda^{(I)}} S^{(I)}\lraf{\pi^{(I)}}
\cpx{Q}{^{(I)}}\lra R^{(I)}[1]$, one can show that
$\Hom_{\D{R}}(\cpx{Q},\cpx{Q}\,^{(I)}[m])=0 \;\;\mbox{for}\;\;
m\in\mathbb{Z}\setminus\{-1,0\}.$ Furthermore, we shall show that
the condition $(1)$ in Lemma \ref{der} implies also that
$\Hom_{\D{R}}(\cpx{Q},\cpx{Q}\,^{(I)}[-1])=0$: Clearly,
$\Hom_R\big(S,\Ker({\lambda})^{I}\,\big)$
$\simeq\Hom_R\big(S,\Ker({\lambda})\big)^{I}$ = $0$, where
$\Ker({\lambda})^{I}$ stands for the direct product of $I$ copies of
$\Ker({\lambda})$. Since $\Ker({\lambda})^{I}$ contains
$\Ker({\lambda})^{(I)}$ as a submodule, we infer that
$\Hom_R\big(S,\Ker({\lambda})^{(I)}\,\big)=0$ and
Ker$\big(\Hom_R(S,\lambda^{(I)})\big)$ $\simeq$
$\Hom_R\big(S, \Ker(\lambda)^{(I)}\big)$ = $0$. Now, it follows from the
following exact commutative diagram:
$$ \small{\xymatrix{
0\ar[r]&\Hom_{\D{R}}(\cpx{Q},\cpx{Q}\,^{(I)}[-1])
\ar@{_{(}-->}[d]\ar[r]^-{(\nu[-1])^*}&\Hom_{\D{R}}(\cpx{Q},R\,^{(I)})
\ar@{_{(}->}[d]^-{\pi_*}\ar[r]^-{{(\,\lambda^{(I)}\,)^*}}&\Hom_{\D{R}}(\cpx{Q},S\,^{(I)})\ar@{_{(}->}[d]^-{\pi_*}\\
0\ar[r]&\Ker\big(\Hom_R(S,\lambda^{(I)})\big)\ar[r]&\Hom_R(S,R^{(I)})\ar[r]^-{(\,\lambda^{(I)}\,)^*}&\Hom_R(S,S^{(I)})}}
$$
that
$\Ker\big(\Hom_R(S,\lambda^{(I)})\big)\simeq\Hom_R\big(S,\Ker({\lambda})^{(I)}\,\big)=0$,
and therefore $\Hom_{\D{R}}(\cpx{Q},\cpx{Q}\,^{(I)}[-1])=0$. Thus,
$$\Hom_{\D{R}}(\cpx{Q},\cpx{Q}\,^{(I)}[n])=0 \;\;\mbox{for}\;\;
n\neq 0.$$

It remains to prove that $\Hom_{\D{R}}(S,\cpx{Q}\,^{(I)}[n])=0
\;\;\mbox{for}\;\; n\neq 0.$ Actually, applying $\Hom_{\D{R}}(S,-)$
to the triangle $R^{(I)}\lraf{\lambda^{(I)}} S^{(I)}\lraf{\pi^{(I)}}
\cpx{Q}{^{(I)}}\lra R^{(I)}[1]$, we have the following long exact
sequence: {\small
$$\cdots\ra\Hom_{\D{R}}(S,S^{(I)}[j])\lra\Hom_{\D{R}}(S,\cpx{Q}{^{(I)}}[j])
\lra\Hom_{\D R}(S,R^{(I)}[j+1])\lraf{(\,\lambda^{(I)}\,)^*} \Hom_{\D
R}(S,S^{(I)}[j+1])\ra\cdots$$}for $j\in\mathbb{Z}$. Since
$\Hom_{\D{R}}(S,S^{(I)}[r])=0$ for $0\neq r\in \mathbb{Z}$ and
$\Hom_{\D R}(S,R^{(I)}[t])=0$ for  $t\in\mathbb{Z}\setminus\{0,1\}$,
we see that $\Hom_{\D{R}}(S,\cpx{Q}{^{(I)}}[j])=0$ for
$j\in\mathbb{Z}\setminus\{-1,0\}$ and that
$\Hom_{\D{R}}(S,\cpx{Q}{^{(I)}}[-1])\simeq\Ker\big(\Hom_R(S,\lambda^{(I)})\big)=0$.
It follows that $\Hom_{\D{R}}(S,\cpx{Q}\,^{(I)}[n])=0
\;\;\mbox{for}\;\; n\neq 0.$ Hence $\Hom_{\D
R}(\cpx{U},\cpx{U}\,^{(I)}[n])=0$ for any $n\neq 0$. This finishes
the proof of the sufficiency.

As to the second part of Lemma \ref{der}, we observe the following:
The complex $\cpx{U}$ over $R$ is a generator of $\D{R}$, that is,
${\rm{Tria}}(\cpx{U})=\D{R}$, since $R\in {\rm{Tria}}(\cpx{U})$ by
the triangle $(\ast\ast)$. Moreover, the complex $\cpx{U}$ is a
tilting complex in $\D{R}$ if and only if it is self-orthogonal, and
$_RS$ has a projective resolution of finite length consisting of
finitely generated projective $R$-modules. Furthermore, if $_RS$ has
finite projective dimension and $ \Ext^{i+1}_R(S,R^ {(I)})=0$ for
any $i\ge 1$, then $_RS$ does have projective dimension at most $1$.
Thus, by the first part of Lemma
\ref{der}, we can show the second part of Lemma \ref{der}. $\square$

\medskip
{\bf Proof of Corollary \ref{tilting}}

(1) By Corollary \ref{new}, the map $\lambda_{\pi^*}: \Lambda \ra
\Lambda_{\pi^*}$ is homological if and only if $\Tor_j^R(S',S)=0$
for all $j\ge 1$. Now, we assume that $_RS$ has projective dimension at most $1$.
Then $\Tor_i^R(S',S)=0$ for all $i\ge 2$. This implies that $\lambda_{\pi^*}$ is homological if and only if
$\Tor_1^R(S',S)=0$. Since $Be_2=S'\oplus S\otimes_RS'$ as right $R$-modules, it
suffices to show that $\Tor_1^R(Be_2,S)=0$.

In fact, from Lemma \ref{sum} (3), we obtain a triple
$(j_!,\,j^!,\,j_*)$ of adjoint functors. Let $\eta:
Id_{\D{B}}\to j_*j^! $ be the unit adjunction with respect to the
adjoint pair $(j^!,j_*)$. Then we have the following fact:

For any $\cpx{X}\in\D{B}$, there exists a canonical triangle in
$\D{B}$:
$$i_*i^!(\cpx{X})\lra\cpx{X}\lraf{\eta_{\cpx{X}}}j_*j^!(\cpx{X})\lra
i_*i^!(\cpx{X})[1],$$ where
$j_*j^!(\cpx{X})={\mathbb{R}}\Hom_R(\cpx{P}{^*},\cpx{\Hom}_B(\cpx{P},\cpx{X}))$.
For the other  triple $(i^*,\,i_*,\,i^!)$ of adjoint triangle
functors, we refer the reader to Lemma \ref{sum} (3).

Let
$$ 0 \lra P^{-1}\lraf{\delta} P^0\lra {}_RS \lra 0$$ be a projective
resolution of $_RS$ with all $P^j$ projective $R$-modules. This
exact sequence gives rise to a triangle $P^{-1}\ra P^0\ra S\ra
P^{-1}[1]$ in $\D{R}$. Then we see from the recollement
$(\star)$ in Lemma \ref{sum} (3) that there is the following
exact commutative diagram: {\small$$ \xymatrix{
i_*i^!(Be_2\otimes_RP^{-1})\ar[d]\ar[r]
&Be_2\otimes_RP^{-1}\ar[r]^-{\small{\eta_{Be_2\otimes_RP^{-1}}}}\ar[d]_-{1\otimes
\delta}&j_*j^!(Be_2\otimes_RP^{-1})\ar[r]\ar[d]_-{j_*j^!(1\otimes
\delta)} &i_*i^!(Be_2\otimes_RP^{-1})[1]\ar[d]\\
i_*i^!(Be_2\otimes_RP^{0})\ar[r]\ar[d]
&Be_2\otimes_RP^0\ar[d]\ar[r]^-{\small{\eta_{Be_2\otimes_RP^{0}}}}
&j_*j^!(Be_2\otimes_RP^{0})\ar[d]\ar[r] &i_*i^!(Be_2\otimes_RP^{0})[1]\ar[d]\\
i_*i^!(Be_2\otimesL_RS)\ar[r]\ar[d]&Be_2\otimesL_RS\ar[r]^-{\eta_{Be_2\otimesL_RS}}\ar[d]
&j_*j^!(Be_2\otimesL_RS)\ar[r]\ar[d] &i_*i^!(Be_2\otimesL_RS)[1]\ar[d] \\
i_*i^!(Be_2\otimes_RP^{-1})[1]\ar[r]&Be_2\otimes_RP^{-1}[1]\ar[r]
&j_*j^!(Be_2\otimes_RP^{-1})[1]\ar[r]&i_*i^!(Be_2\otimes_RP^{-1})[2]}
$$}Since $i_*i^*(Be_1)\simeq Be_2\otimesL_RS$ in $\D{B}$ by Lemma
\ref{calculation} (3), we know that $j_*j^!(Be_2\otimesL_RS)\simeq
j_*j^!i_*i^*(Be_1)=0$, due to  $j^!i_*=0$ in the recollement
$(\star)$. It follows that $j_*j^!(1\otimes \delta)$ is an
isomorphism, and so is $H^0(j_*j^!(1\otimes \delta))$.

Suppose that $H^0(\eta_{P}):P\to H^0\big(j_*j^!(P)\big)$ is
injective for any projective $B$-module $P$. Then
$H^0(\eta_{Be_2\otimes_RP^{-1}})$ is injective since $_RP^{-1}$ is
projective. It follows from the isomorphism $H^0(j_*j^!(1\otimes
\delta))$ that the map $1\otimes\delta:Be_2\otimes_RP^{-1}\to
Be_2\otimes_RP^0$ is injective. This implies that
$\Tor^R_1(Be_2,S)=0$, as desired.

Thus, in the following, we shall prove that  $H^0(\eta_{P}):P\to
H^0\big(j_*j^!(P)\big)$ is injective for any projective $B$-module
$P$.

First, we point out that $H^0(\eta_{P})$ is injective if and only if
$\Hom_{\D{B}}(B,P)\lraf{j^!}\Hom_{\D{R}}\big(j^!(B),j^!(P)\big)$ is
injective. To see this, we consider the the following composite of
maps:
$$\omega^n_{\cpx{X}}:\Hom_{\D{B}}(B,\cpx{X}[n])\lraf{j^!}\Hom_{\D{R}}\big(j^!(B),j^!(\cpx{X})[n]\big)
\lraf{\simeq}\Hom_{\D{B}}(B,j_*j^!(\cpx{X})[n])
$$ for each $n\in\mathbb{Z}$, where the second map is an isomorphism induced by the adjoint pair $(j^!,j_*)$.
Then, one can check directly that $\omega^n_{\cpx{X}}=\Hom_{\D
B}(B,\eta_{\cpx{X}[n]})$. It is known that the $n$-th cohomology
functor $H^{n}(-):\D{B}\to B\Modcat$ is naturally isomorphic to the
Hom-functor $\Hom_{\D{B}}(B,\,-[n])$. So, under this identification,
the map $\omega^n_{\cpx{X}}$ coincides with
$H^{n}(\eta_{\cpx{X}}):H^{n}(\cpx{X})\to H^{n}(j_*j^!(\cpx{X}))$. It
follows that  $H^0(\eta_{P})$ is injective if and only if so is the
map
$\Hom_{\D{B}}(B,P)\lraf{j^!}\Hom_{\D{R}}\big(j^!(B),j^!(P)\big)$.

Second, we claim that if $\Hom_{\D{B}}(i_*i^*(B),P)=0$, then
$\Hom_{\D{B}}(B,P)\lraf{j^!}\Hom_{\D{R}}\big(j^!(B),j^!(P)\big)$ is
injective.

Let $\varepsilon: j_!j^!\to Id_{\D{B}}$ be the counit adjunction
with respect to the adjoint pair $(j_!,j^!)$. Then, for each
$\cpx{X}\in \D{B}$, there exists a canonical triangle in $\D{B}$:
$$j_!j^!(\cpx{X})\lraf{\varepsilon_{\cpx{X}}}\cpx{X}\lra i_*i^*(\cpx{X})\lra
j_!j^!(\cpx{X})[1].$$ Now, we consider the following morphisms:
$$\Hom_{\D{B}}(B,\cpx{X}[m])\lraf{j^!}\Hom_{\D{R}}\big(j^!(B),j^!(\cpx{X})[m]\big)
\lraf{\simeq}\Hom_{\D{B}}\big(j_!j^!(B),\cpx{X}[m]\big)$$ for any
$m\in\mathbb{Z}$, where the last map is an isomorphism given by the
adjoint pair $(j_!,j^!)$. One can check that the composite of the
above two morphisms is the map
$\Hom_{\D{B}}(\varepsilon_{B},\cpx{X}[m]).$ This means that, to show
that
$\Hom_{\D{B}}(B,P)\lraf{j^!}\Hom_{\D{R}}\big(j^!(B),j^!(P)\big)$ is
injective, it suffices to show that
$\Hom_{\D{B}}(\varepsilon_{B},P)$ is injective. For this aim, we
apply $\Hom_{\D{B}}(-,P)$ to the triangle
$$j_!j^!(B)\lraf{\varepsilon_{B}}B\lra i_*i^*(B)\lra j_!j^!(B)[1]$$
and get the following exact sequence of abelian groups:
$$ \xymatrix{\Hom_{\D B}(i_*i^*(B),P)\ar[r]&\Hom_{\D
B}(B,P)\ar[rr]^-{\Hom_{\D{B}}(\varepsilon_{B},P)}&& \Hom_{\D
B}(j_!j^!(B),P)}.$$ If $\Hom_{\D{B}}(i_*i^*(B),P)=0$, then
$\Hom_{\D{B}}(\varepsilon_{B},P)$ is injective, and therefore the
map $j^!: \Hom_{\D{B}}(B,P)\to\Hom_{\D{R}}\big(j^!(B),j^!(P)\big)$
is injective, as desired.

Third, we  show that if $\Hom_R(S,S')=0$, then $\Hom_{\D{B}}(
i_*i^*(B),P)=0$ for any projective $B$-module $P$.

In fact, due to Lemma \ref{calculation} (1) and $(3)$, we have
$i_*i^*(Be_2)\simeq i_*i^*(Be_1)\simeq Be_2\otimesL_RS$ in $\D{B}$.
This implies that $\Hom_{\D{B}}( i_*i^*(B),P)=0$ if and only if
$\Hom_{\D{B}}(Be_2\otimesL_RS,P)=0$. Now, we consider the following
isomorphisms
$$\Hom_{\D{B}}(Be_2\otimesL_RS, P)\simeq\Hom_{\D{R}}(S,
{\mathbb{R}}\Hom_B(Be_2,P))\simeq\Hom_{\D{R}}(S,
e_2P)\simeq\Hom_R(S, e_2P).$$ Since $e_2B\simeq S'$ as $R$-modules,
we have $\Hom_R(S, e_2B)\simeq \Hom_R(S,S')=0$. Note that
$P\in\Add(_BB)$ and $e_2P\in\Add(_RS')$. Thus there is an index set
$I$ such that $e_2P$ is a direct summand of $(S')^{(I)}$. Since
$(S')^{(I)}$ is a submodule of the product $(S')^I$ of $S'$, it
follows that $\Hom_R(S,(S')^{(I)})$ is a subgroup of
$\Hom_R(S,(S')^I)$ which is isomorphic to $\Hom_R(S,S')^I$. Hence
$\Hom_R(S,(S')^{(I)})=0$, $\Hom_R(S,e_2P)=0$ and
$\Hom_{\D{B}}(i_*i^*(B),P)=0$, as desired.

Now, it remains to show that $\Hom_{R}(S,S')=0$. In the following,
we shall prove a stronger statement, namely,
$\Hom_{\D{R}}(S,S'[n])=0$ for any $n\in\mathbb{Z}$.

Since $\lambda$ is a ring epimorphism with $\Tor^R_1(S,S)=0$, we
know from \cite[Theorem 4.8]{Sch} that
$$\Ext^1_R(S,S^{(I)})\simeq\Ext^1_S(S,S^{(I)})=0$$ for any set $I$.
As $_RS$ is of projective dimension at most $1$, we can apply Lemma
\ref{der} to the complex $\cpx{U}:=S\oplus \cpx{Q}$, and get
$\Hom_{\D R}(\cpx{U},\cpx{U}[m])=0$ for $m\ne 0$. This implies that
$\Hom_{\D R}(\cpx{Q},\cpx{Q}[m])=0$ for $m\ne 0$, and that
$$H^m(\mathbb{R}\Hom_R(\cpx{Q},\cpx{Q}))\simeq
\Hom_{\D{R}}(\cpx{Q},\cpx{Q}[m])=\left\{\begin{array}{ll} 0 & \mbox{if}\; m\neq 0,\\
S' & \mbox{if }\; m= 0.\end{array} \right.$$ Thus the complex
$\mathbb{R}\Hom_R(\cpx{Q},\cpx{Q})$ is isomorphic in $\D{R}$ to the
stalk complex $S'$. On the one hand, by the adjoint pair
$\big(\cpx{Q}\otimesL_R-, \mathbb{R}\Hom_R(\cpx{Q},-)\big)$ of the
triangle functors, we have {\small
$$\Hom_{\D{R}}(S,S'[n])\simeq\Hom_{\D{R}}\big(S,\mathbb{R}\Hom_R(\cpx{Q},\cpx{Q})[n]\big)\simeq
\Hom_{\D{R}}\big(S,\mathbb{R}\Hom_R(\cpx{Q},\cpx{Q}[n])\big)\simeq
\mathbb{R}\Hom_R(\cpx{Q}\otimesL_RS,\cpx{Q}[n])$$}for any $n\in
\mathbb{Z}$. On the other hand, since $\lambda$ is homological by
assumption, the homomorphism $\lambda\otimesL_RS: R\otimesL_RS \ra
S\otimesL_RS$ is an isomorphism in $\D R$. It follows from the
triangle
$R\otimesL_RS \lraf{\lambda\otimesL_RS} S\otimesL_RS\lra \cpx{Q}\otimesL_RS\lra
R\otimesL_RS[1]$ that $ \cpx{Q}\otimesL_RS=0$. Hence
$\Hom_{\D{R}}(S,S'[n])\simeq\mathbb{R}\Hom_R(\cpx{Q}\otimesL_RS,\cpx{Q}[n])=0$
for any $n\in \mathbb{Z}.$

Thus, we have proved that, for any projective $B$-module $P$, the
homomorphism  $H^{0}(\eta_{P}):P\to H^{0}(j_*j^!(P))$ is injective
in $B\Modcat$. This finishes the proof of Corollary
\ref{tilting} (1).

(2) From Lemma \ref{sum} (3), we see
that the ring $\Lambda_{\pi^*}$ is zero if and only if the functor
$j^!$ induces a triangle equivalence from $\D{B}$ to $\D{R}$. This
is equivalent to the statement that $j^!(B)$ is a tilting complex
over $R$. Note that $j^!(B)\simeq \cpx{U}[-1]$.
Thus, the ring $\Lambda_{\pi^*}$ is zero if and only
if $\cpx{U}$ is a tilting complex over $R$. Now,
Corollary \ref{tilting} (2) follows directly from Lemma \ref{der}.
$\square$.

\smallskip
{\bf Proof of Corollary \ref{com1}}

$(1)$ Let $\lambda:R\to S$ be the inclusion, $\pi:S\to S/R$ the
canonical surjection and $\lambda':R\to S'$ the induced map by right multiplication.
Since $\lambda$ is injective, we know from Lemma \ref{general case}
that $\big(\lambda, \lambda', \Hom_R(S,S/R),\pi \big)$ is an exact context. If
$_RS$ is flat, then $\Tor^R_i(S', S)=0$ for all $i\geq 1$. Now, $(1)$ follows from Theorem \ref{th1}.

$(2)$ Since $\lambda$ is an injective ring epimorphism, the pair $(\lambda,\lambda')$ is exact by Lemma \ref{general case}. Since the ring $R$ is commutative and $\lambda$ is homological, the ring $S'$ is also commutative by \cite[Lemma 6.5 (5)]{xc1}. Consequently, the noncommutative tensor product $S'\boxtimes_RS$ coincides with the usual tensor product $S'\otimes_RS$ of $S'$ and $S$ over $R$ (see Section \ref{sect3.1}).

By Lemma \ref{gen}, we know that $\Tor^{R}_i(S,S')=0$ for any $i>0$.
Since $R$, $S$ and $S'$ are commutative rings, we have
$\Tor^{R}_i(S',S)\simeq\Tor^{R}_i(S,S')=0$. Thus the assertion $(3)$ in Corollary \ref{new}
is satisfied. It follows from Corollary \ref{new} that the ring homomorphism
$S'\otimes_R\lambda: S'\to S'\otimes_RS$ is a ring epimorphism.
This implies that $\End_{S'}(S'\otimes_RS)\simeq\End_{S'\otimes_RS}(S'\otimes_RS)\simeq S\otimes_RS'$ as rings. Note that $S\simeq \End_R(S)$ as rings and that $\Hom_R(S/R,S)=0$ since $\lambda$ is a ring epimorphism. Thus $\Lambda\simeq\End_R(S\oplus S/R)\simeq B$ as rings. Now, (2) is an immediate consequence of Corollary \ref{new}. $\square$

\medskip
{\bf Proof of Corollary \ref{com2}}

For a commutative ring $R$ and a multiplicative set $\Phi$ of $R$,
the localization map $\lambda: R\to S:=\Phi^{-1}R$ is always homological
since $_RS$ is flat. Therefore, by Corollary \ref{com1} (2), it
suffices to show that $\Hom_R(S/R,S)=0$ and that $S'\otimes_RS$ is
isomorphic to $\Psi^{-1}S'$. Actually, the former follows from the fact
that $\lambda$ is a ring epimorphism. To check the latter, we verify that the well
defined map $$\alpha: S'\otimes_R\Phi^{-1}R\lra \Psi^{-1}S',\;\; y \otimes\frac{r}{x}\mapsto
\frac{(r)\lambda'\,y}{(x)\lambda'}$$ for $y\in S'$, $r\in R$ and $x\in\Phi$,
is an isomorphism of rings, where $\lambda':R\to S'$ is the right multiplication map. Clearly, $\alpha$ is surjective. To see that this map is injective, we note that the map
$$\beta: \Psi^{-1}S'\lra S'\otimes_R\Phi^{-1}R,\;\;\frac{y}{(x)\lambda'}\mapsto y\otimes\frac{1}{x}$$ for $y\in S'$ and $x\in \Phi$, is a well defined ring homomorphism with
$\alpha\beta=1$. Observe that $\alpha$ preserves the multiplication of $S'\otimes_RS$. This finishes the proof of Corollary \ref{com2}. $\square$

\section{Examples \label{sect6}}

Now we present a few examples to show that some conditions in our
results cannot be dropped or weakened.

\medskip
(1) The condition that $\lambda: R\ra S$ is a homological ring
epimorphism in Corollary \ref{new} cannot be weakened to that
$\lambda: R\ra S$ is a ring epimorphism.

Let $R=\left(\begin{array}{ccc} k & 0          & 0\\
                                  k[x]/(x^2) & k          & 0 \\
                                  k[x]/(x^2) & k[x]/(x^2) & k
                         \end{array}\right),$\vspace{0.2cm}
where $k$ is a field and $k[x]$ is the polynomial algebra over $k$
in one variable $x$. Let $S$ be the $3$ by $3$ matrix ring
$M_3(k[x]/(x^2))$. Then the inclusion $\lambda$ of $R$ into $S$ is a
noncommutative localization of $R$, and therefore a ring epimorphism.
Further, we have $\Tor_1^R(S,S)=0\neq\Tor_2^R(S,S)$ (see
\cite{nrs}). Thus $\lambda$ is not homological. So, $_RS$ cannot
have projective dimension less than or equal to $1$. Moreover, one
can check that the ring homomorphism $\lambda':R\ra S'$ associated
to $\lambda$ is an isomorphism of rings. Recall that the pair $(\lambda,\lambda')$ is exact
by Lemma \ref{general case}, and further, that $S'\boxtimes_RS\simeq S\sqcup_{R}S'=S$ as rings by Remark \ref{coproduct}. In this sense, we have $\phi=(\lambda')^{-1}\,\lambda:S'\to S$.
Consequently, $\phi$ is not homological. However, due to Remark \ref{stronger}, the map $\lambda_{\pi^*}:\Lambda\to \Lambda_{\pi^*}$ is homological since $\Tor^R_i(S',S)\simeq \Tor^R_i(R,S)=0$ for each $i\geq 1$.
Hence, without the `homological' assumption on $\lambda$, the conditions $(1)$ and $(2)$ in Corollary \ref{new} are not equivalent.

(2) The condition that $\lambda$ is homological does not guarantee that the
noncommutative localization $\lambda_{\pi^*}:\Lambda\to \Lambda_{\pi^*}$
of $\Lambda$ at $\pi^*$ in Corollary \ref{new} is always
homological.

In the following, we shall use Corollary \ref{new2} to produce a
counterexample.

Now, take $F=\big\{\left(\begin{array}{lc} a & 0\\
b & a\end{array}\right)\mid a,b\in k\big\}$ and $D=\left(\begin{array}{lc} k & 0\\
k & k\end{array}\right)$ with $k$ a field.  Then one can verify that
the extension $\lambda:R\to S$, defined in Corollary \ref{new2}, is
homological, and that the canonical map $\omega: D\ra D/F$ is a
split epimorphism in $F$-Mod, and therefore $D\simeq F\oplus D/F$ as $F$-modules.
Let $e$ be the idempotent of $E$ corresponding the direct summand
$F$ of the $F$-module $D\oplus D/F$. Then $E_{\omega^*}\simeq
E/EeE\simeq M_2(k)$. Furthermore, the noncommutative localization
$\lambda_{\omega^*}:E\to E_{\omega^*}$ of $E$ at ${\omega^*}$ is
equivalent to the canonical surjection $\tau:E\to E/EeE$. Since
$\Ext^2_E(E/EeE,E/EeE)\neq 0$, the map $\tau$ is not
homological. This implies that $\lambda_{\omega^*}$ is not
homological, too. Thus $\lambda_{\pi^*}:\Lambda\to \Lambda_{\pi^*}$
is not homological by Corollary \ref{new2}, that is, the restriction
functor $D\big((\lambda_{\pi^*})_*\big):\D{\Lambda_{\pi^*}}\to
\D{\Lambda}$ is not fully faithful. In addition, one can check that,
for this extension, the $R$-module $_RS$ has infinite projective
dimension.

(3) In Corollary \ref{tilting} (1), we assume that the projective
dimension of $_RS$ is at most $1$. But there does exist an injective
homological ring epimorphism $\lambda: R\ra S$ such that the
projective dimension of $_RS$ is greater than $1$ and that
$\lambda_{\pi^*}: \Lambda\ra \Lambda_{\pi^*}$ is homological.

Let $R$ be a Pr\"ufer domain which is not a Matlis domain. Recall
that a Matlis domain is an integral domain $R$ for which the projective
dimension of the fractional field $Q$ of $R$ as an $R$-module is at
most $1$. In this case, the inclusion $\lambda: R\ra Q$ is an
injective homological ring epimorphism. By Corollary \ref{new}, the
map $\lambda_{\pi^*}: \Lambda\to \Lambda_{\pi^*}$ is homological.

\medskip
{\bf Acknowledgement.}
The authors thank Federik Marks for some discussions on noncommutative localizations. Parts of the paper were revised during a visit of the authors to the University of Stuttgart in May, 2013. The research work of the corresponding author CCX is partially
supported by the Natural Science Foundation (11331006, KZ201410028033) and the Chinese Ministry of Education.

{\footnotesize
\bigskip Hongxing Chen,

School of Mathematical Sciences, BCMIIS, Capital Normal University, 100048
Beijing, People's Republic of  China

{\tt Email: chx19830818@163.com}

\bigskip
Changchang Xi,

School of Mathematical Sciences, BCMIIS, Capital Normal University, 100048
Beijing, People's Republic of  China

{\tt Email: xicc@cnu.edu.cn}}


\end{document}